\documentclass[12equipped]{article}
\usepackage{}
\usepackage{amssymb}
\usepackage{amsfonts}

\usepackage{amssymb,amsmath,tikz}
\usepackage{graphicx}
\usepackage[plainpages]{hyperref}
\usepackage{enumerate}
\usepackage{verbatim}
\usepackage{mathrsfs}
\usepackage{latexsym}

\usepackage [autostyle, english = american]{csquotes}

\input epsf.tex

\usepackage[all]{xy}\usepackage{}

\usepackage{tikz}
\usetikzlibrary{automata,positioning}
\usepackage[top=1.5in, bottom=1.5in, left=1.5in, right=1.5in]{geometry}

\usetikzlibrary{shapes,shadows}

\usetikzlibrary{arrows,calc}
\tikzset{
>=stealth',
help lines/.style={dashed, thick},
axis/.style={<->},
important line/.style={thick},
connection/.style={thick, dotted},
}
\usetikzlibrary{patterns}
\newlength{\hatchspread}
\newlength{\hatchthickness}
\tikzset{hatchspread/.code={\setlength{\hatchspread}{#1}},
         hatchthickness/.code={\setlength{\hatchthickness}{#1}}}
\tikzset{hatchspread=3pt,
         hatchthickness=0.4pt}
\pgfdeclarepatternformonly[\hatchspread,\hatchthickness]% variables
   {custom north west lines}% name
   {\pgfqpoint{-2\hatchthickness}{-2\hatchthickness}}% lower left corner
   {\pgfqpoint{\dimexpr\hatchspread+2\hatchthickness}{\dimexpr\hatchspread+2\hatchthickness}}% upper right corner
   {\pgfqpoint{\hatchspread}{\hatchspread}}% tile size
   {% shape description
    \pgfsetlinewidth{\hatchthickness}
    \pgfpathmoveto{\pgfqpoint{0pt}{\hatchspread}}
    \pgfpathlineto{\pgfqpoint{\dimexpr\hatchspread+0.15pt}{-0.15pt}}
    \pgfusepath{stroke}
   }
\newcommand{\nc}{\newcommand}
\nc{\rnc}{\renewcommand}
\nc{\bb}[1]{{\mathbb #1}}
\nc{\bbA}{\bb{A}}\nc{\bbB}{\bb{B}}
\nc{\bbD}{\bb{D}}
\nc{\bbE}{\bb{E}}\nc{\bbF}{\bb{F}}\nc{\bbG}{\bb{G}}\nc{\bbH}{\bb{H}}
\nc{\bbI}{\bb{I}}\nc{\bbJ}{\bb{J}}\nc{\bbK}{\bb{K}}\nc{\bbL}{\bb{L}}
\nc{\bbM}{\bb{M}}\nc{\bbN}{{\bf N}}\nc{\bbO}{\bb{O}}\nc{\bbP}{\bb{P}}
\nc{\bbQ}{{\bf Q}}\nc{\bbR}{\bb{R}}\nc{\bbS}{\bb{S}}\nc{\bbT}{\bb{T}}
\nc{\bbU}{\bb{U}}\nc{\bbV}{\bb{V}}\nc{\bbW}{\bb{W}}\nc{\bbX}{\bb{X}}
\nc{\bbY}{\bb{Y}}\nc{\bbZ}{\bb{Z}}
\nc{\mbf}[1]{{\mathbf #1}}
\nc{\bfA}{\mbf{A}}\nc{\bfB}{\mbf{B}}\nc{\bfC}{\mbf{C}}\nc{\bfD}{\mbf{D}}
\nc{\bfE}{\mbf{E}}\nc{\bfF}{\mbf{F}}\nc{\bfG}{\mbf{G}}\nc{\bfH}{\mbf{H}}
\nc{\bfI}{\mbf{I}}\nc{\bfJ}{\mbf{J}}\nc{\bfK}{\mbf{K}}\nc{\bfL}{\mbf{L}}
\nc{\bfM}{\mbf{M}}\nc{\bfN}{\mbf{N}}\nc{\bfO}{\mbf{O}}\nc{\bfP}{\mbf{P}}
\nc{\bfQ}{\mbf{Q}}\nc{\bfR}{\mbf{R}}\nc{\bfS}{\mbf{S}}\nc{\bfT}{\mbf{T}}
\nc{\bfU}{\mbf{U}}\nc{\bfV}{\mbf{V}}\nc{\bfW}{\mbf{W}}\nc{\bfX}{\mbf{X}}
\nc{\bfY}{\mbf{Y}}\nc{\bfZ}{\mbf{Z}}
\nc{\bfa}{\mbf{a}}\nc{\bfb}{\mbf{b}}\nc{\bfc}{\mbf{c}}\nc{\bfd}{\mbf{d}}
\nc{\bfe}{\mbf{e}}\nc{\bff}{\mbf{f}}\nc{\bfg}{\mbf{g}}\nc{\bfh}{\mbf{h}}
\nc{\bfi}{\mbf{i}}\nc{\bfj}{\mbf{j}}\nc{\bfk}{\mbf{k}}\nc{\bfl}{\mbf{l}}
\nc{\bfm}{\mbf{m}}\nc{\bfn}{\mbf{n}}\nc{\bfo}{\mbf{o}}\nc{\bfp}{\mbf{p}}
\nc{\bfq}{\mbf{q}}\nc{\bfr}{\mbf{r}}\nc{\bfs}{\mbf{s}}\nc{\bft}{\mbf{t}}
\nc{\bfu}{\mbf{u}}\nc{\bfv}{\mbf{v}}\nc{\bfw}{\mbf{w}}\nc{\bfx}{\mbf{x}}
\nc{\bfy}{\mbf{y}}\nc{\bfz}{\mbf{z}}

\nc{\mcal}[1]{{\mathcal #1}}
\nc{\calA}{\mcal{A}}\nc{\calB}{\mcal{B}}\nc{\calC}{\mcal{C}}\nc{\calD}{\mcal{D}}
\nc{\calE}{\mcal{E}} \nc{\calF}{\mcal{F}}\nc{\calG}{\mcal{G}}\nc{\calH}{\mcal{H}}
\nc{\calI}{\mcal{I}}\nc{\calJ}{\mcal{J}}\nc{\calK}{\mcal{K}}\nc{\calL}{\mcal{L}}
\nc{\calM}{\mcal{M}}\nc{\calN}{\mcal{N}}\nc{\calO}{\mcal{O}}\nc{\calP}{\mcal{P}}
\nc{\calQ}{\mcal{Q}}\nc{\calR}{\mcal{R}}\nc{\calS}{\mcal{S}}\nc{\calT}{\mcal{T}}
\nc{\calU}{\mcal{U}}\nc{\calV}{\mcal{V}}\nc{\calW}{\mcal{W}}\nc{\calX}{\mcal{X}}
\nc{\calY}{\mcal{Y}}\nc{\calZ}{\mcal{Z}}
\nc{\fA}{\frak{A}}\nc{\fB}{\frak{B}}\nc{\fC}{\frak{C}} \nc{\fD}{\frak{D}}
\nc{\fE}{\frak{E}}\nc{\fF}{\frak{F}}\nc{\fG}{\frak{G}}\nc{\fH}{\frak{H}}
\nc{\fI}{\frak{I}}\nc{\fJ}{\frak{J}}\nc{\fK}{\frak{K}}\nc{\fL}{\frak{L}}\nc{\fM}{\frak{M}}
\nc{\fN}{\frak{N}}\nc{\fO}{\frak{O}}\nc{\fP}{\frak{P}}
\nc{\fQ}{\frak{Q}}\nc{\fR}{\frak{R}}\nc{\fS}{\frak{S}}\nc{\fT}{\frak{T}}
\nc{\fU}{\frak{U}}\nc{\fV}{\frak{V}}\nc{\fW}{\frak{W}}\nc{\fX}{\frak{X}}
\nc{\fY}{\frak{Y}}\nc{\fZ}{\frak{Z}}
\nc{\fa}{\frak{a}}\nc{\fb}{\frak{b}}\nc{\fc}{\frak{c}} \nc{\fd}{\frak{d}}
\nc{\fe}{\frak{e}}\nc{\fFf}{\frak{f}}\nc{\fg}{\frak{g}}\nc{\fh}{\frak{h}}
\nc{\fri}{\frak{i}}\nc{\fj}{\frak{j}}\nc{\fk}{\frak{k}}\nc{\fl}{\frak{l}}
\nc{\fn}{\frak{n}}\nc{\fo}{\frak{o}}\nc{\fp}{\frak{p}}
\nc{\fq}{\frak{q}}\nc{\fr}{\frak{r}}\nc{\fs}{\frak{s}}\nc{\ft}{\frak{t}}
\nc{\fu}{\frak{u}}\nc{\fv}{\frak{v}}\nc{\fw}{\frak{w}}\nc{\fx}{\frak{x}}
\nc{\fy}{\frak{y}}\nc{\fz}{\frak{z}}

\DeclareMathOperator{\rank}{rank}

\DeclareMathOperator{\codim}{codim} \DeclareMathOperator{\id}{id}
\DeclareMathOperator{\Image}{Im} \DeclareMathOperator{\Sym}{Sym}

 \DeclareMathOperator{\GL}{GL}

\DeclareMathOperator{\Hom}{{Hom}}

\DeclareMathOperator{\Hilb}{{Hilb}}

 \DeclareMathOperator{\Lie}{Lie}
 \DeclareMathOperator{\tr}{tr}

 \DeclareMathOperator{\End}{End}

      \DeclareMathOperator{\crit}{crit}

\DeclareMathOperator{\SH}{\textbf{SH}}

\DeclareMathOperator{\Crit}{Crit}

\DeclareMathOperator{\ad}{ad}

\DeclareMathOperator{\Sh}{Sh}

\DeclareMathOperator{\Dr}{Dr}

\newcommand{\surj}{\twoheadrightarrow}
\newcommand{\inj}{\hookrightarrow}

\def\angl#1{{\langle #1\rangle}}

\newcommand{\pt}{\text{pt}}

\newcommand{\bSh}{\textbf{Sh}}

\newcommand{\Q}{\bbQ}
\newcommand{\N}{\bbN}

\newcommand{\affY}{Y_{\hbar_1, \hbar_2, \hbar_3}(\widehat{\mathfrak{gl}(1)})}

\DeclareMathOperator{\fac}{fac}
\DeclareMathOperator{\pr}{pr}

\DeclareMathOperator{\BM}{BM}

\newcount\cols
{\catcode`,=\active\catcode`|=\active
 \gdef\Young(#1){\hbox{$\vcenter
 {\mathcode`,="8000\mathcode`|="8000
  \def,{\global\advance\cols by 1 &}%
  \def|{\cr
        \multispan{\the\cols}\hrulefill\cr
        &\global\cols=2 }%
  \offinterlineskip\everycr{}\tabskip=0pt
  \dimen0=\ht\strutbox \advance\dimen0 by \dp\strutbox
  \halign
   {\vrule height \ht\strutbox depth \dp\strutbox##
    &&\hbox to \dimen0{\hss$##$\hss}\vrule\cr
    \noalign{\hrule}&\global\cols=2 #1\crcr
    \multispan{\the\cols}\hrulefill\cr%
   }
 }$}}
}

\newtheorem{thm}{Theorem}[subsection]
\newtheorem{defn}[thm]{Definition}
\newtheorem{lmm}[thm]{Lemma}
\newtheorem{rmk}[thm]{Remark}
\newtheorem{prp}[thm]{Proposition}
\newtheorem{conj}[thm]{Conjecture}
\newtheorem{exa}[thm]{Example}
\newtheorem{cor}[thm]{Corollary}

\newcommand{\C}{{\bf C}}
\newcommand{\R}{{\bf R}}
\newcommand{\Z}{{\bf Z}}

\newcommand{\A}{A_{\infty}}
\newcommand{\CC}{{\mathcal C}}

\renewcommand{\H}{{\mathcal{H}}}

\newcommand{\W}{{\mathcal{W}}}

\newcommand\OO{{\mathcal O}}

\newcommand{\Gg}{{\mathsf{G}}}

\title{Cohomological Hall algebras, vertex algebras and instantons}

\author {Miroslav Rap\v{c}\'{a}k, Yan Soibelman, Yaping Yang, Gufang Zhao}

\begin{document}

%\author[M.~Rapcak]{Miroslav~Rap\v{c}\'{a}k}
%\address{M.R.: Perimeter Institute for Theoretical Physics, 31 Caroline St N, Waterloo, ON N2L 2Y5, Canada}
%\email{miroslav.rapcak@gmail.com}

%\author[Y.~Soibelman]{Yan~Soibelman}
%\address{Y.S.: Department of Mathematics, Kansas State University, Manhattan, KS 66506, USA}
%\email{soibel@math.ksu.edu}

%\author[Y.~Yang]{Yaping~Yang}
%\address{Y.Y.: School of Mathematics and Statistics, The University of Melbourne, 813 Swanston Street, Parkville VIC 3010, Australia}
%\email{yaping.yang1@unimelb.edu.au}

%\author[G.~Zhao]{Gufang~Zhao}
%\address{G.Zh: Institute of Science and Technology Austria,
%Am Campus, 1,
%Klosterneuburg 3400,
%Austria}
%\email{gufang.zhao@ist.ac.at}

\maketitle

\begin{abstract} We define an action of the (double of) Cohomological Hall algebra of Kontsevich and Soibelman
on the cohomology of the moduli space of spiked instantons of Nekrasov. We identify this action
with the one of the affine Yangian of $\mathfrak{gl}(1)$. Based on that we derive the vertex algebra 
at the corner $\W_{r_1,r_2,r_3}$ of Gaiotto and Rap\v{c}\'{a}k. We conjecture that our approach 
works for a big class of Calabi-Yau categories, including those associated with toric Calabi-Yau $3$-folds.
\end{abstract}

\setcounter{tocdepth}{2}
% Change the above number "2" to "1" if we want to hide subsections in the table of contents. 
\tableofcontents

\section{Introduction}

\subsection{Motivations from mathematics}

Nakajima's construction of an action of the infinite Heisenberg algebra on the 
(equivariant) cohomology of the moduli space 
of $U(1)$-instantons on $\C^2$ (see \cite{Nak97}) is an archetypical example of 
description of a Lie algebra via its action by correspondences on the 
cohomology of a moduli space. This branch of geometric representation theory 
has  long history. 

Nakajima's result can be also interpreted as a geometric description (or 
definition 
if one prefers) of the $\W$-algebra of the affine Lie algebra $\widehat{\mathfrak{gl}(1)}$.

The case of $\W$-algebras $\W_r:=\W(\widehat{\mathfrak{gl}(r)})$ was 
discussed recently in 
several papers in relation to the proof of the AGT conjecture. The closest to our point of view are \cite{MO, SV2}. 
%\footnote{The paper was 
%started in July 2016. Since that time several physics papers containing similar 
%ideas have appeared. The most relevant are \cite{GR} and \cite{PR}.}

Main  motivation for our project is the growing role in QFT of the notion of Cohomological Hall algebra (COHA for short) introduced in  2010 in \cite{KS}. Originally, it was considered by the authors  as the mathematical incarnation of the notion of multiparticle BPS algebra. The idea of the algebra structure on closed BPS states goes back to Harvey and Moore  (see e.g. \cite{HaMo}).  Notice that differently
from the original expectations of \cite{HaMo} the associative algebra structure proposed in \cite{KS} exists only on multiparticle BPS states and does not depend on the central charge of the stability structure.
In order to derive the Lie algebra structure on single-particle BPS states one has to  do more work
(see the original conjecture in \cite{KS} and  further developments in \cite{Dav1}). Then COHA ``looks like" the universal enveloping algebra of this Lie algebra. 

COHA was introduced originally in the framework of
Quillen-smooth associative algebras with potential. It was mentioned in \cite{KS} (see e.g. Section 3.5 in the loc. cit.)  that it can be defined in fact for a sufficiently general class of 
$3$-dimensional Calabi-Yau categories ($3CY$ categories for short) and an appropriate cohomology theory of constructible  dg-stacks. This is sufficient for most  of the applications, including e.g. knot theory (see \cite{S}, \cite{GuSto}).

Roughly, COHA is given by an algebra structure on the cohomology of the stack of objects of a $3CY$-category with coefficients in the sheaf of vanishing cycles of its potential.  Motivic Donaldson-Thomas invariants ($DT$-invariants for short)  introduced in \cite{KS2} and revisited in \cite{KS} can be defined in terms of the virtual Poincar\'e polynomial (a.k.a Serre polynomial) of COHA. Relation to  $DT$-invariants explains the initial
interest of physicists in COHA. Surprisingly the {\it algebra} structure of COHA was not seriously used in physics. 

In this paper we illustrate
the general hope that some (maybe all?) quantum algebras which appear in the interplay between $4d$ gauge theories and $2d$ CFTs come from  appropriate Cohomological Hall algebras. \footnote{We use the term ``Cohomological" even in the case
when we are talking about  versions for $K$-theory or  any other generalized cohomology theory.}
This class of algebras includes as 
special cases different versions  of Yangians as well as $\W$-algebras.
If one keeps in mind that COHA is related to factorization algebras  (see \cite{KS})
then its role in the String Theory and Gauge Theory seems to be universal (see \cite{CosGw}).

More precisely, COHA 
acts naturally on the cohomology (or $K$-theory for the $K$-theoretical version) of the moduli spaces of stable configurations of geometric and algebraic objects
which appear  in the Gauge Theory, very much in the spirit of Nakajima's seminal paper mentioned at the beginning.
The relation between instanton partition functions in $4d$ gauge theories and $D$-branes on Calabi-Yau $3$-folds is the natural
framework in which the representation theory of COHA should appear (cf. \cite{ S, Sz2}). In this paper we will illustrate this idea
in the case of spiked instantons introduced by Nekrasov (see \cite{N, Nek1, Nek3, NP,Ne4,Ne5}).
Mathematically, we will follow the  approach to the representation theory of COHA proposed in \cite{S}
as well as  its generalizations (the reader may find \cite{Fra} useful as well).

From the point of view of  BPS algebras
we will consider the action of the BPS algebra 
of  $D0$-branes in $\C^3$ on the equivariant cohomology of the moduli space of 
stable framed $D0-D4$ states (generalization to stable framed $D0-D2-D4$ is 
also possible for more general  Calabi-Yau manifolds). More precisely,  we 
consider  only $D4$ branes which correspond to  coherent sheaves on $\C^3$ 
supported on the  divisor $\sum_{1\le i\le 
3}r_i\C^2_i$. 
Here $r_i\ge 0$ are given integers, and $\C_i^2, 1\le i\le 3$ 
are the coordinate planes $\{z_i=0\}\subset \C^3$, 
where the vector space $\C^3$ is endowed with the standard 
coordinates $(z_1,z_2,z_3)$. Equivalently,  we have a toric divisor $z_0z_1^{r_1}z_2^{r_2}z_3^{r_3}=0$ in $\C{\bf P}^4$
intersected with the plane $z_0=1$. 

From this point of view it is natural to study more general toric diagrams and corresponding local Calabi-Yau $3$-folds. 
More generally one can start with an arbitrary dimer model, provided the $2$-dimensional faces are ``colored'' by non-negative integers $r_i$. These integers correspond to the ranks of gauge groups. In particular one
can hope for the following result.

\begin{conj} With a toric Calabi-Yau $3$-fold $X$ and a collection of non-negative integers $r_1,...,r_k$ (ranks) assigned to the $2$-dimensional faces of the toric diagram of $X$, one can associate a  VOA
$\W_{X, r_1,...,r_k}$ which we will call the $\W$-algebra associated to $X$ and the ranks $r_i, 1\le i\le k$ . 
\end{conj}

%More precisely, let $(Q,W)$ denote the quiver with potential associated with the toric diagram of $X$, and $(Q^{fr}, W^{fr})$
%denote the corresponding framed quiver with framed potential. The definition of the former is well-known (in fact a pair $(Q,W)$ can be associated with any dimer model.  
%In order to define the corresponding framed quiver, we add a new vertex  for each $2$-dimensional face of the toric diagram, and pair of opposite arrows $i_v, j_v$ for each vertex $v$ of $Q$ (here $i_v$ starts
%at the framing vertex and ends at $v$, while $j_v$ is opposite$). 

If $X=\C^3$ and $D=\sum_{1\le i\le 
3}r_i\C^2_i$, then $\W_{\C^3, r_1,r_2,r_3}:=\W_{r_1,r_2,r_3}$ coincides with the vertex algebra 
$Y_{r_1,r_2,r_3}$ defined in \cite{GR}.  It was later used in \cite{PR,PR2} as a building block
for construction of more complicated vertex algebras
associated with toric $CY$ $3$-folds via the gluing procedure proposed in the loc. cit. This gluing algorithm corresponds to the mathematical notion of conformal extension of VOAs (see for example \cite{F}).

The  vertex algebra $\W_{r_1,r_2,r_3}$ is isomorphic to the quotient of the  vertex algebra $\W_{1+\infty}$ by the 
$2$-sided ideal which corresponds to the curve 
$$\sum_{1\le i\le 3}{r_i\over{\lambda_i}}=1$$
in the space of triples $(\lambda_1,\lambda_2,\lambda_3)$ such that
$$\sum_{1\le i\le 3}{1\over{\lambda_i}}=0$$ 
used in \cite{P1,P2} to parametrize $\W_{1+\infty}$ algebras. For $r_1=r, r_2=r_3=0$ the quotient is isomorphic to the $\W$-algebra 
$\W_r:=\W(\widehat{\mathfrak{gl}(r)})$.

There are at least two ways to approach the construction of $\W_{X, r_1,..., r_k}$ for general toric $X$. We have already mentioned the approach of \cite{PR} which is analogous  to the topological vertex formula in the sense that one 
constructs  vertex algebras  associated to the colored toric diagrams 
starting with the basic one for $\C^3$. The 
idea of gluing complicated vertex algebras from the basic ones 
was discussed in the case of toric {\it surfaces} for example in \cite{DeGuPu} and \cite{FGu}.

Alternatively one can construct $\W_{X, r_1,..., r_k}$ starting with the Cohomological Hall algebras $\H_X$ of the Calabi-Yau $3$-fold $X$. In that case  $\W_{X, r_1,..., r_k}$ is defined as the commutant of the subalgebra of screening operators acting on a highest weight $\H_X$-module.
In this  paper we use this  approach in the  case of $\W_{r_1,r_2,r_3}$.

Speculating further one can hope that there exists a vertex algebra   associated to a general ind-constructible locally 
ind-Artin three-dimensional Calabi-Yau category endowed with a stability 
structure and a framing (see \cite{KS2}, \cite{S} 
for the background 
material). 

\begin{conj} With any  ind-constructible locally ind-Artin 
$3CY$-category endowed with a stability structure  one can associate a class of 
vertex algebras, parametrized by rays in $\C=\R^2$ and a choice of framing 
object (for almost all rays the algebra is trivial).

\end{conj}

The point is that such a category with stability structure gives rise to a COHA. Its ``double'' should act on the moduli space of stable framed objects associated with
a given ray, given framing object and the stability structure. The hope is that this
action underlies a structure of vertex algebra.\footnote{In this paper we use words ``vertex algebras" and ``vertex operator algebras"
synonymously.}

The above discussion naturally leads to a question about different versions of the AGT conjecture. The conventional one proved
for $\W_r$ corresponds to the case $X=\C^3, r_1=r, r_2=r_3=0$.

From the perspective of COHA the appearance of the moduli space of instantons in the conventional AGT conjecture
is due to the ``dimensional reduction'' proposed in \cite{KS}, Section 4.8 (see also a detailed discussion in \cite{D}). Recall that the dimensional reduction allows one to replace COHA associated  with a $3$-dimensional Calabi-Yau 
category  by a simpler algebra associated with a $2$-dimensional Calabi-Yau category.

\begin{rmk} We would like to make a comment about the terminology and clarify the confusion in the literature.
In \cite{SV2} and 
subsequent papers of the same authors they used the term
``Cohomological Hall algebra'' for a special case of general COHA. The former is the ``dimensional reduction" of the one introduced in \cite{KS}. In \cite{YZ16} the same object is called
preprojective COHA.

In order to avoid the confusion and some conflict in terminology we propose to 
call the COHA introduced in the foundational paper \cite{KS} by {\bf 3d COHA}, while its special case  considered
in \cite{SV2,YZ16} will be called {\bf 2d COHA}. The terminology is justified by the observation that in \cite{KS} 
the authors defined COHA in the framework of $3$-dimensional Calabi-Yau categories, 
while in \cite{SV2} the authors deal with a special case of a $2$-dimensional Calabi-Yau category ($2CY$-category for short). 

In the current paper we will often  call the $3d$ COHA simply by  COHA keeping the term $2d$ COHA for
the special case considered in \cite{SV2, YZ16}.

We also remark that although in this paper we use several techniques of \cite{SV2,YZ16}, where the $2d$ COHA was utilized, they should not be sufficient
in general, especially for non-toric CY 3-folds. In those cases the $3d$ COHA  will be necessary to use. 

\end{rmk}

\subsection{About generalization of COHA  to $4$-dimensional Calabi-Yau categories}

It is important to have a  generalization of COHA to a class of $4$-dimensional Calabi-Yau varieties. The main example
should be worked out in the case of $\C^4$. It will allow us to accommodate more general gauge theories introduced by
Nekrasov (see \cite{N, Nek1, Nek3, NP,Ne4,Ne5}). That story  involves non-holomorphic ADHM-type relations. In the current paper we discuss a special case related to the  moduli spaces of {\it spiked instantons}.  Then we deal with $\C^3$ instead of $\C^4$ and all the relations are holomorphic.

It is not clear at the moment how to define an  analog of COHA for  an appropriate class of $4$-dimensional 
Calabi-Yau categories. One can hope to do that in the framework  \cite{N, Nek1, Nek3, NP,Ne4,Ne5} which is  special, since one can still define the perfect obstruction theory.  It gives one a hope that there is a $4$-dimensional COHA defined in terms of the cohomology with coefficients in a constructible sheaf
on the  moduli stack of objects of a $4$-dimensional Calabi-Yau category. 

Let us briefly discuss the origin of the perfect obstruction theory in Nekrasov's story.
Recall that the  $qq$-characters of Nekrasov (see  \cite{N, Nek1, Nek3, NP,Ne4,Ne5}) are defined as integrals of the generating functions of Chern classes of
some natural vector bundles on the moduli spaces of solutions of generalized ADHM equations introduced in the loc. cit. An
important part of the story is that one can obtain a real virtual fundamental class (in fact of odd dimension).

A more general set generalizing the one of  \cite{NP,Nek1,N,  Nek3, Ne4,Ne5} should be the $4$-dimensional non-compact toric Calabi-Yau 
manifold with a toric subscheme. More precisely, the manifold considered in the 
loc.cit. is $\C^4$, 
and the subscheme is given by a collection
of six coordinate planes $\C^2_{ij}, 1\le i,j\le 4$ with multiplicities 
$r_{ij}\in \Z_{\ge 0}$. 

Existence of the virtual fundamental class over which one 
should perform 
the integration is not obvious. For a $4$-dimensional CY 
category endowed with a stability condition,  the moduli space of stable (or 
polystable)
objects of the heart of a $t$-structure does not have to support  a perfect obstruction theory. Indeed
besides of the group $Ext^1(E,E)$ which gives the tangent space $T_E$ at an object $E$ one has 
also groups $Ext^i(E,E), i\ge 2$. Differently from the case of 
$3CY$ categories those groups do not
cancel each other in the virtual tangent space, 
since now we have the Serre duality in dimension $4$ rather than $3$. Notice that by  Serre 
duality, the vector space $Ext^2(E,E)$ carries a non-degenerate complex-valued 
bilinear form $(\bullet,\bullet)$. Let us  choose a {\it real} vector subspace 
$Ext^2_{\R}(E,E)\subset Ext^2(E,E)$ such that the restriction of 
$(\bullet,\bullet)$ to $Ext^2_{\R}(E,E)$ is positive. Then we can obtain a 
perfect obstruction theory associated with the vector
space $Ext^0(E,E)\oplus Ext^1(E,E)\oplus Ext^2_{\R}(E,E)\oplus 
(Ext^2_{\R}(E,E))^{\ast}$. 
The corresponding real vector bundle $(Ext^{\bullet}(E,E),(\bullet,\bullet))$ 
over the moduli space of (poly)stable objects $E$ should be oriented for the 
virtual integration purposes. There is a  condition on $det(Ext^{\bullet}(E,E), 
 (\bullet,\bullet))$ which imposes topological restrictions on the underlying 
moduli space. Similarly to the notion of orientation data in \cite{KS2} the 
determinant  itself is a tensor square, but the choice of square root is an 
additional piece of data (see \cite{BorJoy}, \cite{CaLeu} for some relevant constructions).

\begin{rmk} Technically, the above considerations hold only for compact 
Calabi-Yau manifolds (or CY categories with the compact space of (poly)stable 
objects. This is not the case for $\C^4$ or any toric Calabi-Yau $4$-fold. 
In such cases 
one can use the equivariant version of the above considerations 
and observe that the set of the fixed points of the natural action of the torus  on the 
moduli space of Nekrasov instantons on $\C^4$ is compact. When there are
toric divisors consisting of coordinate planes, one can use 
the $(\C^{\ast})^3$-action on $\C^4$ obtained from the natural
$(\C^{\ast})^4$-action by imposing the condition that the product of weights of 
the action is equal to $1$. This ensures that the standard holomorphic 
Calabi-Yau form on $\C^4$ is preserved. In the end one obtains a real 
equivariant virtual fundamental class which can be used for equivariant 
integration as in \cite{NP,Nek1,N,  Nek3, Ne4,Ne5}.

\end{rmk}

\subsection{Motivations from physics}

Let us rephrase the above discussion in a slightly more physical language. We start by reviewing the physics of the BPS algebra \cite{HaMo,KS} associated 
to Calabi-Yau $3$-folds and its relation to VOAs. Motivated by the standard Alday-Gaiotto-Tachikawa \cite{AGT,W} setup, 
we propose a generalization of the AGT correspondence for spiked instanton configurations of \cite{NP,Nek1,N,  Nek3}  restricted to branes in toric Calabi-Yau
$3$-folds. We finish the introduction by some speculations related to more general configurations.

\subsubsection{BPS algebra and VOA}

A rich class of examples of BPS algebras from \cite{HaMo} should arise from the compactification of type IIA string theory on a Calabi-Yau $3$-fold $X$, 
i.e. superstrings in ${\R}^4\times X$. BPS particles are then associated to $D6-D4-D2-D0$ branes wrapping compact complex cycles inside $X$ and spanning 
one extra direction inside ${\R}^4$. The corresponding BPS algebra is an algebra capturing dynamics of such BPS particles. The physically motivated notion of
BPS algebra was mathematically formalized in  \cite{KS} in the notion of the Cohomological Hall algebra.

Turning on a fixed background of $D6$ and $D4$ branes
(possibly including branes supported on non-compact cycles), one can look at the subalgebra of BPS particles associated to $D2-D0$ branes supported on 
compact cycles. The corresponding configuration can be studied from two different perspectives. First, from the perspective of $D6$ and $D4$ branes, the 
compact $D2-D0$ branes introduce a gauge-field flux and correspond to non-trivial instanton configurations. COHA is then expected to relate configurations 
of different instanton numbers. Secondly, from the perspective of the compact $D2-D0$ branes, the moduli space of instantons has a description in terms of the moduli 
space of vacua of a quiver quantum mechanics. COHA then relates vacua of theories of different ranks.

From the above perspective, it is natural to expect that the equivariant cohomology of the moduli space of instantons of  the $5d$ and $7d$ theories on $D4$-
and $D6$-branes carries the structure of a representation of an appropriate COHA (or its subalgebra associated to the dynamics of $D2-D0$ branes).

 In this paper, we will be mostly interested in configurations concerning non-compact $D4$-branes (see e.g. \cite{AJS,J,NiYaYo}). Lifting the type IIA brane setup to
M-theory by introducing an extra M-theory circle, $D4$-branes lift to $M5$-branes wrapping also the extra circle. This lift leads to the standard configuration
of the $2d-4d$ correspondence associated to $M5$-branes wrapping a complex two-dimensional variety $M_4$ and an extra Riemann surface 
$M_2=S^1\times {\R}$. Roughly speaking, compactification of the M5-brane theory on $M_2$ leads to a four-dimensional gauge theory supported on $M_4$ whereas 
compactification on $M_4$ is expected to lead to a two-dimensional CFT for compact $M_4$ or a chiral algebra for non-compact $M_4$. The duality between these two 
theories is known as the AGT correspondence, $2d-4d$ correspondence or the BPS/CFT correspondence \cite{AGT,FGu,GGP,Nek1}. 
The corresponding VOA[$M_4$] appears as an algebra of chiral operators in the $2d$ CFT, generating symmetries of the theory 
and extending the Virasoro algebra generating conformal transformations.

From the 2d perspective, the BPS algebra of $D0$- and $D2$-branes gives rise to operators on the Hilbert space of the theory that 
can be identified with the equivariant cohomology of the moduli space of instantons 
associated to the divisor $M_4$. It is natural to expect that COHA acting on the equivariant cohomology of the moduli space actually 
leads to the vertex operator algebra  VOA[$M_4$]. The equivariant cohomology
of the moduli space of instantons can be then identified with a generic module for VOA[$M_4$]. 
The Calabi-Yau perspective and the corresponding COHA then provides a way to universally study VOA[$M_4$] for a large class of complex surfaces $M_4$ associated to different divisors  in a 
Calabi-Yau $3$-fold. Let us now discuss some examples starting with the well-known story of $M_4=\C^2$ in the $\Omega$-background and moving towards more exotic configurations.   

\subsubsection{Standard AGT setup}

The simplest example of the above setup is the configuration of Alday-Gaiotto-Tachikawa \cite{AGT,W} relating the Nekrasov partition function \cite{N6} of a $U(r)$ gauge theory on $M_4=\C^2$ in the $\Omega$-background with conformal blocks of the $\mathcal{W}_r$ algebra\footnote{We use the notation $\mathcal{W}_r=\mathcal{W}(\widehat{\mathfrak{gl}(r)})$, i.e. the $\mathcal{W}$-algebra associated to the principal embedding of $\mathfrak{sl}(2)$ inside $\mathfrak{gl}(r)$, instead of $\mathcal{W}_r=\mathcal{W} (\widehat{\mathfrak{sl}(r)})$ used in some of the literature. These two differ by a factor of $\widehat{\mathfrak{gl}(1)}$.} on $M_2$. This configuration can be simply embedded inside our setup by considering $r$ M5-branes wrapping $\C^2_{x_1,x_2}$ inside the Calabi-Yau $3$-fold $\C^3_{x_1,x_2,x_3}$ in the presence of a $B$-field with equivariant parameters $\epsilon_1$, $\epsilon_2,\epsilon_3=-\epsilon_1-\epsilon_2$ associated to the rotations of the $\C_{x_i}$ planes.

A key step in the proof of the AGT correspondence is the construction of the action of $\mathcal{W}_r$ on the equivariant cohomology of the moduli space of instantons on $\C^2$ with equivariant parameters $\epsilon_1, \epsilon_2$ \cite{MO, SV2}. The moduli space has an alternative description in terms of representations of the ADHM quiver \cite{ADHM,Dou1,Dou2}. Physically, the ADHM quiver can be viewed as a quiver of the $U(n)$ gauge theory on $n$ $D0$-branes bound to the stack of $r$ $D4$-branes. The ADHM moduli
can be then identified with the moduli space of vacua of such a supersymmetric quiver quantum mechanics. The dual perspective in terms of a type IIB configuration (to be described below) from \cite{GR,PR,PR2} indeed associates the algebra $\mathcal{W}_r$ to such a simple divisor.

\subsubsection{The $\C^3$ example} 

\begin{figure}
  \centering
      \includegraphics[width=0.54\textwidth]{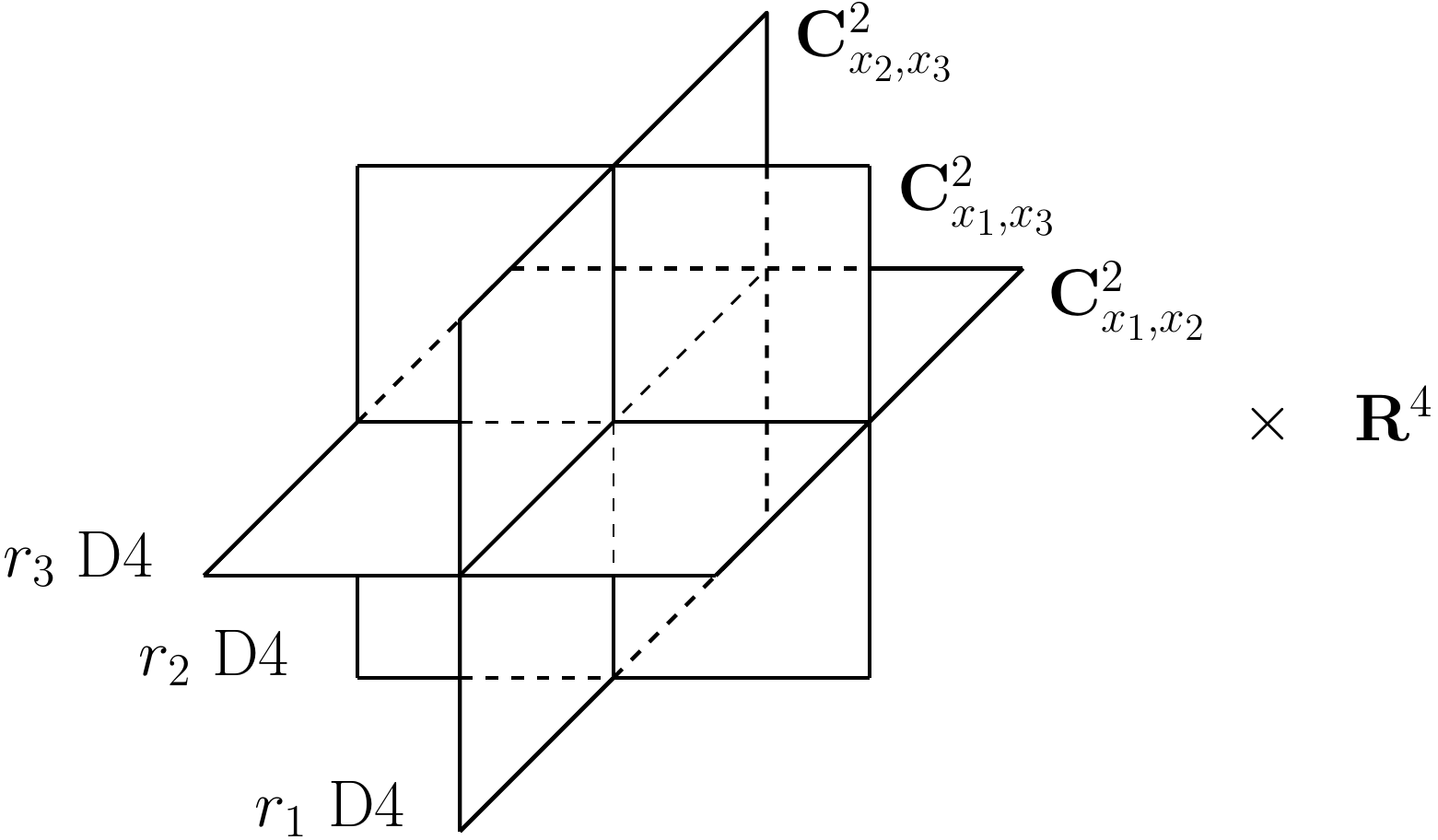}
  \caption{Configuration of branes in type IIA string theory associated to 
spiked instantons. Three stacks of D4-branes span the three four-cycles inside 
$\C^3$ fixed under the $T^2$ action discussed in the main text together with one orthogonal direction in $\R^4$.}
\label{typeIIA}
\end{figure}

\begin{figure}
  \centering
       \includegraphics[width=0.29\textwidth]{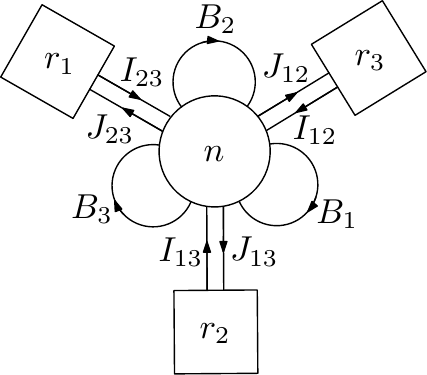}
  \caption{Quiver for the spiked-instanton configuration.}
\label{quiver}
\end{figure}

The configuration above has a natural generalization from the three-dimensional perspective. One can consider three 
stacks of $r_1,r_2$ and $r_3$ $M5$-branes supported on $\C^2_{x_1,x_2}$, $\C^2_{x_1,x_3}$ and $\C^2_{x_2,x_3}$ inside 
$\C^3_{x_1,x_2,x_3}$, i.e. a configuration associated to the divisor
$r_1\C^2_{x_2,x_3}+r_2\C^2_{x_1,x_3}+r_3\C^2_{x_1,x_2}$ 
 with $r_i\in \Z_{\ge 0}$. Compactifying on the extra Riemann surface
shared by all the branes, one obtains a triple of four-dimensional
$U(r_i)$ gauge theories supported at the three irreducible components of the divisor, namely on
$\C^2_{x_1,x_2}$, $\C^2_{x_1,x_3}$ and $\C^2_{x_2,x_3}$, mutually interacting along
their intersections $\C_{x_1}$, $\C_{x_2}$, $\C_{x_3}$ via bi-fundamental 2d fields. This setup can be identified with a restriction of the more general
spiked-instanton setup of $M5$-branes intersecting inside $\C^4$ from \cite{Nek1, N, Nek3}. 
The corresponding moduli space of instantons has a quiver description from the figure \ref{quiver} that can be again derived as the moduli space of vacua for $D0$-branes bound to $D4$-branes from the figure \ref{typeIIA} as shown in \cite{NP}. In the presence of a single stack of $M4$-branes, e.g.  $r_1=r_2=0$, the quiver reduces to the standard ADHM quiver with only two loops and a single framing node of rank $r_3$. 

The above M-theory setup can be related along the lines of \cite{LV,NW} to the configuration from \cite{GR}
using the duality between the M-theory on a torus and type IIB string theory in the presence of 
a web of $(p,q)$-branes. In the example at hand,  $\C^3=\R^6$ endowed with the
standard symplectic structure  has the natural Hamiltonian action of $T^3=U(1)^3$, whose moment map realizes $\C^3$ as a singular Lagrangian $T^3$-
fibration over the first octant in ${\R}^3$. The action of the 2-dimensional subtorus $T^2\subset T^3$ preserving the canonical bundle is generated by the following rotations $(e^{it_1}z_1,z_2,e^{-it_1}z_3)$
and $(z_1,e^{it_2}z_2,e^{-it_2}z_3)$. The moment map of this $T^2$ action from $\C^3$ to ${\R}^2$ is given by
$\mu_1=|z_1|^2-|z_3|^2$ and $\mu_2=|z_2|^2-|z_3|^2$. 
The directions in which the $T^2$ torus fibration, when projected to ${\R}^2$, are as follows.
The $t_1$ action degenerates for $z_1=z_3=0$,
corresponding to the $\mu_1=0,\mu_2>0$, the $t_2$ action degenerates at $\mu_2=0$ and $\mu_1>0$ and
finally $t_1+t_2$ degenerates at $\mu_1=\mu_2<0$. The degeneration of the fibers in the $\mu_1,\mu_2$ plane
is shown in the figure \ref{gluing} on the left. 

\begin{figure}
  \centering
      \includegraphics[width=0.65\textwidth]{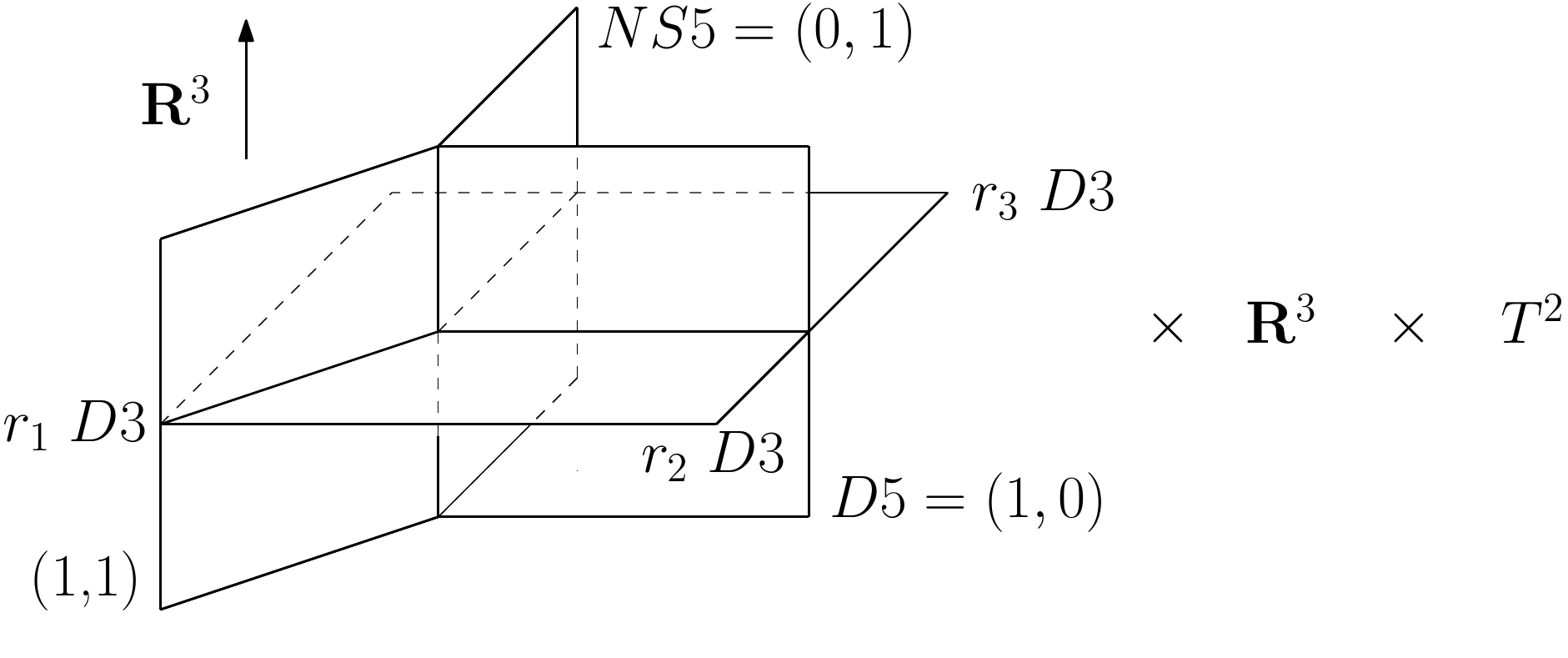}
  \caption{Configuration of branes in type IIB string theory leading to 
$\W_{r_1,r_2,r_3}$. $NS5$, $D5$ and $(1,1)$ branes span directions 01456 and a line inside the 
23-plane as shown in the figure. Three stacks of $D3$-branes are attached to the 
five branes from the three corners and span the extra 01 directions. The two dimensional 
junction of all the branes supports the $\W_{r_1,r_2,r_3}$.}
\label{typeIIB}
\end{figure}

From the dual point of view, one gets a type IIB theory on $\R^8\times T^2$ with one of the cycle $S_1\subset T^2$ coming from the 
toric fibration of the Calabi-Yau $3$-fold and the other cycle corresponding to the M-theory circle $S_2\subset T^2$ from the lift of the 
type IIA configuration discussed above.
Singularities of the torus fibration correspond to $(p,q)$-branes spanning orthogonal directions with $p$ and $q$ labeling the
 degenerating circle. The geometry of the Calabi-Yau $3$-fold thus maps to a web of $(p,q)$-branes.
$M5$-branes associated to the faces in the toric diagram map to $D3$-branes attached to $(p,q)$-branes from the 
three corners as shown in the figure \ref{typeIIB}. This is exactly the setup of \cite{GR} that identified a three-parameter family of VOAs 
$\mathcal{W}_{r_1,r_2,r_3}$ as an algebra of local operators at a two-dimensional junction of interfaces in the four-dimensional theory 
coming from the low-energy dynamics of the D3-branes.

Note that $\mathcal{W}_{r_1,r_2,r_3}$ becomes the standard $\mathcal{W}_{r_i}$
algebra if two of the remaining parameters $r_j=0$ for $j\ne i$ vanish. Since the quiver description of spiked instantons reduces
to the standard ADHM quiver in this case, the duality explains the standard AGT correspondence. Motivated by the above duality 
and this observation, it is natural to expect that $\mathcal{W}_{r_1,r_2,r_3}$ should act on the equivariant cohomology of the 
moduli space of spiked instantons in a greater generality. Most of the paper is devoted to the proof of this proposal.

\begin{rmk}
Technically speaking, the modules (the generic modules) we construct in the present paper are not those (the specialized modules) from Gaiotto-Rapcak \cite{GR},
but they are different type of module over the same vertex algebra $\mathcal{W}_{r_1,r_2,r_3}$.
The generic module in the present paper coincides with the free field realization of  $\mathcal{W}_{r_1,r_2,r_3}$ in \cite{PR2}. 
The basis of the generic module is expect to be given by $(r_1+r_2+r_3)$-tuples of planar partitions, whose character would be a product of eta functions. 
On the other hand, the specialized module from \cite{GR} is defined in terms of some quantum Hamiltonian reduction of Kac-Moody algebras, which resulting the "maximally degenerate" modules of $\mathcal{W}_{r_1,r_2,r_3}$. 
The basis of the specialized module  is expect to be given by the 3d partitions that restricted to lie in a certain region, whose character would be the 
MacMahon function when $r_1, r_2, r_3\to \infty$. 
\end{rmk}

\subsubsection{Further remarks and possible generalizations}

\begin{figure}
  \centering
      \includegraphics[width=1\textwidth]{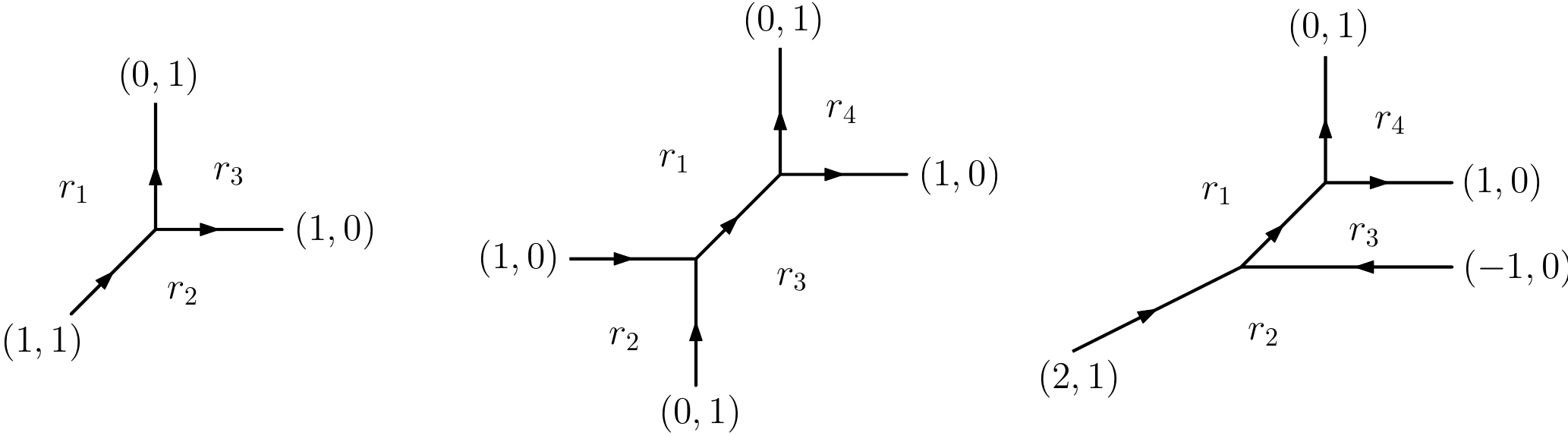}
  \caption{Toric diagram associated to ${\C}^3$ (left), $\mathcal{O}(-1)\oplus \mathcal{O}(-1)\rightarrow {\bf P}^1$ (middle) and  $\mathcal{O}(-2)\oplus \mathcal{O}\rightarrow {\bf P}^1$ (right). The lines show loci where $(p,q)$-cycles of the torus $T^2$ degenerate.}
\label{gluing}
\end{figure}

The example of $\C^3$ has a natural generalization for an arbitrary toric Calabi-Yau $3$-fold given by a toric diagram
specifying loci where the torus cycles degenerate. The two simplest examples are shown in the figure \ref{gluing} and correspond to the bundles $\mathcal{O}(-1)\oplus \mathcal{O}(-1)\rightarrow {\bf P}^1$ and  $\mathcal{O}(-2)\oplus \mathcal{O}\rightarrow {\bf P}^1$ respectively. From the type IIB
perspective, the toric diagram can be again identified with  a nontrivial web of $(p,q)$-branes. M5-branes associated to each
face in the diagram map to stacks of D3-branes attached to the $(p,q)$-branes from various corners. Vertex operator algebras
associated to such configurations were identified in \cite{PR,PR2} with various extensions\footnote{Note also closely related gluing at the level of affine Yangians \cite{GLPZ,GLP}, quantum toroidal algebras \cite{FJMM} or minimal models \cite{HM}.} of tensor products of $\mathcal{W}_{r_1,r_2,r_3}$ associated to each
trivalent junction of the diagram by bi-fundamental fields associated to internal lines of the diagram. 

The free field realization of $\mathcal{W}_{r_1,r_2,r_3}$ and their modules from \cite{F,PR2} gives an explicit realization of the bi-fundamental fields in terms of exponential vertex operators of the free boson. This generalizes the well-known constructions of lattice vertex operator algebras \cite{B,FLM} from the $\mathcal{W}_{0,0,1}=\widehat{\mathfrak{gl}(1)}$  case to an arbitrary $\mathcal{W}_{r_1,r_2,r_3}$ with the necessity to include contour integrals of screening charges along the lines of \cite{DF,Fel}. 

The duality reviewed above suggests a generalization of the AGT correspondence for divisors in toric Calabi-Yau $3$-folds which agrees with our Conjecture 1.1.1. The gauge theories supported on smooth components of the divisor $D$ fixed under the toric action and interacting along the intersection of the smooth components correspond to $(p,q)$-web VOAs. 
In particular, 
based on the analysis of the vacuum character of the glued algebras from \cite{PR}, we can conjecture that the Drinfeld double of the equivariant spherical COHA associated to the last two configurations in figure \ref{gluing} can be identified with shifted affine Yangians of $\mathfrak{gl}(1|1)$ and $\mathfrak{gl}(2)$ respectively. Shifts are determined by the intersection number of the corresponding divisor $D$ with the ${\bf P}^1$ associated to the internal line.

Note also a different proposal of \cite{FGu} for the algebras associated to the two examples above. \cite{FGu} conjectures both configurations to lead to different shifts of the $\widehat{\mathfrak{gl}(2)}$ Yangian. 
On the other hand, calculations from \cite{PR} suggest an appearance of the Yangian of $\widehat{\mathfrak{gl}(1|1)}$ in the $\mathcal{O}(-1)\oplus \mathcal{O}(-1)\rightarrow {\bf P}^1$ case and the Yangian of $\widehat{\mathfrak{gl}(2)}$  in the $\mathcal{O}(-2)\oplus \mathcal{O}\rightarrow {\bf P}^1$ case.  The shifts are determined by the intersection number of the corresponding ${\bf P}^1$ with the divisor. This can be seen from the calculation of vacuum characters in the large $r_i$ limit and the appearance of the $\widehat{\mathfrak{gl}(2)}$ and $\widehat{\mathfrak{gl}(1|1)}$ Kac-Moody algebras for a particular choice of divisor in the two cases.

\subsection{Contents of the paper}
Let us briefly discuss the contents of the paper.

In \S~\ref{sec:COHA_tor_sheaves}, we discuss the relation of COHA in our main example of the quiver with potential  and its relation to torsion sheaves on $\C^3$. We also introduce the equivariant and spherical versions of COHA. In 
\S~\ref{sec:framed_quiver} we discuss a particular framing of our quiver with potential as well as the corresponding framed stable representations. In \S~\ref{sec:Main result} after a more technical reminder on  COHA, we state our main result, and discuss the strategy we use to prove the main result. This part can also be considered as an over view of the rest of the paper. 
In \S~\ref{sec:CritCOHA_act}, we show that the COHA naturally acts on the cohomology of the moduli space of  stable framed representations.

In \S~\ref{sec:CartanDrinfeld}, we define a Cartan algebra, which acts on the COHA and the cohomology of the moduli space of  stable framed representations via characteristic classes of the tautological bundles. The COHA, with the Cartan algebra added, admits a Drinfeld coproduct. 
In \S~\ref{sec:double} we prove that the Drinfeld double associated to the coproduct from \S~\ref{sec:CartanDrinfeld} acts on the cohomology of the moduli space of  stable framed representations. In particular, the Drinfeld double is isomorphic to the affine Yangian $\affY$

In \S~\ref{sec:free-field} we define a central extension of the double COHA, and define a more interesting coproduct on the centrally extended algebra. This coproduct geometrically comes from hyperbolic localizations on the moduli space of  stable framed representations, therefore has the usual $R$-matrix formalism as in \cite{MO}.
Finally, after recalling on the VOA at the corner in \S~\ref{sec:VOA_W}, we prove, in \S~\ref{sec:Y110} that the action of the centrally extended COHA on the  cohomology of the moduli space of  stable framed representations gives rise to the VOA. Moreover, the coproduct from \S~\ref{sec:free-field} gives the free field realizations of the VOA.

{\it Acknowledgments:} Part of the work was done when Y.Y. and G.Z. were visiting the Perimeter Institute for Theoretical Physics. The research of G.Z. at IST Austria, Hausel group, is supported by the
Advanced Grant Arithmetic and Physics of Higgs moduli spaces No. 320593 of the European Research Council. 
The research of Y.S. was partially 
supported by an NSF grant and Munson-Simu Faculty Star Award of KSU.
The research
of M.R. was supported by the Perimeter Institute for Theoretical Physics, which is in turn supported by the Government of Canada through the Department
of Innovation, Science and Economic Development and by the Province of Ontario
through the Ministry of Research, \& Innovation and Science.

We thank to I. Cherednik, K. Costello, E. Diaconescu, D. Gaiotto, 
V. Gorbounov, S. Gukov, P. Koroteev, F. Malikov, N. Nekrasov, A. Okounkov, V. 
Pestun, T. Proch\'{a}zka, J. Ren, O. Schiffmann, J. Yagi for useful discussions and correspondences. 
We thank the anonymous referee for careful reading of the paper.
Y.S. is 
grateful to IHES, MSRI and Perimeter Institute for Theoretical Physics for excellent research conditions.

\section{COHA and torsion sheaves}\label{sec:COHA_tor_sheaves}

We will work over the field of complex numbers $\C$, although many results 
hold over any field of characteristic zero.

%\subsection{$3CY$ and $2CY$ categories from quivers with potential}
\subsection{Equivaraint COHA, spherical COHA and the dimensional reduction}
\label{sec:2.1}

Technical details and formal definition of COHA will be recalled in \S\ref{subsec:app_KS_crit_coh} in detail. Here we just briefly summarize a  
few facts. 

We are going to use the following notation. For a quiver $Q$ with the set 
of vertices $I$ we denote by $\overline{Q}$ 
the {\it double quiver} obtained from $Q$ by adding an opposite arrow $a^*$ for any arrow $a$ of $Q$. 
We denote by $\widehat{Q}$ the {\it triple quiver} obtained from $\widehat{Q}$ by adding a new loop $l_i$ 
for each vertex $i\in I$ of $Q$. A {\it potential} is an element $W$ of the
vector space $\C Q/[\C Q,\C Q]$, where $\C Q$ is the path
algebra of the quiver $Q$. In other words $W$ is a cyclically 
invariant non-commutative polynomial in arrows of $Q$.

Let $W$ be a {\it potential} for $\widehat{Q}$ given by the formula
$W=\sum_{a,i}l_i[a,a^*]$, where the summation is taken over all vertices $i$ and all arrows $a$ of $Q$.

\begin{exa} 
\label{exa:QW}Let $Q=J$ be the Jordan quiver (one vertex and one loop). The triple 
quiver $Q_3:=\widehat{J}$ has three loops $B_1,B_2,B_3$.
The above potential has the form $W_3=B_3[B_1,B_2]$ (plus cyclic permutations, 
which we always skip in the notation).

\end{exa}

It is well-known that for any quiver $Q$  the pair 
$(\widehat{Q}, W)$ gives rise to a $3CY$-category 
$\CC_{(\widehat{Q}, W)}$ 
which has a $t$-structure 
with the heart given by  finite-dimensional representations of the 
Jacobi algebra $\C\widehat{Q}/\langle \partial W\rangle$
(i.e. we consider the quotient of the path algebra by the $2$-sided ideal 
$\langle \partial W\rangle$ generated by cyclic derivatives of $W$).

E.g. in the above example of the pair $(Q_3,W_3)$, the $t$-structure in 
question has the heart given by the 
category of finite-dimensional representations
of the polynomial algebra $\C[B_1,B_2,B_3]$, i.e. it is the category 
$Tors(\C^3)_0$ of torsion sheaves on $\C^3$ with zero-dimensional support.
Considering $(\C^{\ast})^3$-equivariant torsion sheaves which correspond to 
cyclic modules we see that they are
enumerated by $3d$ partitions.

The dimensional reduction (see \cite{KS}, Section 4.8) associates with 
the $3CY$-category $\CC_{(\widehat{Q}, W)}$ a $2CY$-category 
$\CC_{(\overline{Q}, \partial W)}$
which has a $t$-structure with the heart given by finite-dimensional representations
of the {\it preprojective algebra} $\Pi_Q=\C\overline{Q}/\langle \sum_a[a,a^*]\rangle$. 
Intuitively, the relation between $\CC_{(\widehat{Q}, W)}$
and $\CC_{(\overline{Q}, \partial W)}$ can be thought of as a relation between two Calabi-Yau manifolds: 
a $3$-dimensional one and a $2$-dimensional one,
so that the former is the product of the latter by an affine line. More 
generally, 
one can upgrade a (triangulated $\A$) 
category of homological
dimension $2$ to a $3CY$ category by the categorical analog of the geometric construction 
of 
taking the total space of canonical bundle. Detailed explanation of this construction requires more
of the techniques of derived algebraic geometry in the spirit of  \cite{KoSo2006} ,
which is far from the subject of this paper.

%More generally, the $3d$ critical COHA was defined in \cite{KS} \footnote{ Recall that we call it $3d$ COHA, or simply 
%COHA.} such as follows:
%$$\H^{(Q,W)}=\oplus_{\gamma\in \Z_{\ge 0}^I}\H_\gamma^{(Q,W)}=
%\oplus_{\gamma\in \Z_{\ge 
%0}^I}(H_{c,\Gg_\gamma}^{\bullet}(\M_\gamma,W_\gamma))^{\vee}.$$
%Here $\M_\gamma$ 
%is the algebraic variety of finite-dimensional representations of the quiver $Q$ of 
%dimension
%$\gamma=(\gamma^i)_{i\in I}\in \Z_{\ge 0}^I$ in the coordinate space $\prod_{i\in I}\C^{\gamma^i}$, 
%gauge group $\Gg_\gamma:=\prod_{i\in I}GL_{\gamma^i}$, where $GL_m:=GL(m,\C)$ 
%is the general linear group, acts naturally by 
%changing the basis, $W_\gamma:=Tr(W): \M_\gamma\to \C$,
%and the dual is taken to the compactly supported equivariant cohomology of 
%$\M_\gamma$ with coefficients in the sheaf of vanishing cycles
%of the regular $\Gg_\gamma$-invariant function $W_\gamma$. Alternatively we can 
%work with the rapid decay cohomology or any other cohomology theory (see \cite{KS}).
%For the convenience of the readers, we include a quick reminder of some basic facts about the critical cohomology in \S~\ref{subsec:reminder_crit_coh}

%\begin{rmk}
%1) Notice that one can consider representations of quivers with potential in any abelian category, 
%e.g. in the category of coherent sheaves on a smooth algebraic variety. 
%Although in that case the arising category is not necessary $3$-dimensional 
%Calabi-Yau, the dimensional reduction still works.

COHA $\H^{(Q,W)}$ is in general graded by dimension vectors of the representations of $Q$.
We can (and will)  consider equivariant versions $\H^{(Q,W),{\bf T}_{m+1}}$ of COHA with respect to the torus ${\bf T}_{m+1}:=(\C^{\ast})^{m+1}$ rescaling the arrows of $Q$ with arbitrary weights. Here $(m+1)$ is the number of arrows of $Q$.
This version of COHA was also introduced in \cite{KS}.  

We will also consider a ${\bf T}_m$-equivariant version $\H^{(Q,W),{\bf T}_m}$ of $3d$ COHA.
Graded components $\H_{{\bf T}_m,\gamma}$ are the equivariant Borel-Moore homology with respect 
to the Cartesian product of the torus ${\bf T}_m:=(\C^{\ast})^m$ and the group $\Gg_\gamma$.
The weights $\varepsilon_i$ of the ${\bf T}_m$-action are chosen 
in such a way that $\prod_i\varepsilon_i=1$ (Calabi-Yau condition). 

%3) There are many papers in which the action of COHA can be computed explicitly, see e.g.
%\cite{Neg,Kor}.
%\end{rmk}

The {\it spherical COHA} is defined  as the subalgebra
$S\H^{(Q,W)}\subset \H^{(Q,W)}$ generated by representations with dimension vectors
$e_i=(0,...,1,...,0), 1\le i\le |I|$. There is an obvious equivariant version of the spherical COHA.

\subsection{$2$-dimensional COHA, algebra ${\bf SH^c}$ and torsion sheaves on $\C^2$}
\label{sec:2.2}

Technical details of the $2d$ COHA will be discussed in \S \ref{subsec:2dCOHA} in detail.
Here we recall few basics facts in order to make a comparison with the $3d$ COHA and
torsion sheaves.

In \cite{SV2} the authors defined an associative algebra ${\bf SH^c}$ over the 
field of 
rational functions in one variable $\C(k)$. It was proved in the loc. cit. that this algebra is 
$\Z$-graded, countably filtered and admits a ``triangular'' decomposition into 
the tensor product (as graded vector spaces) of three subalgebras
$${\bf SH^c}={\bf SH}^{>}\otimes {\bf SH^{c,0}}\otimes {\bf SH}^{<},$$
where the Hilbert-Poincar\'e series of the   
factor ${\bf SH}^{>}$ is given by $\prod_{r>0, l\ge 0}{1\over{1-t^rq^l}}$, 
while for ${\bf SH}^{<}$ 
it is given by  $\prod_{r<0, l\ge 0}{1\over{1-t^rq^l}}$.

%This infinite product almost coincides  with the motivic DT-series of $(Q_3,W_3)$ computed in Proposition 7 from \cite{KS}.
%In the next subsection we consider the spherical COHA $S\H^{(Q_3,W_3)}$  which
%is isomorphic to ${\bf SH}^{>}$ after extension of scalars to the field of rational functions in two variables.
%For the spherical COHA the corresponding motivic DT-series coincides with the above infinite product.

The algebra ${\bf SH}^{>}$ 
is closely related to  the
$2d$ COHA ${\mathcal H}({Coh_{{\bf T}_2,0}(\C^2)})$ 
of the category $Coh_{{\bf T}_2,0}(\C^2)$ of ${\bf 
T}_2:=(\C^{\ast})^2$-equivariant coherent sheaves on 
$\C^2$ which have zero-dimensional support. 
In particular, as   graded vector space the $2d$ COHA is isomorphic to 
$\oplus_{d\ge 0}H_{BM}(Coh_{{\bf T}_2,0}(\C^2,d))$.
Here $Coh_{{\bf T}_2,0}(\C^2,d))$ denotes the stack of equivariant coherent 
sheaves 
with $0$-dimensional support having degree $d$, and 
$H_{BM}$ denotes 
the  Borel-Moore homology of this stack. Since our stack is the 
quotient stack, its (co)homology is defined as the 
equivariant (co)homology of the corresponding algebraic variety.
\footnote{More accurately, our  BM-homology  is defined as the dual to the critical compactly supported cohomology
from \cite{KS}, so they should be  called  critical Borel-Moore homology.}

The algebra ${\bf SH^{c,0}}$ is obtained from the graded commutative algebra $H_{BM}^{GL(\infty)\times {\bf T}_2}(pt)$ of equivariant
Borel-Moore homology  by adding infinitely many
central elements ${\bf c}=(c_1,c_2,...)$ (see loc. cit. for the details).

The associative product on the 
$2d$ COHA   was defined in \cite{SV2} by observing 
that the stack $Coh_{{\bf T}_2,0}(\C^2,d)$ is isomorphic to the stack  
of pairs of commuting $d\times d$ matrices. 
The algebra product is defined in term of correspondences  in the usual 
Nakajima (or Hall algebra)
style.

It is known that the algebra structure on the $2d$ COHA ${\mathcal 
H}({Coh(\C^2)_0})$ 
(no equivariancy conditions)
agrees with the one obtained as a result of the dimensional reduction of  
$\H^{(Q_3,W_3)}$ (see Appendix to \cite{ReSo}). Similar result holds for equivariant and spherical versions.

The spherical  $2d$ 
COHA $S{\mathcal H}({Coh_{{\bf T}_2,0 }(\C^2)})\subset {\mathcal H}({Coh_{{\bf 
T}_2,0}(\C^2)})$
is the subalgebra generated by the equivariant Borel-Moore homology of 
$1\times 1$ matrices. 
Both algebras, $2d$ COHA and spherical $2d$ COHA  are   the  dimensional 
reductions along one
of the loops of the corresponding versions of 3d COHA in the sense of \cite{KS}. 
It was shown in  \cite{SV2} (see Corollary 6.4) that after extension of scalars, 
the spherical Hall algebra  $S{\mathcal H}({Coh(\C^2)_0})$
and the algebra ${\bf SH}^{>}$ become isomorphic.  

The structure of the algebra ${\bf SH^c}$ 
was studied in many papers, both mathematical and physical
(to mention just a few:  \cite{AS}, \cite{BouMaZh}, \cite{GabGoLPe, P1}).
Furthermore the $\W$-algebra $\W_r$ is closely related (see \cite{SV2}) 
to the algebra ${\bf SH^c}$. From the point of view of the above discussion,
it becomes clear that $\W_r$ is the ``dimensional reduction" of some
``$3$-dimensional VOA". This agrees with the previously made remark that $\W_r\simeq\W_{r, r_2=0, r_3=0}$.

\subsection{$3d$ COHA and torsion sheaves on $\C^3$}

The stack of finite-dimensional representations of the algebra $\C[x_1,x_2]$
 contains the 
Hilbert scheme $Hilb(\C^2)=\sqcup_{n\ge 0}Hilb_n(\C^2)$ consisting of cyclic 
finite-dimensional representations of $\C[x_1,x_2]$.
As we know since Nakajima's work  this Hilbert scheme is isomorphic to the 
moduli space of stable 
framed representations of a certain quiver, and there is a similar description 
for the rank $r$ torsion free sheaves $F$ on ${\bf P}^2$ which are framed at the 
line ${\bf P}^1_\infty$ at infinity via an 
isomorphism $F_{|{\bf P}^1_\infty}\simeq \OO_{{\bf P}^1}^r$ (see \cite{Nak1}).

It is known (see e.g. \cite{SchVas2}) that the   $2d$ 
COHA ${\mathcal H}({Coh_{{\bf T}_2, 0}(\C^2)})$ is isomorphic to the positive part 
of 
the affine Yangian $\affY$. More conceptual approach to the 
quiver Yangians of \cite{MO}
based on $3d$ COHA
was proposed in \cite{D}. 

Furthermore, for each pair $(r,n)\in \Z_{\ge 1}\times \Z$ one can define the 
moduli space $\mathfrak{M}_{r}(n)$ of rank $r$ 
framed torsion-free sheaves 
on ${\bf P}^2$ which have the second Chern class $n$. 
Let $V_r=\oplus_n H_{BM}^{(\C^{\ast})^2\times GL_r}(\mathfrak{M}_{r}(n))$ 
be the equivariant Borel-Moore homology 
of the disjoint union of such moduli spaces. 
Then the  $2d$ COHA ${\mathcal H}({Coh_{{\bf T}_2, 0}(\C^2)})$
acts faithfully on $V_r$ by correspondences. 
This gives a generalization of the classical results of Nakajima (see \cite{Nak1}).

One of the  aims of this paper is to explain the $3$-dimensional avatar of this observation.
Although we consider only the case of the quiver with potential $(Q_3,W_3)$ there is little doubt that
the results hold for any quiver with potential coming from a dimer model (hence for any toric Calabi-Yau $3$-fold).

Thus we define a spherical $3d$ COHA for a general quiver with potential $(Q,W)$ 
in the obvious way: it is  the subalgebra $S{\mathcal H}^{Q,W}\subset 
{\mathcal H}^{Q,W}$ generated by $1$-dimensional representations. The 
equivariant version is defined similarly. 

The $3d$ COHA is  $\H^{(Q_3,W_3)}={\mathcal 
H}({Coh(\C^3)_0})$, i.e. it is the COHA of the category of 
$0$-dimensional 
torsion sheaves on $\C^3$. Then we have the corresponding spherical $3d$ COHA 
$S{\mathcal H}({Coh(\C^3)_0})$ as well as the equivariant versions $S{\mathcal 
 H}({Coh_{{\bf T}_3, 0}(\C^3)})\subset
\mathcal{H}({Coh_{{\bf T}_3, 0}(\C^3)})$ which are COHA and spherical COHA of 
the category of ${\bf T}_3=(\C^{\ast})^3$-equivariant $0$-dimensional torsion 
sheaves on $\C^3$.

In general we will denote by $Coh(X)_0$ the category 
of  sheaves on a toric Calabi-Yau $3$-fold which have a 
$0$-dimensional support. We will be mostly interested in the case $X=\C^3$.

\begin{prp}  The category $Coh(\C^3)_0$ is equivalent to the category 
of finite-dimensional modules over the 
Jacobi algebra $\C Q_3/\langle \partial W_3\rangle$, 
where $\C Q_3$ is the path algebra of the quiver $Q_3$ (see 
Example \ref{exa:QW}).
 \end{prp}
 
{\it Proof.} It follows from the observation that the corresponding Jacobi 
algebra is isomorphic to 
the algebra of polynomials $\C[B_1,B_2, B_3]$. $\blacksquare$

Then  COHA 
${\mathcal H}({Coh_{{\bf 
T}_3, 0}(\C^3)})=\oplus_{d\ge 0}H_{BM}^{{\bf T}_3}
(Coh_{{\bf  T}_3, 0}(\C^3,d))$, where $Coh_{{\bf 
T}_3, 0}(\C^3,d))$  is the stack of 
degree $d$ equivariant sheaves 
on $\C^3$ which have $0$-dimensional support. This COHA is isomorphic to the 
equivariant version of COHA 
$ {\mathcal H}^{(Q_3,W_3)}$.

\section{Framed quiver and its stable representations}

\label{sec:framed_quiver}\subsection{Framed quiver}
%Let $(Q_3, W_3)$ be the quiver with potential considered in Example \ref{exa:QW}. That is, $Q_3$ is the quiver with one vertex labelled by $0$, and three loops $B_1, B_2, B_3$. The potential is $W_3=B_3[B_1,B_2]$ (plus cyclic permutations). 

In this subsection we are going to define the framed quiver with potential  
$(Q_3^{fr}, W^{fr}_3)$. The framed quiver
$Q_3^{fr}$ is obtained by adding to $Q_3$ three new vertices 
(framing vertices) with three 
pairs of opposite arrows.

More precisely, the quiver $Q_3^{fr}$ is defined such as follows:
\begin{itemize}
\item The set of vertices is \{0, 1, 2, 3\}, where the vertices $1,2,3$ are 
framing vertices.
\item There are two types of arrows. 
\begin{enumerate}
\item
We have three  loops $B_1, B_2, B_3$ at the vertex $0$.  Hence $(0, 
B_1,B_2, B_3)$ is the quiver $Q_3$
considered previously.
\item
We also have ``framing" arrows
$I_{12}: 3\to 0$, $J_{12}: 0\to 3$; 
$I_{13}: 2\to 0$, $J_{13}: 0\to 2$;
$I_{23}: 1\to 0$, $J_{23}: 0\to 1$. 
\end{enumerate}
\end{itemize}

We define the framing potential $W_3^{fr}$ by the formula
\begin{align*}
W^{fr}_3&:=B_3([B_1, B_2]+I_{12}J_{12})+B_2([B_1, B_3]+I_{13}J_{13})+B_1([B_2, 
B_3]+I_{23}J_{23})\\
&=
W_3+B_1I_{23}J_{23}+B_2I_{13}J_{13}+B_3I_{12}J_{12}.
\end{align*}

%\begin{tikzpicture}

%\tikzset{every loop/.style={min distance=30mm,looseness=10}}

%\draw[step=1cm,gray,very thin] (0,0) grid (8,8);
%\draw[fill] (8,4) circle [radius=0.1];
%\node [above] (inf1) at (8,4) {};
%\node [below] (inf2) at (8,4) {};
%\node [right] (inf1r) at (8,4) {};
%\node [left] (inf2l) at (8,4) {};

%\draw (3,4) ellipse (1 cm and 0.1 cm);
%\draw (4,5) ellipse (0.1 cm and 1 cm);
%\draw (5,4) ellipse (1 cm and 0.1 cm);

%\draw[fill](6,6) circle [radius=0.1];
%\draw[fill](10,6) circle [radius=0.1];
%\draw[fill](8,2) circle [radius=0.1];

%\node [above] (inf3) at (6,6) {};
%\node [below] (inf4) at (6,6) {};

%\node [above] (inf5) at (10,6) {};
%\node [below] (inf6) at (10,6) {};

%\node [right] (inf7) at (8,2) {};
%\node [left] (inf8) at (8,2) {};

%\draw [thick, ->] (inf1)-- node[above]{$j_1$}(inf3);
%\draw [thick, <-] (inf2)-- node[below]{$i_1$}(inf4);

%\draw [thick, ->] (inf1)-- node[above]{$j_2$}(inf5);
%\draw [thick, <-] (inf2)-- node[below]{$i_2$}(inf6);

%\draw [thick, ->] (inf1r)-- node[right]{$j_3$}(inf7);
%\draw [thick, <-] (inf2l)-- node[left]{$i_3$}(inf8);

%\draw [thick, ->] (inf1)  edge[loop] node[above]{$B_3$}(inf1);
%\draw [thick, ->] (inf2l)  edge[loop left] node[left]{$B_2$}(inf2l);
%\draw [thick, ->] (inf1r)  edge[loop right] node[right]{$B_1$}(inf1r);
%\end{tikzpicture}
%{\bf Figure: Framed  quiver $Q_3^{fr}$ }
%\vspace{3mm}

In order to define a representation of $(Q_3^{fr}, W_3^{fr})$ let us fix a 
dimension vector $(n, r_1, r_2, r_3)$ at the vertices $\{0, 1, 2, 3\}$ 
respectively. Then for a representation of $Q_3^{fr}$ of this dimension we have 
the linear maps
\begin{align*}
&B_i: \C^n\to  \C^n, i\in \{1, 2, 3\}  \\
&I_{12}: \C^{r_3}\to \C^n, \,\ J_{12}: \C^n\to \C^{r_3};\\ 
&I_{13}: \C^{r_2}\to \C^n,\,\  J_{13}: \C^n\to \C^{r_2};\\ 
&I_{23}: \C^{r_1}\to \C^n,\,\  J_{23}: \C^n\to \C^{r_1}. 
\end{align*}

%\begin{center}
 % \centering
 %     \includegraphics[width=0.29\textwidth]{quiver.pdf}
 % \caption{}
%\end{center}

We will use the shorthand notation $\vec{r}:=(r_1, r_2, r_3)$, and $\C^{\vec{r}}:=\C^{r_1}\oplus 
\C^{r_2}\oplus \C^{r_3}$. The space of representations of $Q_3^{fr}$ with the dimension vector
$(n,\vec{r})$ is
\begin{align*}
\calM_{\vec{r}}(n):=&\{ (B_i, I_{ab}, J_{ab} )\mid i\in \underline{3}, 
ab\in  \overline{3}\}\\
=&\Hom(\C^n, \C^n)^3
\oplus \Hom(\C^n, \C^{\vec{r}})\oplus \Hom(\C^{\vec{r}}, \C^n), 
\end{align*}
where $\underline{3}:=\{1, 2, 3\}$, $\overline{3}:=\big\{\{1, 2\}, \{1, 3\}, \{2, 3\}\big\}$. 

The groups $\GL_n(\C)$ and $\GL_{\vec{r}}:=\GL_{r_1}(\C)\times \GL_{r_2}(\C) \times \GL_{r_3}(\C)$ act on $\calM_{\vec{r}}(n)$ by conjugation. 
Let ${{\bf T}_2}$ be the two dimensional torus. We think ${{\bf T}_2}$ as a subtorus in $(\C^{*})^3$ given by the equation $t_1t_2t_3=1$. 
Then, ${{\bf T}_2}$ acts on $\calM_{\vec{r}}(n)$  by
\begin{equation*}
(t_1, t_2, t_3)\cdot \{B_i, I_{ab}, J_{ab}\}
=(t_1B_1, t_2B_2, t_3B_3,  I_{12}, t_1t_2 J_{12}, 
I_{13}, t_1t_3 J_{13},  I_{23}, t_2t_3J_{23}),
\end{equation*}
where $t_1t_2t_3=1$, and $\{B_i, I_{ab}, J_{ab}\}\in \calM_{\vec{r}}(n)$. 
Under this ${{\bf T}_2}$-action, the potential $W_3^{fr}$ 
is invariant.

Let us introduce the  ``moment maps" $\mu_{ab}: \calM_{\vec{r}}(n) \to 
\mathfrak{gl}_n$, for $ab\in \overline{3}$ given by
\begin{align*}
&\mu_{12}=[B_1, B_2]+I_{12}J_{12}, \,\ \mu_{13}=[B_1, B_3]+I_{13}J_{13},\,\ \mu_{23}=[B_2, B_3]+I_{23}J_{23}.
\end{align*}
Then  $W^{fr}_3=B_3\mu_{12}+B_2\mu_{13}+B_1\mu_{23}$.
\begin{lmm}
 For  $1\leq a, b, c \leq 3$, and $a, b, c$ are distinct, 
we have 
\begin{align*}
&\frac{\partial W^{fr}_3}{ \partial B_a}=\mu_{bc}, \,\ 
\frac{\partial W^{fr}_3}{ \partial I_{ab}}=J_{ab} B_c, \,\
\frac{\partial W^{fr}_3}{ \partial J_{ab}}=B_c I_{ab}. 
\end{align*}
\end{lmm}
{\it Proof.}
We have
\begin{align*}
\frac{\partial W^{fr}_3}{ \partial B_1}
=&\frac{\partial (B_3\mu_{12}+B_2\mu_{13}+B_1\mu_{23})}{ \partial B_1}
=[B_2, B_3] +[B_3, B_2]+\mu_{23}=\mu_{23}. 
\end{align*}
The other identities follow from a straightforward calculation. 
The lemma is proven. $\blacksquare$

\subsection{Stable representations of the framed quiver}

We say (cf. \cite{N}) that a representation $(B_i, I_{ab}, J_{ab})_{i\in \underline{3}, 
ab\in  \overline{3}}$ belonging to $\calM_{\vec{r}}(n)$ is stable if the following
condition is satisfied:

\begin{equation}\label{eq:stabN}
 \C\angl{B_1, B_2, B_3} I_{12}(\C^{{r_3}})+
\C\angl{B_1, B_2, B_3} I_{13}(\C^{r_2})+
\C\angl{B_1, B_2, B_3} I_{23}(\C^{r_1})
=\C^{n}, 
\end{equation}
where $\C\angl{B_1, B_2, B_3}$ is the ring of non-commutative polynomials in the variables $B_1, B_2, B_3$.  
In other words, a representation is stable, if the non-commutative polynomials in $B_1, B_2, B_3$ applied to the image of $I_{12}, I_{13}, I_{23}$ generate  $\C^n$. 
\begin{exa}
If $\vec{r}=(0, 0, 1)$, the stability condition \eqref{eq:stabN} becomes 
$\C\angl{B_1, B_2, B_3} I_{12}(\C^{1})=\C^{n}$, which means $I_{12}$ is a cyclic vector. 
\end{exa}

One can also impose a stronger stability condition by saying that the representation is stable if the following condition is 
satisfied:

\begin{equation}\label{eq:stabD}
\C\angl{B_1, B_2} I_{12}(\C^{r_3})
+\C\angl{B_1, B_3} I_{13}(\C^{r_2})
+ \C\angl{ B_2, B_3} I_{23}(\C^{r_1})=\C^n. 
\end{equation}
Note that \eqref{eq:stabD} implies \eqref{eq:stabN}. 

\begin{prp}\label{prop:critical_locus}
Assume $(B_i, I_{ab}, J_{ab})_{i\in \underline{3}, ab\in \overline{3}}$ is in 
the critical locus of $W^{fr}_3$. That is, $(B_i, I_{ab}, J_{ab})$ satisfies the following equations
\begin{align*}
&\mu_{12}=[B_1, B_2]+I_{12}J_{12}=0,\,\, \mu_{13}=[B_1, B_3]+I_{13}J_{13}=0,\\
&\mu_{23}=[B_2, B_3]+I_{23}J_{23}=0,\\
&J_{12} B_3=0, \,\ J_{13} B_2=0, \,\ J_{23} B_1=0,\\
&B_3 I_{12}=0, \,\  B_2 I_{13}=0, \,\ B_1 I_{23}=0. 
\end{align*}
Then, the condition \eqref{eq:stabN} is equivalent to \eqref{eq:stabD}. 
\end{prp}
{\it Proof.}
We start with a representation $(B_i, I_{ab}, J_{ab})\in \calM_{\vec{r}}(n)$, which is stable under condition \eqref{eq:stabN}: 
\[
\C^n=\C\angl{B_1, B_2, B_3}I_{12} (\C^{r_3})+ \C\angl{B_1, B_2, 
B_3}I_{13}(\C^{r_2})+ \C\angl{B_1, B_2, B_3}I_{23}(\C^{r_1}). 
\] 
%Assume it satisfies the equations $
%B_3 I_{12}=0,  B_1 I_{23}=0,  B_2 I_{13}=0. 
%$
% As $B_3 I_{12}=0$, we have

Using the assumptions we have:
\begin{align*}
&\C\angl{B_1, B_2, B_3}I_{12} (\C^{r_3})\\
=&\C\angl{B_1, B_2}I_{12} (\C^{r_3})
+ \C\angl{B_1, B_2, B_3} [B_1,B_3] \C\angl{B_1, B_2, B_3} I_{12}(\C^{r_3})\\
&+ \C\angl{B_1, B_2, B_3} [B_2,B_3] \C\angl{B_1, B_2, B_3} I_{12}(\C^{r_3})\\
=&\C\angl{B_1, B_2}I_{12} (\C^{r_3})
+ \C\angl{B_1, B_2, B_3} I_{13}J_{13} \C\angl{B_1, B_2, B_3} I_{12}(\C^{r_3})\\
&+ \C\angl{B_1, B_2, B_3} I_{23}J_{23} \C\angl{B_1, B_2, B_3} I_{12}(\C^{r_3})\\
\subset 
&\C\angl{B_1, B_2}I_{12} (\C^{r_3})
+ \C\angl{B_1, B_2, B_3} I_{13}(\C^{r_2})
+ \C\angl{B_1, B_2, B_3} I_{23}(\C^{r_1}). 
\end{align*} 
As a consequence,  the condition \eqref{eq:stabN} can be simplified to
 \[
 \C^n=\C\angl{B_1, B_2}I_{12} (\C^{r_3})+ \C\angl{B_1, B_2, 
B_3}I_{13}(\C^{r_2})+ \C\angl{B_1, B_2, B_3}I_{23}(\C^{r_1}).
 \] 
Moreover, for any homogenous $f(B_1, B_2, B_3)  \in \C\angl{B_1, B_2, B_3}$, such that \begin{eqnarray*}
&&f(B_1, B_2, B_3) I_{12}(\C^{r_3})\\
&=&f_1(B_1, B_2) I_{12}(\C^{r_3})
+g_1(B_1, B_2, B_3) I_{13}(\C^{r_2})+h_1(B_1, B_2, B_3) I_{23}(\C^{r_1}), 
\end{eqnarray*}
we could assume $
\deg_{B_3}(f) > \deg_{B_3}(g_1), \,\  \deg_{B_3}(f) > \deg_{B_3}(h_1)$. 
Similarly, we have
\begin{eqnarray*}
&&g_1(B_1, B_2, B_3) I_{13}(\C^{r_2})\\
&=&f_2(B_1, B_2, B_3) I_{12}(\C^{r_3})
+g_2(B_1, B_3)  I_{13}(\C^{r_2})
+h_2(B_1, B_2, B_3) I_{23}(\C^{r_1}), 
\end{eqnarray*}
with $\deg_{B_3}(f) >\deg_{B_3}(f_2)$. 

One can thus iteratively get  
 \[
 \C^n=\C\angl{B_1, B_2}I_{12}(\C^{r_3})+ \C\angl{B_1,  B_3}I_{13}(\C^{r_2})+ 
\C\angl{B_1, B_2, B_3}I_{23}(\C^{r_1}).
 \] 
Repeatedly using the equations of the critical locus of $W^{fr}_3$, we can also get rid of $B_1$ from the last factor. This shows the two 
conditions are equivalent. 
 The proposition is proven. $\blacksquare$

\begin{rmk}
%1) Note that condition \eqref{eq:stabN} and \eqref{eq:stabD} are different if we 
%do not impose the conditions $\{B_a I_{bc}=0\}$. 

%2) 
Alternative proof of Proposition \ref{prop:critical_locus} follow from  the 
non-holomorphic generalization of ADHM equations 
proposed by Nekrasov. 
The corresponding discussion of stability can be found in \cite{N}, Section 
3.4. 
Roughly, the idea is to study critical points of the function 
$f=Tr(\sum_i M_iM^{\ast}_i)$, where $M_i$ are matrices in the 
representations of the left hand sides of the relations from Proposition \ref{prop:critical_locus}
and $M_i^{\ast}$ are their Hermitian conjugates. One shows that if at a 
critical point of $f$ the stability assumption \eqref{eq:stabN} holds, then it is (up to the 
gauge group action) is a point of absolute minimum of $f$. On the one hand at 
such points all $M_i=0$, i.e. the equations from the Proposition are satisfied.
On the other hand they are stable representations, since they are points of 
absolute minimum. The details can be found in the loc.cit. in a bigger 
generality.
\end{rmk}

Let
\begin{equation}
\calM_{\vec{r}}(n)^{st}
=\Big\{( B_i, I_{ab}, J_{ab})\in  \calM_{\vec{r}}(n) \mid \sum_{ab\in\overline{3}}\C\langle B_1,B_2,B_3\rangle I_{ab}(\C^{r_{ab}})=\C^n\Big\}
\end{equation}
be the stable locus of $\calM_{\vec{r}}(n)$. It consists of 
the representations of $(Q^{fr}_3, W^{fr}_3)$ of dimension vector $(n, \vec{r})$ that satisfy the equation \eqref{eq:stabN}. 
The quotient $\calM_{\vec{r}}(n)^{st}/\GL_n$ is denoted by $\fM_{\vec{r}}(n)$. We call this moduli space of solutions to the equation \eqref{eq:stabN} the moduli space of spiked instantons.

%\textcolor{blue}{NEXT PARAGRAPH SHOULD BE REWRITTEN: WE ARE TALKING ABOUT REPS OF NON-FRAMED 
%QUIVER IN FRAMED MODULI SPACE}
\begin{rmk}\label{rmk_framestable_general}
The above representations can be thought of as generalizations of stable framed 
representations of the quiver with potential
discussed in \cite{S}. In fact in the loc. cit.  the general notion of stable 
framed object in the framework of triangulated categories endowed with stability 
structure is proposed. 

The difference of our case discussion with the one in \cite{S} is that here we would like to consider representations of COHA
for the quiver $(Q_3,W_3)$ in the cohomology of the moduli space of stable representations
of the quiver with potential
$(Q_3^{fr}, W_3^{fr})$. Hence the relations for the stable framed representations
are given by another potential. 

%The rest of the property to be stable framed is the same as in the loc. cit. Namely,
%we consider representations of the Jacobi algebra $\C\langle B_1, B_2, 
%B_3\rangle/\langle\partial W_3^{fr}\rangle$ in the $n$-dimensional
%vector space $V\simeq \C^n$ such that
%there is no 
%proper subrepresentation which contains images of all maps $I_{ab}$.
%In other words 
%\[
%\sum_{ab\in\overline{3}}\C\langle B_1,B_2,B_3\rangle I_{ab}(E_{ab})=V,\]
%where $E_{ab}\simeq \C^{r_{ab}}$ are framing vector spaces.

\end{rmk}

%\begin{rmk}
%We have at out disposal also the following stability condition
%\[
%\C\angl{B_1, B_2} I_{12}(\C^{r_3})=\C^n, \text{or}\,\ 
%\C\angl{B_1, B_3} I_{13}(\C^{r_2})=\C^n, \text{or}\,\ 
%\C\angl{B_2, B_3} I_{23}(\C^{r_1})=\C^n. 
%\]
%It corresponds to $(a, b, c)=(n, b, c)$, or $(a, b, c)=(a, n, c)$, or  $(a, b, 
%c)=(a, b, n)$. 
%
%This condition is strictly stronger than the one we used above. To distinguish 
%these two stability conditions, the semi-stable (STABLE?) locus with respect 
%to this stability condition is denoted by $\Rep(Q_3^{fr},\vec{r})^{s}$. 
%\end{rmk}

%\begin{prp} Stable framed representation of $(Q_3^{fr}, W^{fr}_3)$ is 
%given by a finite-dimensional representation of the quiver $Q_3^{fr}$ 
%with relations
%$$[B_1,B_2]+i_3j_3=0, [B_2,B_3]+i_1j_1=0, [B_1,B_3]+i_2j_2=0,$$
%$$j_kB_k=0, B_ki_k=0, 1\le k\le 3,$$
%subject to the following stability condition: there is no 
%proper subrepresentation which contains images of all maps $i_1,i_2,i_3$.
%In other words $\sum_{1\le k\le 3}\C[B_1,B_2,B_3]i_k(W_k)=V$,
%where $W_k\simeq \C^{r_k}, 1\le k\le 3$ are framing vector spaces,
%and $V\simeq \C^n$ is a vector space associated to the only vertex of the 
%quiver $Q_3$.

%\end{prp}

%{\it Proof.} The proof is essentially a generalization of the one in the case 
%when only $B_1,B_2, i_3,j_3$ are non-trivial (see \cite{Nak1}, Chap. 2). 
%$\blacksquare$

Similarly to \cite{S} one can prove the following result.

\begin{prp}\label{prp:autstable} 
\begin{enumerate}
	\item The group of automorphisms of a stable framed representation is 
	trivial.
	\item The isomorphism classes of stable framed representations naturally form an algebraic variety, 
	called the moduli space of stable framed representations.
	\item The variety $\fM_{\vec{r}}(n)$ is smooth.
\end{enumerate}\end{prp}
For completeness, we prove the last assertion. 
Let $(B_i,I_{jk},J_{jk})$ be the representative of an isomorphism class of representations in $\fM_{\vec{r}}(n)$.  We need to prove that its automorphism group is trivial (cf. Proposition \ref{prp:autstable}). 
We have the following 3 vector subspaces on $\C^n$, 
\[
V_1=\{f(B_2,B_3)I_{23}(\C^{r_1})\mid f\hbox{ a non-commutative polynomial in 2 variables}\}
\] and similarly $V_2$ and $V_3$. The point $(B_i,I_{jk},J_{jk})$ being stable means that $\C^n=V_1+V_2+V_3$.

\begin{rmk}
Recall the geometric interpretation of the moduli of stable representations of $(Q_3^{fr},W_3^{fr})$ in the case
when $r_2=r_3=0$, and hence $i_k=0, j_k=0, k=2,3$, and $r_1=1$. It is  the moduli space
of torsion sheaves on $\C^3$ such that their restriction to any plane $x_3=c$ is a cyclic module over the
algebra $\C[B_1,B_2,B_3=c]=\C[B_1,B_2]$. Such a module is the same as a point of the Hilbert scheme $Hilb(\C^2)$,
but it is better to think that $\C^2$ is placed in $\C^3$ as a plane  $x_3=c$.

Generally, if  $r_1=r\ge 1$ we obtain the moduli space of rank $r$ instantons on $\C{\bf P}^2$, equivalently, torsion-free
rank $r$ sheaves on $\C{\bf P}^2$ endowed with an isomorphism with $\OO^{r}$ on the line $l_{\infty}\simeq \C{\bf P}^1=
\C{\bf P}^2-\C^2$ at
infinity.

In the subsequent paper we plan  to study an analogous geometric description for an arbitrary triple $(r_1,r_2,r_3)$.

\end{rmk}

\section{Reminder on COHA for $(Q_3,W_3)$  and main result}
\label{sec:Main result}

\subsection{Reminder on the critical cohomology}\label{subsec:reminder_crit_coh}
Let $D^b(X)$ be the bounded derived category of constructible sheaves of $\bbQ$-vector spaces on an algebraic variety $X$ over a field of characteristic zero, and $\mathbb{D}_X$ be the Verdier duality functor for $D^b(X)$. In particular, $\mathbb{D}_{\pt}$ is the vector space dual. Denote by $H_c(X)^{\vee}$ the Verdier dual of the compactly supported cohomology $H_c(X)$ of $X$. 
Let $t: X\to \pt$ be the structure map. Then we have 
$H_{c}^*(X)^\vee= \mathbb{D}_{\pt}(t_{!} \Q_{X})$. Note that this is the same as $t_*\mathbb{D}_X\Q_{X}$, which is the  Borel-Moore homology of $X$ in the usual sense.

If $X$ carries a $G$-action, where $G$ is an algebraic group, we denote by $H_{c, G}^*(X)^\vee$ the Verdier dual  to the corresponding equivariant compactly supported cohomology of $X$. More generally, we consider cohomology valued in a weakly equivariant sheaf \cite{GKM}, in order to apply the general machineries of \cite{GKM,Bra}. However, in the situation we are interested in, $A$ will be taken to be the constant sheaf or its images under functors well-define for weakly equivariant sheave, hence will have a natural structure of a weakly equivariant sheaf. For any weakly equivariant complex of constructible sheaves $A$ on $X$, we define $H_{c,G}(X,A)^\vee:=\mathbb{D}_{\pt}p_!A$. This is a module over $H^*_G(\pt)$. The latter is the ring of functions on $\fg^*:=Lie(G)^\ast$ if $G$ is abelian, and the ring of conjugation invariant functions on $\fg^*$ in general. 

We will be considering the critical cohomology in the following setup. 
Let $G$ be a complex linear algebraic group.
Suppose $X$ is a smooth complex algebraic variety with a $G$-action, endowed with a $G$-invariant regular function $f$. 
As in \cite[Section 2.4]{D}, we assume every $x\in X$ is contained in a $G$-invariant open affine neighborhood. 
We have the functor of vanishing cycles $\varphi_f$ (see e. g. \cite[Section 7.2]{KS}), 
which applied to any weakly $G$-equivariant complex of sheaves on $X$, resulting in a weakly equivariant complex of sheaves. 
%(??? DO WE CONSIDER COMPLEXES OF EQUIVARIANT SHEAVES?). 
We are primarily interested in the complex $\varphi_f\Q_X$, which is supported on the critical locus of $f$. The {\it critical cohomology} of $X$ is defined to be $H_{c,G}(X,\varphi_f\Q_X)^\vee$. Abusing the notation we will sometimes denote it by $H_{c,G}(\Crit(f),\varphi_f\Q_X)^\vee$.

The critical cohomology enjoys natural functoriality (see e.g. \cite{D}, Section 2.7). In particular, assume that we are given  an algebraic  $G$-variety $Y$ which is endowed with a $G$-invariant function $g$. Let $F:X\to Y$ be a $G$-equivariant map of algebraic $G$-varieties.  Let $f=gF$.  Then  there is a pullback map $F^*:H_{c,G}(Y,\varphi_g\Q_Y)^\vee\to H_{c,G}(X,\varphi_f\Q_X)^\vee$. If furthermore $F$ is a
proper map then
there is a pushforward map $F_*:H_{c,G}(X,\varphi_f\Q_X)^\vee\to H_{c,G}(Y,\varphi_g\Q_Y)^\vee$. If $F: X\to Y$ is an affine bundle, then the pullback $F^*$ induces an isomorphism.
It follows that if $\pi:\calV\to X$ is a $G$-equivariant vector bundle, then for every $G$-invariant regular function $f$ on $X$, there is a well-defined {\it Euler class}, defined as an operator $Eu(\calV):=Eu_f(\calV): H_{c,G}(X,\varphi_f\Q_X)^\vee\to H_{c,G}(X,\varphi_f\Q_X)^\vee$. More precisely, $Eu(\calV)=s^*s_{*}$ where
\begin{align*}
s_{*}: H_{c, G}(X, \varphi_f \Q_X)^{\vee} \to H_{c, G}(\calV, \varphi_{f\pi}\Q_{\calV})^{\vee} 
\end{align*}
and 
\begin{align*}
s^{*}:  H_{c, G}(\calV, \varphi_{f\pi}\Q_{\calV})^{\vee} \to H_{c, G}(X, \varphi_f \Q_X)^{\vee} 
\end{align*}
where  $s: X \inj \calV$ is the zero section of $\calV$, and the composition $f\pi: \calV\to \mathbb{C}$ is considered as a $G$-invariant regular function on the total space of $\calV$. 
This operator has an alternative description given in \cite{D}, Proposition~2.16, which we will not use in this paper. 

Note that taking $f=0$, we obtain an operator $H_{c,G}(X,\Q_X)^\vee\to H_{c,G}(X,\Q_X)^\vee$. This works even in the case when $X$ is not necessarily smooth. We keep the notation $Eu(\calV)$ for this operator. In the case when $X$ is smooth, there is a fundamental class $1$ in $H_{c,G}(X,\Q_X)^\vee$, which is the multiplicative unity for the ring structure.  Let $Eu(\calV)\cdot 1 \in H_{c,G}(X,\Q_X)^\vee$ be the element obtained by applying $Eu(\calV)$ to this identity element. Then, $Eu(\calV)\cdot 1$ is the usual Euler class in $H_{c,G}(X,\Q_X)^\vee$ in the classical sense.
%(WHAT THE NOTATION $\cdot$ MEANS HERE?) 

In general, the  Euler class operator satisfies the  Whitney sum formula. That is, if $\calV$ and $\calW$ are two equivariant vector bundles on $X$, and $f$ is a $G$-invariant regular function on $X$, then for the corresponding Euler operators we have $Eu_f(\calV)\circ Eu_f(\calW)=Eu_f(\calW)\circ Eu_f(\calV)=Eu_f(\calW\oplus\calV): H_{c,G}(X,\varphi_f \Q_X)^\vee\to H_{c,G}(X,\varphi_f\Q_X)^\vee$.

Let  
$\pi: \calL \to X$ be a $G$-equivariant line bundle on $X$. Then we will sometimes call the Euler class operator associated with $\pi$ by {\it first Chern class operator}, or simply {\it first Chern class}.
By the splitting principle, we can define in the standard way the Chern roots of any equivariant vector bundle. Any symmetric polynomials in the Chern roots of $\calV$ are well-defined operators on $H_{c,G}(X,\varphi_f\Q_X)^\vee$.

Applying to the special case when $G=\GL_n(\C)$ and $\calV\to X$ is the trivial vector bundle $X\times\C^n\to X$ with the natural diagonal action, we see that
the equivariant Chern roots are the same as the equivariant parameters with respect to the natural action of the maximal torus. More precisely, any symmetric polynomial in the equivariant parameters is an element in $H^*_G(\pt)$. Its natural action on $H_{c,G}(X,\varphi_f\Q_X)^\vee$ coincides with the action
of the corresponding symmetric polynomial in Chern roots.

Let $\calV \to X$ be a $G$-equivariant vector bundle of rank $n$ and $f$ a $G$-invariant regular function on $X$. 
Let  $z$  be the standard coordinate of the Lie algebra $Lie(\C^*)=\C$, and let $q$ be the natural 1-dimensional representation of $\C^*$. Let $\C^*$ act trivially on $X$. Then, $H_{c, G\times\C^*}(X,  \varphi_f)^{\vee}\cong H_{c, G}(X,  \varphi_f)^{\vee}[z]$ with $z$ acting trivially on $H_G(X,  \varphi_f)^{\vee}$. The Euler class \[Eu(\calV\otimes q)\] is a well-defined operator on $H_{c, G}(X,  \varphi_f)^{\vee}[z]$. 
Let us expand the operator $Eu(\calV\otimes q)$ in the powers of $z$ at $z=0$. Then 
the coefficient of $z^n$ for any $n\in \mathbb N$ is a symmetric polynomial in the Chern roots of $\calV$, hence a well-defined operator on $H_{c, G}(X,  \varphi_f)^{\vee}$. Let $\lambda_1, \cdots, \lambda_n  \in H_G(X, \varphi_f)^{\vee}$ be the Chern roots of $\calV$. 
Then we have
\[
Eu(\calV\otimes q)=\prod_{i=1}^n (\lambda_i+z).
\] 
In particular, $Eu(\calV\otimes q)$ defines an element in $\End(H_{c, G}(X,  \varphi_f)^{\vee})[z]$.
Note that as a polynomial in $z$, the highest degree of $Eu(\calV\otimes q)$ is $z^{\rank(\calV)}$.

Following the usual convention, we define the total Chern polynomial of $\calV$, denoted by $\lambda_{1/z}$, to be the following polynomial of operators 
on $H_{c,G}(X,\varphi_f\Q_X)^\vee[z^{-1}]$ in the variable $1/z$
\[
\lambda_{1/z}(\calV)=\frac{Eu(\calV\otimes q)}{z^{\rank(\calV)}}=\prod_{i=1}^n (1+\lambda_i /z).
\] 
The leading term of $\lambda_{1/z}(\calV)$ is $1$, and the coefficients of $z^{-j}$ commute with each other. Thus, it has the inverse, which is also given by a power series in $1/z$ with coefficients being operators on $H_{c,G}(X,\varphi_f\Q_X)^\vee$. 
By the Whitney sum formula, the total Chern polynomial is well-defined for any equivariant virtual vector bundles, and hence we have a family of mutually commuting operators defined by the map
\[
\lambda_{1/z}: K_G(X)\to \End(H_{c,G}(X,\varphi_f\Q_X)^\vee)[[z^{-1}]].\]
Here $K_G(X)$ is the Grothendieck group of the category of $G$-equivariant vector bundles
on $X$.

The following fact is standard and will be used.
\begin{lmm}\label{lem:restriction to open}
Assume $X$ is a $G$-variety, with $f: X\to \C$ a $G$-invariant function on $X$. Let 
$\varphi_f$ be the vanishing cycle sheaf for $f$. 
Let $j: U\subset X$ be an open embedding which is stable under the $G$ action. 
Denote by $f|_{U}: U\to \C$ be the restriction of $f$ on $U$. 
Then we have the natural homomorphism of graded vector spaces:
\[
j^*: H_{c, G}^*(X,  \varphi_{f})^\vee \to H_{c, G}^*(U,  \varphi_{f|_U})^\vee. 
\]
\end{lmm}
{\it Proof.}
By definition, we have 
$H_{c}^*(X,  \varphi_{f}):=p_{X!} \varphi_{f}[-1](\bbQ_X)$, where $p_{X}: X\to 
\pt$ is the projection to a point. 
As $j$ is an open embedding, we have $j^{!}=j^*$. This gives the following isomorphism 
\begin{align*}
H_{c}^*(U,  \varphi_{f|_{U}})
:=p_{U!} \varphi_{f|_{U}}[-1](\bbQ_U)
=p_{X!}j_{!}   \varphi_{f \circ j}[-1] j^*\bbQ_X
\cong &p_{X!}j_{!} j^!\varphi_{f}[-1] \bbQ_X, 
\end{align*}
where the last isomorphism uses the fact $\varphi_{f \circ j} j^* \cong j^* 
\varphi_{f} $ (See \cite[Page 10]{D}). 
The natural map $j_{!} j^!\to \id$ induces $H_{c}^*(U,  \varphi_{f|_{U}})\to 
H_{c}^*(X,  \varphi_{f})$. 
The desired map is obtained by taking the Verdier duality. The argument is extended to the equivariant setting the same way as \cite[\S~2.4]{D}. 
$\blacksquare$

\subsection{Reminder on the 3d COHA}\label{subsec:app_KS_crit_coh}

In this section, we focus on the quiver $Q_3$ with potential $W_3$  and discuss 
in detail the critical $3d$ COHA defined in \cite{KS} 
of this quiver with potential. The general setup is the one of \S\ref{sec:2.1}. 
As we have already mentioned, to save the 
terminology we will call ``critical $3d$ COHA'' simply COHA,
in case if does not lead to a confusion. Nevertheless in this section we will
often keep the adjective ``critical" in order to make it easier for the reader 
to compare our exposition with the foundational paper \cite{KS}.

%Recall that the  quiver $Q_3$ has only one vertex, and three loops $B_1, B_2, B_3$. The potential $W_3$  is given by $W_3=B_1[B_2, B_3]$. 

For any  $n\in \bbN$, let $V$ be a $n$-dimensional vector space. 
Then the  space of representations of $Q_3$ with dimension vector $n$ consists 
of triple of $n\times n$ matrices: 
\[
\{(B_1, B_2, B_3)\mid B_i\in \mathfrak{gl}_n, i\in 
\underline{3}\}=\mathfrak{gl}_n^3. 
\]
The group $\GL_{n}$ acts on the representation space $\mathfrak{gl}_n^3$ via 
simultaneous conjugation. 

Denote by $\tr(W_3)_n$ the trace of $W_3$ on $\mathfrak{gl}_n^3$. Note that 
we have
\[
\tr(W_3)_n=\tr(B_1[B_2, B_3])=\tr( B_2[B_3, B_1] )=\tr(B_3[B_1, B_2]). 
\]
Let $\Crit(\tr (W_3)_n)$ be the critical locus of $\tr (W_3)_n$. It is easy to 
see that the following holds.
\begin{lmm}\label{lem:crit}
\[
\Crit(\tr (W_3)_n)=\Big\{(B_1, B_2, B_3)\in  \mathfrak{gl}_n^3\mid [B_1, B_2]=0, 
[B_2, B_3]=0, [B_1, B_3]=0\Big\}. 
\]
\end{lmm}
The vanishing cycle complex on 
$\mathfrak{gl}_n^3$, denoted by $\varphi_{(W_3)_{n}}$, is 
supported on $\Crit(\tr (W_3)_n)$.

Consider the graded vector space
\[
{\mathcal H}^{(Q_3,W_3)}:=
\bigoplus_{n\in \N}{\mathcal H}^{(Q_3,W_3)}(n)=
\bigoplus_{n\in \N}H_{c, 
\GL_n}^*(\mathfrak{gl}_n^3, \varphi_{(W_3)_n})^\vee. \]

For $n_1,n_2\in \bbN$ such that $n=n_1+n_2$, let $V_1\subset V $ be a subspace 
of $V$ of dimension $n_1$. 
We write $G_n:=\GL_{n}$ for short. Let $P\subset \GL_n$ be the parabolic 
subgroup preserving the subspace $V_1$ and $L:=G_{n_1}\times G_{n_2}$ be the 
Levi subgroup of  $P$. We have the corresponding Lie algebras 
\[
\Lie(G_n)=\mathfrak{gl}_n,  \,\  
\Lie(L)=\mathfrak{l}=\mathfrak{gl}_{n_1}\times \mathfrak{gl}_{n_2}.
\]
We write $\Lie(P)= \mathfrak{p}_{n_1, n_2}$ or $\mathfrak{p}$ for simplicity. Its nilpotent radical  is denoted by $\mathfrak{n}_{n_1, n_2}$. 
We have $\mathfrak{p}_{n_1, n_2}^3=\{(B_i)_{i\in \underline{3}}\in 
\mathfrak{gl}_{n}^3\mid B_i(V_1)\subset V_1, \,\ \text{for all $i\in 
\underline{3}$}\}$. We have the following correspondence of $L$-varieties. 
\begin{equation}\label{basic corresp}
\xymatrix
{\mathfrak{gl}_{n_1}^3\times \mathfrak{gl}_{n_2}^3&\mathfrak{p}_{n_1, n_2}^3 
\ar[l]_(0.4){p}\ar[r]^(0.4){\eta} &\mathfrak{gl}_{n_1+n_2}^3,
}\end{equation}
where $p$ is the natural projection induced from $\mathfrak{p}\to \mathfrak{l}$, and $\eta$ is the embedding. 
The trace functions $\tr (W_3)_{n_i}$ on $\mathfrak{gl}_{n_i}^3$ induce a 
function $\tr (W_3)_{n_1}\boxplus \tr (W_3)_{n_2}$
on the product $\mathfrak{gl}_{n_1}^3\times \mathfrak{gl}_{n_2}^3$.
We define $\tr(W_3)_{n_1, n_2}$ on $\mathfrak{p}_{n_1, n_2}^3$ to be
\[
\tr(W_3)_{n_1, n_2}:=p^*(\tr (W_3)_{n_1}\boxplus \tr (W_3)_{n_2})=\eta^*(\tr 
(W_3)_{n_1+n_2}).
\]
Note that we have
$p^{-1} (\Crit(\tr (W_3)_{n_1})\times \Crit(\tr (W_3)_{n_2}))
\supsetneqq
\eta^{-1}(\Crit(\tr (W_3)_{n_1+n_2}))$. Indeed, 
\begin{align*}
&p^{-1} (\Crit(\tr (W_3)_{n_1})\times \Crit(\tr (W_3)_{n_2}))=\{(B_i)_{i\in 
\underline{3}}\in \mathfrak{p}^3_{n_1, n_2}\mid [\pr(B_i), \pr(B_j)]=0, i, j\in 
\underline{3}\}\\
&\text{while} \,\ \eta^{-1}(\Crit(\tr (W_3)_{n_1+n_2}))=\{(B_i)_{i\in 
\underline{3}}\in \mathfrak{p}^3_{n_1, n_2}\mid [B_i, B_j]=0, i, j\in 
\underline{3}\}
\end{align*}
By Lemma \ref{lem:crit}, we have $\Crit(\tr (W_3)_n)\subseteq (\tr 
(W_3)_n)^{-1}(0)$.  Therefore, we have a Thom-Sebastiani isomorphism \cite{M01}.

The multiplication $m^{\crit}: {\mathcal H}^{(Q_3,W_3)}(n_1)\otimes {\mathcal H}^{(Q_3,W_3)}(n_2)\to {\mathcal H}^{(Q_3,W_3)}(n)$ is the composition of 
the following \cite{KS}. For simplicity, we omit the shifting of cohomological degrees in what follows.
\begin{enumerate}
\item The Thom-Sebastiani isomorphism   
\begin{align*}
&{\mathcal H}^{(Q_3,W_3)}(n_1)\otimes {\mathcal H}^{(Q_3,W_3)}(n_2)
\cong 
H_{c, L}^*(\mathfrak{l}_{n_1, n_2}^3, \varphi_{(W_3)_{n_1} \boxplus \tr 
(W_3)_{n_2}} )^\vee.
\end{align*}

%---------------------------------------------------------------------------
\item 
Using the fact that $\mathfrak{p}_{n_1, n_2}^3$ is an affine bundle over 
$\mathfrak{l}_{n_1, n_2}^3$ with fiber $\mathfrak{n}_{n_1, n_2}^3$, and
$\tr (W_3)_{n_1, n_2}$ is the pullback 
of $\tr (W_3)_{n_1} \boxplus \tr (W_3)_{n_2}$, we have

\begin{align*}
p^*:H_{c, L}^*(\mathfrak{l}_{n_1, n_2}^3,
 \varphi_{(W_3)_{n_1} \boxplus \tr (W_3)_{n_2}})^\vee
\cong 
&H_{c, L}^*(\mathfrak{p}_{n_1, n_2}^3, 
\varphi_{(W_3)_{n_1, n_2} } )^{\vee}
% \\ \cong &H_{c, G}^*(G\times_{\bold{M}_{\Gamma, n_1, n_2}, \varphi_{
%W_{n_1, n_2} } )^{\vee}.
\end{align*}
%---------------------------------------------------------------------------
\item
Using the fact $\tr (W_3)_{n_1, n_2}$ is the restriction of $\tr (W_3)_{n}$ to
$\mathfrak{p}_{n_1, n_2}^3$, we have
\begin{align*}
\eta_*:H_{c, L}^*(\mathfrak{p}_{n_1, n_2}^3, 
\varphi_{(W_3)_{n_1, n_2} } )^{\vee}\to &
H_{c, L}^*(\mathfrak{gl}_n^3, 
\varphi_{(W_3)_{n} } )^{\vee}.
\end{align*}
\item 
Pushforward along 
$G\times_{P} \mathfrak{gl}_n^3\to \mathfrak{gl}_n^3, (g, m)\mapsto gmg^{-1}$,
 we obtain 
\begin{align*}
H_{c,L}(\mathfrak{gl}_n^3, 
\varphi_{(W_3)_{n} } )^{\vee}&\cong H_{c,P}(\mathfrak{gl}_n^3, 
\varphi_{(W_3)_{n} } )^{\vee}\cong H_{c, G}^*(G\times_{P} 
\mathfrak{gl}_n^3, 
\varphi_{(W_3)_{n} } )^{\vee}\\
&\to 
H_{c, G}^*(\mathfrak{gl}_n^3, 
\varphi_{(W_3)_{n} } )^{\vee}={\mathcal H}^{(Q_3,W_3)}(n).
\end{align*}
\end{enumerate}

It was proved in \cite{KS} that the multiplication $m^{\crit}$  is associative.

\begin{defn} \cite[\S 7.6]{KS} \label{def:crit COHA in KS}
The COHA of the quiver with potential $(Q_3, W_3)$ is defined as an associative algebra
given by the graded vector space
\[
{\mathcal H}^{(Q_3,W_3)}:=\bigoplus_{n\in \N}H_{c, 
\GL_n}^*(\mathfrak{gl}_n^3, \varphi_{(W_3)_n})^\vee,\] endowed with the 
multiplication $m^{\crit}$ described above.
\end{defn}

\subsection{Main result}
Let $\tr W_3^{fr}$ be the trace function on $\calM_{\vec{r}}(n)$ of the quiver with potential $(Q^{fr}_3, W^{fr}_{3})$. Let $\varphi_{W_3^{fr}}$ be the  vanishing cycle complex on $\calM_{\vec{r}}(n)$ with support on the critical locus of $\tr W_3^{fr}$. 

Recall the group $\GL_n\times \GL_{\vec{r}} \times {{\bf T}_2}$ acts on $\calM_{\vec{r}}(n)$. 
Consider the equivariant cohomology 
\[
V_{r_1,r_2,r_3}(n):=H_{c, \GL_n\times \GL_{\vec{r}} \times {{\bf T}_2}}(\calM_{r_1, r_2, r_3}(n)^{st},\varphi_{W_3^{fr}})^{\vee}
\]
of $\calM_{r_1, r_2, r_3}(n)^{st}$ with values in $\varphi_{W_3^{fr}}$, and 
\[
V_{r_1, r_2, r_3}:=\bigoplus_{n\in 
\bbN}V_{r_1,r_2,r_3}(n). 
\]

The main theorem of this paper is the following. 
\begin{thm}\label{thm:main}
The vertex operator algebra $\calW_{r_1,r_2,r_3}$ acts on  $V_{r_1, r_2, r_3}$. 
\end{thm}

We prove Theorem \ref{thm:main} through the following steps. 

Step 1: We show that $\H^{(Q_3, W_3)}$ acts on 
$V_{r_1,r_2,r_3}$. 
This action is naturally obtained via a geometric construction, which is similar to that of the Hall multiplication of the 3d COHA in \cite{KS}. 
This is done in 
\S~\ref{sec:CritCOHA_act}.

Step 2: We extend the action in Step 1 to an action of the Drinfeld double 
$D(\mathcal{SH}^{(Q_3, W_3)})\cong \mathcal{SH}^{(Q_3, W_3)}\otimes \calH^0\otimes \mathcal{SH}^{(Q_3, W_3)}$ 
of the spherical 3d COHA $\mathcal{SH}^{(Q_3, W_3)}$. The Drinfeld double $D(\mathcal{SH}^{(Q_3, W_3)})$ will be shown to be isomorphic to the affine Yangian $\affY$. 
This is done in \S~\ref{sec:CartanDrinfeld} and \S~\ref{sec:double} by constructing a Drinfeld coproduct on 
$D(\mathcal{SH}^{(Q_3, W_3)})$ and checking its compatibility with certain restriction map $V_{r_1, r_2, r_3}\to (V_{1, 0, 0})^{\otimes r_1} \otimes (V_{0, 1, 0})^{\otimes r_2}\otimes (V_{0, 0, 1})^{\otimes r_3}$ (Proposition \ref{prop:coprod}). 

Step 3: 
We introduce a central extension of $D(\mathcal{SH}^{(Q_3, W_3)})\cong \affY$ by a polynomial ring with infinitely many variables 
$\C[c_l^{(1)}, c_l^{(2)}, c_l^{(3)}: l\geq 0]$. 
Denote this central extension by $\SH^{\vec{c}}$. 
We construct a coproduct $\Delta^{\vec{c}}$ (different than the Drinfeld coproduct) on $\SH^{\vec{c}}$, and 
deduce a free field realization of $\SH^{\vec{c}}$, i.e., an action on the Fock space 
$(V_{1, 0, 0})^{\otimes r_1} \otimes (V_{0, 1, 0})^{\otimes r_2}\otimes (V_{0, 0, 1})^{\otimes r_3}$. 
Furthermore, this free field realization is compatible with the hyperbolic restriction (defined in \S\ref{subsec:hypLoc}) of the moduli of spiked instantons.
That is, we have the following diagram
\[
\xymatrix@C=4.5em@R=0.5em{
\SH^{\vec{c}}\ar[r]^{(\Delta^{\vec{c}})^{r_1+r_2+r_3}}&(\SH^{\vec{c}})^{\otimes 
(r_1+r_2+r_3)}\\
&\\
\ar@(ul,ur)[]&\ar@(ul,ur)[]\\
V_{r_1,r_2,r_3}\ar[r]^(0.3){\text{hyperbolic}}_(0.3){\text{restriction}}&(V_{1, 0, 0})^{\otimes r_1} \otimes (V_{0, 1, 0})^{\otimes r_2}\otimes (V_{0, 0, 1})^{\otimes r_3}
}
\]

Step 4: We show that the action of $\SH^{\vec{c}}$ on $V_{r_1,r_2,r_3}$ 
factors through $\calW_{r_1, r_2, r_3}$. This is done by proving that 
the action on the free field realization factors though the screening operators 
from Definition~\ref{def:screen}. This will be done in \S~\ref{sec:Y110}.

\begin{rmk}
As we mentioned in the introduction, examples illustrating our main result should be derived from toric Calabi-Yau $3$-folds
or from dimer models (a.k.a. brane tilings in physics \cite{HK,FHKVW,HV,OY,NiYaYo}). Recall the framework for the latter.
In general a dimer model can be understood as a bipartite graph with vertices colored in two colors $B$ (black) and $W$ (white) on the $2$-dimensional torus $(S^1)^2$ such that
the complement is a disjoint union of finitely many contractible subsets (faces). With such geometric data one can associate a quiver with potential.
Vertices of the quiver are in bijection with faces, edges correspond to the edges of the dual graph, and the orientation is chosen in such a way that if an edge connects two adjacent faces then the vertex $B$ is on the right and $W$ is on the left. 
The potential is a sum (with signs) of paths about vertices.
Framing ranks $r_1,...,r_k$ needed for the notion of framed representation is an additional piece of data. 

The relation of  dimer models, toric Calabi-Yau $3$-folds and stable framed representations of quivers is known in many classes of examples (see e.g. \cite{MozRei}).
E.g. in the case of the resolved conifold (see \cite{Sz1}) dimer models parametrize stable framed representations of
the corresponding Klebanov-Witten quiver with potential. These representations are fixed points for the action
of the two-dimensional torus ${\bf T}_2$. By localization theorem they generate the equivariant cohomology of the moduli spaces of stable framed representations. Our main result claims that it is acted by the corresponding equivariant COHA, its spherical version and their doubles.

\end{rmk}

\section{Action on the cohomology of the moduli of stable framed 
representations}\label{sec:CritCOHA_act}

\subsection{The equivariant COHA}
\label{subsec:torus action}
Recall $q_i$ is the natural 1-dimensional representation of the $i$ th $\C^*$ of $\mathbf{T}_3$. The $\mathbf{T}_2 \subset \mathbf{T}_3$ subtorus satisfies $\prod_{i} q_i=1$. The torus $\mathbf{T}_3$ acts on 
$\mathfrak{gl}_n^3$ component-wise in the straightforward way. Writing the coordinates of ${{\bf T}_2}$ as $(t_1,t_2)$, then ${{\bf T}_2}$ acts on $\mathfrak{gl}_n^3$
 by 
\begin{equation} \label{T3action}
(t_1, t_2)\cdot ( B_1, B_2, B_3)=(t_1B_1, t_2B_2, (t_1t_2)^{-1}B_3). 
\end{equation}
Under the action, the potential $W_3=\tr(B_3[B_1, B_2])$ is invariant. 
The definition of $\H^{(Q_3, W_3)}$ can be modified to encode the ${{\bf T}_2}$ action. 
\begin{align*}
&\H^{(Q_3, W_3),_{{\bf T}_2}}:=\bigoplus_{n\in \N}H_{c, \GL_n\times  
	{{\bf T}_2}}^*(\mathfrak{gl}_n^3, \varphi_{(W_3)_n})^\vee 
%\\ \text{and } & V_{r_1, r_2, r_3}:=\bigoplus_{n\in \bbN} H_{c, \GL_n\times \GL_{\vec{r}}\times 	{{\bf T}_2}}^* (\calM_{\vec{r}}(n)^{st}, \varphi_{W_3^{fr}})^{\vee}. 
\end{align*}
This gives a ${\bf T}_2$-equivariant version of COHA.

As the equivariant COHA is more relevant to the present paper, in lieu of the non-equivariant one, from now on we will not consider the latter, and hence denote the former by $\H^{(Q_3, W_3)}$ instead. 

\subsection{Action of the 3d COHA}
\label{subsec:Action3dCOHA}
The general construction of representations of COHA from 
stable framed objects can be found in \cite[Section 4]{S}. 
In this section, we deal with the special stable framed object $\calM_{\vec{r}}(n)^{st}$ of the quiver with potential $(Q_3, W_3)$. Due to the difference between \cite[Section 4]{S} and here \ref{rmk_framestable_general}, we include the proof of the following.
\begin{thm}\label{thm:critaction}
The algebra   ${\mathcal H}^{(Q_3,W_3)}$ acts on 
\[
V_{r_1, r_2, r_3}=\bigoplus_{n\in \bbN} H_{c, \GL_n\times \GL_{\vec{r}}}^*(\calM_{\vec{r}}^{st}(n), \varphi_{W_3^{fr}})^{\vee}, \]
for any $(r_1, r_2, r_3)\in \bbN^3$. 
\end{thm}

	More precisely,
	there are two  natural actions by correspondences of the equivariant COHA 
	on $V_{r_1,r_2,r_3}$. 
	One is given by ``creation operators'', while the other one is given by 
	``annihilation operators'' (see \cite{S}).
	Combining together these two actions one can define the action of an associative  
	algebra  denoted by $D(\H^{(Q_3,W_3), {\bf T}_2})$ or $D({\mathcal H}({Coh_{{\bf 
				T}_2}}(\C^3)_0)$, which we will call the (equivariant) double COHA (see 
	\cite{S}).

Below, instead of constructing the equivariant double by the  Nakajima-style construction mentioned  above,
we are going to utilize a more familiar approach by introducing a Hopf algebra structure on the equivariant spherical
COHA and making the corresponding Drinfeld double.
More precisely, let $\mathcal{SH}^{(Q_3,W_3), {\bf T}_2}$ be the spherical subalgebra of $\mathcal{H}^{(Q_3,W_3), {\bf T}_2}$ generated by elements in $\mathcal{H}^{(Q_3,W_3), {\bf T}_2}(1)=H_{c, \GL_1\times  
	{{\bf T}_2}}^*(\mathfrak{gl}_1^3, \C)^\vee$. 
Let $D(\mathcal{SH}^{(Q_3,W_3), {\bf T}_2})$ be the Drinfeld double of $\mathcal{SH}^{(Q_3,W_3), {\bf T}_2}$. 
In Sections \ref{sec:CartanDrinfeld} and \ref{sec:double}, we spell out the action of this double $D(\mathcal{SH}^{(Q_3,W_3), {\bf T}_2})$ on $V_{r_1,r_2,r_3}$, and identify it with the affine Yangian of $\widehat{\mathfrak{gl}(1)}$.

\subsubsection{A correspondence}
We fix a flag $F_1\subset  F_2 \subset V$ of linear subspaces of $V$, with $\dim(F_1)=n'$, and $\dim(F_2)=n'+n''$. Let $F_3:=V/F_2$ be the quotient with $\dim(F_3)=n'''$. Therefore, $\dim(V)=n=n'+n''+n'''$, and we have the exact sequence:
\[
0\to F_1\to F_2 \to V\to F_3\to 0
\]
Similarly, we fix a flag $(E_{r_i})_1\subset  (E_{r_i})_2 \subset E_{r_i}$ of linear subspaces of $E_{r_i}$, with $\dim(E_{r_i})_1)=r_i'$, and $\dim(E_{r_i})_2=r_i'+r_i''$.  We have the exact sequence
\[
0\to (E_{r_i})_1\to (E_{r_i})_2  \to E_{r_i} \to  (E_{r_i})_3\to 0.
\]
Define the following subvariety of $\calM_{ \vec{r}}(n)$. 
\begin{align*}
F_{(\vec{r}', \vec{r}'', \vec{r}''')}(n', n'', n''' ):=&
\{(B_i, I_{ab}, J_{ab})\in \calM_{ \vec{r}}(n)\mid  i\in \underline{3}, \{a, b, c\}=\{1, 2, 3\},  \\ &
B_i(F_1)\subset F_1, I_{ab}((E_{r_c})_1)\subset F_1, J_{ab}(F_1 )\subset (E_{r_c})_1, \\&
B_i(F_2)\subset F_2,  I_{ab}((E_{r_c})_2)\subset F_2, J_{ab}(F_2 )\subset (E_{r_c})_2\}
\end{align*}
We have the commutative diagram
\begin{equation}\label{eq:asso}
\xymatrix@C=0.1em @R=0.5em{
	&\calM_{ \vec{r}}(n)&\\
	&F_{(\vec{r}', \vec{r}'', \vec{r}''')}(n', n'', n''' )\ar[ld]_{p_{1, 23}}\ar[rd]^{p_{12, 3}} \ar[dd]^{p}
	\ar@{}[u]|-*[@]{\subset}_{\eta}
	&\\
	\calM_{ \vec{r}'}(n')\times
	Z_{(\vec{r}'', \vec{r}''')}(n'', n''') \ar[rd]^{p_{23}}  & &  
	Z_{(\vec{r}', \vec{r}'' )}(n', n'')\times 
	\calM_{ \vec{r}'''}(n''') \ar[ld]_{p_{12}}
	\\
	&\calM_{ \vec{r}'}(n')\times \calM_{ \vec{r}''}(n'')\times  \calM_{ \vec{r}'''}(n''')&
}
\end{equation}
where $p_{1, 23}$ is obtained by restriction to the subspace $F_1$, and projection to the quotient $V/F_{1}$. 
The map $p_{23}$ is defined using the short exact sequence
\[
0\to F_2/F_1\to V/ F_1\to F_3\to 0. 
\]
The map $p_{12, 3}$ is obtained by restriction to the subspace $F_2$, and projection to the quotient $F_3$.
The map $p_{23}$ is defined using the short exact sequence
\[
0\to F_1\to F_2\to F_2/F_1\to 0. 
\]
\begin{exa}\label{ex_corresp_special}
As a special case when $F_2=V$ and $(E_{r_i})_2=E_{r_i}$ $V_1\subset V$, $i=1, 2, 3$, we write
	\[
	\dim(V_1)=n_1, \,\ \dim(E_{r_i})_1=r_i', \,\ i=1, 2, 3.
	\]
The variety $F_{(\vec{r}', \vec{r}'', \vec{r}''')}(n', n'', n''' )$ will be written as $	Z_{(\vec{r}', \vec{r}'' )}(n', n'')\subseteq \calM_{ \vec{r}}(n)$. 
	It gives rise to the following correspondence
	\begin{equation}\label{eq:simple corresp}
	\xymatrix{
		\calM_{ \vec{r}'}(n') \times \calM_{\vec{r}''}(n'')
		& 
		Z_{(\vec{r}', \vec{r}'' )}(n', n'')\ar[l]_{p}\ar[r]^{\eta}
		&
		\calM_{\vec{r}}(n),
	}\end{equation}
with maps given by
	\begin{align*}
	&p: (B_i, I_{ab}, J_{ab}) \mapsto \Big( \big((B_i)_{1}, (I_{ab})_1, (J_{ab})_{1}\big), 
	\big((B_i)_{2}, (I_{ab})_2, (J_{ab})_{2}\big)\Big), \\
	&\eta:  (B_i, I_{ab}, J_{ab}) \mapsto (B_i, I_{ab}, J_{ab}) \,\ \text{is the natural inclusion. }
	\end{align*}
	
	\end{exa}
\subsubsection{Stability conditions}
\begin{lmm}\label{lem:stablity}
	\leavevmode
	\begin{enumerate}
		\item
		Notations are as in \eqref{eq:simple corresp}. We have
		\[
		p^{-1}(\calM_{ \vec{r}'}(n')^{st} \times \calM_{\vec{r}''}(n'')^{st})
		\subset \eta^{-1} (\calM_{ \vec{r}}(n)^{st})
		\subset p^{-1}(\calM_{ \vec{r}'}(n') \times \calM_{\vec{r}''}(n'')^{st}). 
		\]
		\item Notations are as in \eqref{eq:asso}. We have
		\[
		p^{-1}(\calM_{ \vec{r}'}(n')^{st} \times \calM_{\vec{r}''}(n'')^{st}
		\times \calM_{\vec{r}'''}(n''')^{st}
		)
		\subset \eta^{-1} (\calM_{ \vec{r}}(n)^{st})\]
	\end{enumerate}
\end{lmm}
{\it Proof. }
We show the statement (1). 
For any element $(B_i, I_{ab}, J_{ab})_{i\in \underline{3}, ab\in \overline{3}}\in p^{-1}(\calM_{ \vec{r}'}(n')^{st} \times \calM_{\vec{r}''}(n'')^{st})$, by definition, 
\begin{align}
&\C\angl{(\vec{B})_{1}} Im((\vec{I})_{1})=V_1, \label{eq:sub}\\
&\C\angl{(\vec{B})_{2}} Im((\vec{I})_{2})=V_2 \label{eq:quot}
\end{align}
Consider the following commutative diagram 
\begin{equation}\label{eq:commu}
\xymatrix{
	0\ar[r]& V_1 \ar[r]&V \ar[r]^{\pi}& V_2 \ar[r]& 0\\
	0\ar[r]& \vec{E_1} \ar[r] \ar[u]_{(\vec{I})_{1}}&\vec{E} \ar[r] \ar[u]_{\vec{I}}& \vec{E_2} \ar[r] \ar[u]_{(\vec{I})_{2}}& 0
}
\end{equation}
where $\vec{E}=(E_{r_1}, E_{r_2}, E_{r_3}))$, $\vec{E_1}=((E_{r_1})_1, (E_{r_2})_1, (E_{r_3})_1)$, 
and $\vec{E_2}$ is the corresponding quotient. 

For any $v\in V$, let $v_2\in V_2$ be the image of $v$ under $\pi$. 
By \eqref{eq:quot}, there exists a function $f$, and a vector $\vec{e_2}\in \vec{E}_2$ such that
$f((\vec{B})_2) (\vec{I})_2 (\vec{e_2})=v_2$. 

Let $\vec{e}\in \vec{E}$ be any lifting of the vector $\vec{e_2}\in \vec{E}_2$. Consider the element 
\[
f(\vec{B}) \vec{I}(\vec{e}) \in V. 
\]
We have $\pi( f(\vec{B}) \vec{I}(\vec{e}))=f((\vec{B})_2) \vec{I}_2(\vec{e}_2)=v_2$ by the commutativity of diagram \eqref{eq:commu}. 
Therefore, there exists an element $v_1 \in V_1$, such that $v=v_1+f(\vec{B}) \vec{I}(\vec{e})$. 
By condition \eqref{eq:sub}, there exists a function $g$, and a vector $\vec{e_1}\in \vec{E}_1$ such that
$g((\vec{B})_1) (\vec{I})_1 (\vec{e_1})=v_1$. Therefore, 
\[
v=g((\vec{B})_1) (\vec{I})_1 (\vec{e_1})+f(\vec{B}) \vec{I}(\vec{e})=
g(\vec{B}) \vec{I} (\vec{e_1})+f(\vec{B}) \vec{I}(\vec{e}). 
\]
This implies the inclusion $p^{-1}(\calM_{ \vec{r}'}(n')^{st} \times \calM_{\vec{r}''}(n'')^{st})
\subset \eta^{-1} (\calM_{ \vec{r}}(n)^{st})$.

The inclusion $\eta^{-1} (\calM_{ \vec{r}}(n)^{st})
\subset p^{-1}(\calM_{ \vec{r}'}(n') \times \calM_{\vec{r}''}(n'')^{st})$ is clear. 
Indeed, by projecting to the quotient spaces $\vec{E}_2,  V_2$, the condition 
$\C\angl{(\vec{B})} Im((\vec{I}))=V$ implies $\C\angl{(\vec{B})_{2}} Im((\vec{I})_{2})=V_2$. 

The assertion (2) follows from a similar argument. This completes the proof. 
$\blacksquare$

\begin{defn}
	\leavevmode
	\begin{enumerate}
		\item When $\vec{E}_1=0$, and $V_1\neq 0$, we have $\calM_{ \vec{0}}(n')^{st}=\emptyset$, as $\C\angl{(\vec{B})_{1}} Im((\vec{I})_{1})=\{0\} \neq V_1$. In this case, we define
		\[
		Z_{(\vec{0}, \vec{r}'')}(n', n'')^{st}:=\eta^{-1} (\calM_{ \vec{r}}(n)^{st}). 
		\] 
		Similarly, we define 
		$
		F_{(\vec{0}, \vec{r}'', \vec{r}''')}(n', n'', n''')^{st}:=p^{-1}(\calM_{0}(n') \times \calM_{\vec{r}''}(n'')^{st} \times \calM_{\vec{r}'''}(n''')^{st} )\cap \eta^{-1} (\calM_{ \vec{r}}(n)^{st}). 
		$
		\item Otherwise, we define
		\[
		Z_{(\vec{r}', \vec{r}'' )}(n', n'')^{st}:=p^{-1}(\calM_{ \vec{r}'}(n')^{st} \times \calM_{\vec{r}''}(n'')^{st}). 
		\]
		Similarly, we define 
		\[
		F_{ (\vec{r}', \vec{r}'', \vec{r}''')}(n', n'', n''')^{st}:=p^{-1}(\calM_{ \vec{r}'}(n')^{st} \times \calM_{\vec{r}''}(n'')^{st} \times \calM_{\vec{r}'''}(n''')^{st}),  
		\] 
		if $\vec{r}'\neq 0$.
	\end{enumerate}
\end{defn}

This gives rise to the following correspondence of $L\times  {\bf T}_2$-varieties. 
\begin{equation}\label{eq:corresp stab}
\xymatrix@R=1.5em
{
	\calM_{ \vec{r}'}(n')^{st}\times
	\calM_{\vec{r}''}(n'')^{st}
	\ar@{^{(}->}[d] & 
	Z_{(\vec{r}', \vec{r}'' )}(n', n'')^{st} 
	\ar[l]_(0.4){p^{st}}\ar[r]^(0.6){\eta^{st}}
	\ar@{^{(}->}[d] &
	\calM_{\vec{r}}(n)^{st}\ar@{^{(}->}[d]\\
	\calM_{ \vec{r}'}(n') \times \calM_{\vec{r}''}(n'') & Z_{(\vec{r}', \vec{r}'' )}(n', n'') 
	\ar[l]_(0.4){p}\ar[r]^(0.6){\eta} & \calM_{\vec{r}}(n) 
}\end{equation}

\subsection{Proof of Theorem~\ref{thm:critaction}}
As a further special case of Example~\ref{ex_corresp_special}, when $\vec{r}'=0$, we have $\vec{r}''=\vec{r}$.
For consistency we write $n'=n_1$ and $n''=n_2$.
 In this case, $(E_{r_c})_1=\{0\}$, and 
	$Z_{(\vec{r}', \vec{r}'' )}(n', n'')$ will be simply denoted by
	$Z_{\mathfrak{p}}$.
	
Taking into account the stability conditions, the correspondence \eqref{eq:corresp stab} induces the following diagram of $L\times {\bf T}_2$-varieties. 
\begin{equation}\label{basic corresp}
\xymatrix@R=1.5em
{
&p^{-1}(\mathfrak{gl}_{n_1}^3\times 
\calM_{\vec{r}}(n_2)^{st}) \ar[ld]_p &\\
\mathfrak{gl}_{n_1}^3\times 
\calM_{\vec{r}}(n_2)^{st} \ar@{^{(}->}[d] & 
 Z_{ \mathfrak{p}} ^{st} \ar[l]_(0.4){p^{st}}\ar[r]^{\eta^{st}}
 \ar@{^{(}->}[d] \ar@{}[u]|-*[@]{\subset}&
 \calM_{\vec{r}}(n)^{st}\ar@{^{(}->}[d]\\
\mathfrak{gl}_{n_1}^3\times 
\calM_{\vec{r}}(n_2) & Z_{ \mathfrak{p}} \ar[l]_(0.4){p}\ar[r]^(0.4){\eta} & 
\calM_{\vec{r}}(n)}\end{equation}
Note that $p:  Z_{ \mathfrak{p}} \to 
\mathfrak{gl}_{n_1}^3\times 
\calM_{\vec{r}}(n_2)$ is an affine bundle, and we have an open embedding $Z_{ \mathfrak{p}}^{st}\subset  p^{-1}(\mathfrak{gl}_{n_1}^3\times 
\calM_{\vec{r}}(n_2)^{st})$.

We now construct an action map
\begin{equation}\label{eq:action}
a^{\crit}: \H^{(Q_3, W_3)}(n_1) \otimes_{H_{ {\bf T}_2}(\pt)} 
V_{r_1, r_2, r_3}(n_2)
\to V_{r_1, r_2, r_3}(n). 
\end{equation} as the composition of the following 
morphisms.
\begin{enumerate}
\item The Thom-Sebastiani isomorphism   
\begin{align*}
& \H^{(Q_3, W_3)}(n_1) \otimes 
V_{r_1, r_2, r_3}(n_2)  \cong 
H_{c, L\times \GL_{\vec{r}}\times {\bf T}_2}^*(\mathfrak{gl}_{n_1}^3\times 
\calM_{\vec{r}}(n_2)^{st}, \varphi_{
(W_3)_{n_1} \boxplus (W_3^{fr})_{n_2}} )^\vee.
\end{align*}
%---------------------------------------------------------------------------
\item 
Using the fact that $p^{-1}(\mathfrak{gl}_{n_1}^3\times 
\calM_{\vec{r}}(n_2)^{st})$ is an affine bundle over 
$\mathfrak{gl}_{n_1}^3\times 
\calM_{\vec{r}}(n_2)^{st}$, and
$\tr (W_3^{fr})_{n_1, n_2}$ is the pullback 
of $\tr (W_3)_{n_1} \boxplus \tr (W_3^{fr})_{n_2}$, we have
\begin{align*}
p^*: &H_{c, L\times \GL_{\vec{r}}\times {\bf T}_2}^*(\mathfrak{gl}_{n_1}^3\times 
\calM_{\vec{r}}(n_2)^{st}, \varphi_{
(W_3^{fr})_{n_1} \boxplus (W_3)_{n_2}} )^\vee \\
\cong
&H_{c, L\times \GL_{\vec{r}}\times {\bf T}_2}^*(p^{-1}(\mathfrak{gl}_{n_1}^3\times 
\calM_{\vec{r}}(n_2)^{st}), 
\varphi_{(W_3^{fr})_{n_1, n_2}} )^{\vee}
\end{align*}
%---------------------------------------------------------------------------
\item By Lemma \ref{lem:restriction to open}, let $j: Z_{ 
\mathfrak{p}}^{st}\subset  p^{-1}(\mathfrak{gl}_{n_1}^3\times 
\calM_{\vec{r}}(n_2)^{st})$ be the open 
embedding, we then have a natural map by restriction
\[
j^*: H_{c, L\times \GL_{\vec{r}}\times {\bf T}_2}^*(p^{-1}(\mathfrak{gl}_{n_1}^3\times 
\calM_{\vec{r}}(n_2)^{st}), 
\varphi_{(W_3^{fr})_{n_1, n_2} } )^{\vee}
\to H_{c, L\times \GL_{\vec{r}}\times {\bf T}_2}^*(Z_{ \mathfrak{p}}^{st}, 
\varphi_{(W_3^{fr})_{n_1, n_2}} )^{\vee}
\]
%---------------------------------------------------------------------------
\item
Using the fact $\tr (W_3^{fr})_{n_1, n_2}$ is the restriction of $\tr 
(W_3^{fr})_{n}$ on $\calM_{\vec{r}}(n)^{st}$ to
$Z_{ \mathfrak{p}}^{st}$. We have
\begin{align*}
\eta^{st}_*:  H_{c, L\times \GL_{\vec{r}}\times {\bf T}_2}^*(Z_{ \mathfrak{p}}^{st}, 
\varphi_{(W_3^{fr})_{n_1, n_2}} )^{\vee}
\to &
H_{c, L\times \GL_{\vec{r}}\times {\bf T}_2}^*(\calM_{\vec{r}}(n)^{st}, 
\varphi_{(W_3^{fr})_{n}} )^{\vee}.
\end{align*}
\item 
Pushforward along 
$G\times_{P} \calM_{\vec{r}}(n)^{st}\to \calM_{\vec{r}}(n)^{st}, (g, 
m)\mapsto gmg^{-1}$,
 we get 
\begin{align*}
&H_{c, L\times \GL_{\vec{r}}\times {\bf T}_2}^*(\calM_{\vec{r}}(n)^{st}, \varphi_{
(W_3^{fr})_{n}} )^{\vee}
\cong H_{c, P\times \GL_{\vec{r}}\times {\bf T}_2}^*(\calM_{\vec{r}}(n)^{st}, \varphi_{
(W_3^{fr})_{n}} )^{\vee}\\
\cong &H_{c, G\times \GL_{\vec{r}}\times {\bf T}_2}^*(G\times_{P} 
\calM_{\vec{r}}(n)^{st}, \varphi_{(W_3^{fr})_{n} } )^{\vee}\to 
H_{c, G\times \GL_{\vec{r}}\times {\bf T}_2}^*(\calM_{\vec{r}}(n)^{st} , 
\varphi_{(W_3^{fr})_{n} } )^{\vee}.
\end{align*}
\end{enumerate}
The map $a^{\crit}$ gives rise to an action of $\H^{(Q_3, W_3)}$ on 
$V_{r_1, r_2, r_3}$ by a similar argument 
as the proof of associativity of the product on $\H^{(Q_3, W_3)}$ given in 
\cite{KS}. 

This completes the proof of Theorem~\ref{thm:critaction}.

\section{The extended COHA and Drinfeld coproduct}\label{sec:CartanDrinfeld}
In this section we introduce an extended 3d COHA $\calH^0\ltimes \H^{(Q_3, W_3)}$ by adding to $\H^{(Q_3, W_3)}$ a polynomial algebra $\calH^0$ over $H_{c,{\bf T}_2}(\pt)$ in infinitely many variables. 
We call $\calH^0$ the Cartan algebra, as it will be the commutative subalgebra in the triangular decomposition of the double COHA $D(\mathcal{SH}^{(Q_3, W_3)})=\mathcal{SH}^{(Q_3, W_3)}\otimes \calH^0\otimes \mathcal{SH}^{(Q_3, W_3)}$. 
Furthermore, we extend the action of $\H^{(Q_3, W_3)}$ on  $V_{r_1, r_2, r_3}$ 
to an action of the extended COHA $\calH^0\ltimes \H^{(Q_3, W_3)}$ . 
In particular, the action of $\calH^0$ is obtained using the tautological bundles on $\calM_{\vec{r}}(n)^{st}$. 

\subsection{The tautological bundles}
\label{sub:taut bundle}

\label{subsec:taut}
Recall we have $\calM_{\vec{r}}^{st}=\coprod_{n\in \N} \calM_{\vec{r}}^{st}(n)$ for the space of 
stable representations of $(Q_3^{fr}, W^{fr}_3)$. 
Let $\calV_n$ be the $\GL_n$-equivariant vector bundle 
\[
\calV_n:=\calM_{\vec{r}}^{st}(n)\times \C^n
\] on $\calM_{\vec{r}}^{st}$, where $\C^n$ is the standard 
representation of $\GL_n$, and $\GL_n$ acts on $\calM_{\vec{r}}^{st}(n)\times \C^n$ diagonally. 
The equivariant vector bundle  $\calV_n$ descents to a vector bundle on $\calM_{\vec{r}}^{st}(n)/\GL_n$, which is still denoted by $\calV_n$ without raising confusions. We will consider its K-theory class in $K_{\GL_n}(\calM_{\vec{r}}^{st}(n))\cong K(\calM_{\vec{r}}^{st}(n)/\GL_n)$. The collection $\calV:=(\calV_n)_{n\ge 0}$ can be thought of as a single bundle on $\calM_{\vec{r}}^{st}$.

%Equivalently, $\calV_n$ is the following tautological bundle 
%\[\calV_n=\calM_{\vec{r}}^{st}(n) \times_{\GL_n}\C^n
%\] on $\calM_{\vec{r}}^{st}(n) /\GL_n$. 

For  $i=1, 2, 3$, let $\calE_{r_i}:=\calM_{\vec{r}}^{st} \times \C^{r_i} $ be the trivial $\C^{r_i}$ bundles on $\calM_{\vec{r}}^{st}$ with the natural (non-trivial) action of $\GL_{\vec{r}}$. 

In this section, we use the torus ${{\bf T}_2}$-action as in Section \ref{subsec:torus 
action}.  For integers $(m_1, m_2, m_3)\in \Z^3$, we define a ${{\bf T}_2}$-module structure on $\C$ 
by
\[
(t_1, t_2, t_3) \cdot v:=t_1^{m_1}t_2^{m_2}t_3^{m_3}v, \,\ (t_1, t_2, t_3)\in 
{{\bf T}_3}  \,\ \text{with} \,\ t_1t_2t_3=1, v\in \C, 
\]
and denote this $1$-dimensional ${{\bf T}_2}$-module by $L(m_1, m_2, m_3)$. 
In general for a ${{\bf T}_2}$-module $V$, we use the notational convention (see (2.8.2) of 
\cite{Nak99}) $
q_1^{m_1}q_2^{m_2}q_3^{m_3}V:=L(m_1, m_2, m_3)\otimes V$.  In particular, $q_i$ is the natural 1-dimensional representation of the $i$th $\C^*$ in $\mathbf{T}_3$.
This notation admits an obvious generalization to any $\C$-linear symmetric 
monoidal category endowed with a ${{\bf T}_2}$-action, in particular to the 
category of ${{\bf T}_2}$-equivariant vector bundles on a stack or a moduli space.

We have the following sequence of morphisms  of $(\GL_n\times {{\bf T}_2})_{n\ge 0}$-equivariant vector 
bundles on the connected components $\calM_{\vec{r}}^{st} (n)$ of $\calM_{\vec{r}}^{st} $:
\[
\xymatrix{
\calV \ar[r]^(0.2){d_1}&
{\begin{matrix}
(q_1 q_2)^{-1}\calV
\oplus (q_1 q_3)^{-1} \calV
\oplus (q_2 q_3)^{-1} \calV\\
\oplus q_1^{-1} \calE_{r_1}
\oplus q_2^{-1} \calE_{r_2}
\oplus q_3^{-1} \calE_{r_3}
\end{matrix}}
\ar[r]^(0.6){d_2} &
{\begin{matrix}
q_1^{-1}\calV
\oplus q_2^{-1} \calV
\oplus q_3^{-1} \calV\\
\oplus \calE_{r_1}
\oplus  \calE_{r_2}
\oplus  \calE_{r_3}
\end{matrix}}\ar[r]^(0.7){d_3} & \calV
}
\]
where the morphisms are given by
\[
d_1=\begin{bmatrix} 
B_3\\
-B_2\\
B_1\\
J_{23}\\
-J_{13}\\
J_{12}
\end{bmatrix} \,\
d_2=
\begin{bmatrix} 
B_2 & B_3 &0& I_{23} &0 &0 \\
-B_1 & 0 &B_3& 0& I_{13} &0 \\
0 & -B_1 &-B_2& 0&0& I_{12}  \\
J_{12} & 0 &0 & 0&0& 0  \\
0 & J_{13} &0 & 0&0& 0  \\
0 & 0 &J_{23} & 0&0& 0  
\end{bmatrix}, 
\]
and 
$
d_3=\begin{bmatrix} 
B_1& B_2&B_3&I_{12}& I_{13} &I_{23}
\end{bmatrix} . 
$
\begin{lmm}
The ADHM equations 
\[
[B_a, B_b]+I_{ab}J_{ab}=0, B_aI_{bc}=0, J_{ab} B_{c}=0, \,\ \text{where}\,\ \{a, b, c\}=\{1, 2, 3\}
\]
imply
$
d_2 \circ d_1=0, \,\ \text{and}\,\ d_3 \circ d_2=0. 
$
\end{lmm}

\begin{lmm}
The stability condition \eqref{eq:stabN} implies the map $d_3$ is surjective. 
\end{lmm}
{\it Proof.}
For any $v\in \C^n$, by \eqref{eq:stabN}, we have
\[
v=f_1(B_1, B_2, B_3) I_{12} (v_3)
+f_2(B_1, B_2, B_3) I_{13} (v_2)
+f_3(B_1, B_2, B_3) I_{23} (v_1), 
\]
for some polynomials $f_1, f_2, f_3$ and $v_i\in \C^{r_i}$, $i=1, 2, 3$.  
Therefore, $v\in \Image(d_3)$. 
This implies that $d_3$ is surjective. 
$\blacksquare$

Note that the map $d_1$ is not injective in general. 
%See \cite[Section 2.9]{Nak99} for example. 

Consider the following tautological element in 
$K_{\GL_n\times \GL_{\vec{r}}\times 
{{\bf T}_2}}(\calM_{r_1, r_2, 
r_3}^{st}(n))$, which is obtained from the alternating sum of the above 
complex. 
\[
\calF(\calV_n, \calE_{\vec{r}}):= 
(q_1-q_1^{-1}+q_2-q_2^{-1}+q_3-q_3^{-1})\calV_n
+(q_1^{-1}-1)\calE_{r_1}
+(q_2^{-1}-1)\calE_{r_2}
+(q_3^{-1}-1)\calE_{r_3}\]
Recall the following correspondence from \eqref{basic corresp}, which is used to define the action of  $\H^{(Q_3, W_3)}$ on $V_{r_1, r_2, r_3}$, where $n=n_1+n_2$:
\[
\xymatrix@R=1.5em{
\mathfrak{gl}_{n_1}^3\times \calM_{\vec{r}}^{st}(n_2)
& Z_{ \mathfrak{p}}^{st} \ar[l]_(0.3){p^{st}}\ar[r]^(0.3){\eta^{st}}
&\calM_{\vec{r}}^{st}(n)
}\]
Let $V_{n_1}:=\mathfrak{gl}_{n_1}^{3}\times \C^{n_1}$ be the  
$\GL_{n_1}$-equivariant vector bundle on $\mathfrak{gl}_{n_1}^{3}$. 
The following lemma can be proved by a straightforward calculation.
\begin{lmm}\label{Lem:F_L}
On $K_{\GL_{n_1}\times \GL_{n_2} \times 
\GL_{\vec{r}}\times {{\bf T}_2}
}(Z_{ \mathfrak{p}} ^{st})$, we have the following equality: 
\begin{align*}
&(\eta^{st})^*\calF(\calV_n, \calE_{\vec{r}})-
(p^{st})^*(1\boxtimes \calF(\calV_{n_2}, \calE_{\vec{r}}))\\
=&
(p^{st})^* \Big((q_1-q_1^{-1}+q_2-q_2^{-1}+q_3-q_3^{-1}) (V_{n_1} \boxtimes 1) \Big). 
 \end{align*} 
\end{lmm}

\subsection{The extended 3d COHA}
We will use the notions of Chern classes in critical cohomology recollected in \S~\ref{subsec:reminder_crit_coh}.
Recall $q_i$ is the natural 1-dimensional representation of the $i$ th $\C^*$ of $\mathbf{T}_3$. Let 
$\hbar_i$ be the equivariant first Chern class of $q_i$, $i=1, 2, 3$. The $\mathbf{T}_2 \subset \mathbf{T}_3$ subtorus satisfies $\prod_{i} q_i=1$, therefore, taking Chern class, we have $\sum_{i}\hbar_i=0$. 
Recall $V_{n}=\mathfrak{gl}_{n}^3 \times \C^n$ is the $\GL_n$-equivariant vector bundle on $\mathfrak{gl}_{n}^3$. 
Let $\nu_1, \cdots, \nu_n \in H_{c, 
\GL_n}^*(\mathfrak{gl}_n^3, \varphi_{(W_3)_n})^\vee$ be the Chern roots of $V_{n}=\mathfrak{gl}_{n}^3 \times \C^n$. 

Let $\calH^0$ be a polynomial algebra over $H_{{\bf T}_2}(\pt)$ in the free variables $\{\psi_j\mid j\in \bbN\}$. 
Consider the generating function
\begin{equation}\label{psi(z)}
\psi(z):=1-(\hbar_1\hbar_2\hbar_3) \sum_{j\geq 0} \psi_{j} z^{-j-1}\in 
\calH^0[[z^{-1}]]. 
\end{equation}
Here $\calH^0[[z^{-1}]]$ consists of infinite power series in $z^{-1}$ with each coefficient in $\calH^0$. 

We define an $\calH^0$-action on the 3d COHA $\H^{(Q_3, W_3)}$ as follows. 
For any $f\in \H^{(Q_3, W_3)}(n)=H_{c, 
\GL_n}^*(\mathfrak{gl}_n^3, \varphi_{(W_3)_n})^\vee$ in the degree $n$ piece of  $\H^{(Q_3, W_3)}$, and the generating series $\psi(z)$ as in \eqref{psi(z)}, 
the action of $\psi(z)$ on $f$ is defined as 
\begin{align}\label{eq:action alg}
\psi(z)\bullet f:=
f \left(\prod_{t=1}^n
\frac{
(z-\nu_t-  \hbar_1)
(z-\nu_t-  \hbar_2)
(z- \nu_t-\hbar_3)}
{(z-\nu_t+  \hbar_1)
(z-\nu_t+  \hbar_2)
(z- \nu_t+\hbar_3)}\right)^- \in  \H^{(Q_3, W_3)}(n)[[z^{-1}]], 
\end{align}
where $^-$ stands for the expansion of the rational function at $z=\infty$. 
In the right hand side of \eqref{eq:action alg} we use the natural structure of $H_{c, 
\GL_n\times \mathbf{T}_3}^*(\mathfrak{gl}_n^3, \varphi_{(W_3)_n})^\vee$ as a module over  $\C[\hbar_1,\hbar_2,\hbar_3,\nu_1,\dots,\nu_n]$. 
In particular, the action of $-(\hbar_1\hbar_2\hbar_3)\psi_j$ on $f\in \H^{(Q_3, W_3)}(n)$ is given by the multiplication of the coefficient of  $z^{-j-1}$ in the expansion of 
$\prod_{t=1}^n
\frac{
(z-\nu_t-  \hbar_1)
(z-\nu_t-  \hbar_2)
(z- \nu_t-\hbar_3)}
{(z-\nu_t+  \hbar_1)
(z-\nu_t+  \hbar_2)
(z- \nu_t+\hbar_3)}$ at $z=\infty$. 
Note that, after taking expansion, each coefficient of $z^{-j-1}$ is a polynomial in the variables $\hbar_1,\hbar_2,\hbar_3,\nu_1,\dots,\nu_n$, thus when applied to $f$ gives a well-defined element in $\H^{(Q_3, W_3)}(n)$. 
Here $\H^{(Q_3, W_3)}(n)[[z^{-1}]]$ consists of infinite power series in $z^{-1}$ with each coefficient in $\H^{(Q_3, W_3)}(n)$ (without any completion). 

\begin{defn}
We define the extended 3d COHA as
\[
\calH^0\ltimes \H^{(Q_3, W_3)}. 
\]
\end{defn}
Recall that $\calH^0\ltimes \H^{(Q_3, W_3)}$ has underlying vector space $\calH^0\otimes_{H_{{\bf T}_2}(\pt)} \H^{(Q_3, W_3)}$ with the ring structure determined by the properties that both $\calH^0$ and $ \H^{(Q_3, W_3)}$ are subalgebras, and that  for any $f\in \H^{(Q_3, W_3)}$ we have $\psi(z) f \psi(z)^{-1}= \psi(z)\bullet f$ given by \eqref{eq:action alg}.

We now follow the  notations of Section \ref{subsec:taut}. 
%Let $\calV_n=\calM_{\vec{r}}^{st}(n)\times \C^n$ be the $\GL_n$-equivariant vector bundle on $\calM_{\vec{r}}^{st}$.  and $\calE_{r_i}=\calM_{\vec{r}}^{st} \times \C^{r_i} $ be the trivial $\C^{r_i}$ bundles on $\calM_{\vec{r}}^{st}$, $i=1, 2, 3$. 
Let $\lambda_1, \cdots, \lambda_n \in V_{r_1, r_2, r_3}(n)$ be the Chern roots of $\calV_n$, and 
$\mu_1^{(i)}, \cdots, \mu_{r_i}^{(i)} \in V_{r_1, r_2, r_3}(n)$ be the Chern roots of $\calE_{r_i}$, $i=1, 2, 3$.

We now define an $\calH^0$-action on $V_{r_1, r_2, r_3}=\bigoplus_{n\in \bbN} V_{r_1, r_2, r_3}(n)$ as follows. 
For any $m\in V_{r_1, r_2, r_3}(n)$, and $\psi(z)$ as in \eqref{psi(z)}, 
\begin{align}
\psi(z)&\cdot m 
:=\lambda_{-1/z}
(\calF(\calV_{n}, \calE_{\vec{r}}))m \notag\\
%=&
%\Big(\prod_{a=1}^{r_1}
%\frac{1-\frac{1}{z}(\mu_a-\hbar_1)}{1-\frac{1}{z}\mu_a}
%\prod_{b=1}^{r_2}\frac{1-\frac{1}{z}(\mu_a-\hbar_2)}{1-\frac{1}{z}\mu_b}
%\prod_{b=1}^{r_3}\frac{1-\frac{1}{z}(\mu_a-\hbar_3)}{1-\frac{1}{z}\mu_c}\\
%\cdot&
%\prod_{d=1}^{n}\frac{
%1-\frac{1}{z}(\lambda_d+\hbar_1)
%}{1-\frac{1}{z}(\lambda_d-\hbar_1)}
%\frac{
%1-\frac{1}{z}(\lambda_d+\hbar_2)
%}{1-\frac{1}{z}(\lambda_d-\hbar_2)}
%\frac{
%1-\frac{1}{z}(\lambda_d+\hbar_3)
%}{1-\frac{1}{z}(\lambda_d-\hbar_3)}\Big) m
%\\
=&
\Big(\prod_{a=1}^{r_1}\frac{z-\mu_a^{(1)}+\hbar_1}{z-\mu_a^{(1)}}
\prod_{b=1}^{r_2}\frac{z-\mu_b^{(2)}+\hbar_2}{z-\mu_b^{(2)}}
\prod_{b=1}^{r_3}\frac{z-\mu_c^{(3)}+\hbar_3}{z-\mu_c^{(3)}} \notag\\
\cdot&
\prod_{d=1}^{n}\frac{z-\lambda_d-\hbar_1}{z-\lambda_d+\hbar_1}
\frac{z-\lambda_d-\hbar_2}{z-\lambda_d+\hbar_2}
\frac{z-\lambda_d-\hbar_3}{z-\lambda_d+\hbar_3}\Big) m \in V_{r_1, r_2, r_3}(n)[[z^{-1}]]. 
\label{action of Y0}
\end{align}
Again, the expansion is taken at $z=\infty$ as above. For simplicity, we omit the $^-$ from the notation. 

 \begin{prp}\label{pro:action of 2/3}
The actions of $ \H^{(Q_3, W_3)}$ and $\calH^0$ on $V_{r_1, r_2, r_3}$ 
can be extended to an action of the extended 3d COHA $\calH^0\ltimes \H^{(Q_3, W_3)}$ on $V_{r_1, r_2, r_3}$.
 \end{prp}
{\it Proof.}
For any $f\in \H^{(Q_3, W_3)}({n_1})$, and $m\in V_{r_1, r_2, r_3}(n_2)$, with $n=n_1+n_2$, we need to show the following equality in 
$V_{r_1, r_2, r_3}(n)$: 
\[
\psi(z)\Big(f (\psi(z)^{-1} m)\Big)=( \psi(z)\bullet f) m. 
\]
By \eqref{eq:action alg}, it suffices to show 
 \begin{equation}\label{eq:action of ext}
 \psi_{n}(z) \Big(f (\psi_{n_2}(z)^{-1} \cdot m)\Big)=\Big(f \prod_{t=1}^{n_1}
\frac{
(z-\nu_t-  \hbar_1)
(z-\nu_t-  \hbar_2)
(z- \nu_t-\hbar_3)}
{(z-\nu_t+  \hbar_1)
(z-\nu_t+  \hbar_2)
(z- \nu_t+\hbar_3)}\Big)(m), 
\end{equation}
where we write a subscript $n$ in $\psi_{n}(z) \cdot v$ to emphasis the action of $\psi(z)\in \calH^0[[z^{-1}]]$ on the element $v \in V_{r_1, r_2, r_3}(n)$. 
Both sides of \eqref{eq:action of ext} are elements in $V_{r_1, r_2, r_3}(n_1+n_2)=V_{r_1, r_2, r_3}(n)$. 
We have the following correspondence from \eqref{basic corresp}: 
\[
\xymatrix@R=1.5em{
\mathfrak{gl}_{n_1}^3\times \calM_{\vec{r}}^{st}(n_2)
& Z_{ \mathfrak{p}}^{st} \ar[l]_(0.3){p^{st}}\ar[r]^(0.3){\eta^{st}}
&\calM_{\vec{r}}^{st}(n)
}\]
By the action of  $\H^{(Q_3, W_3)}$  (used in the first equality) and the projection formula (used in the third equality as follows), the left hand side of \eqref{eq:action of ext} is
\begin{align*}
%-------- line 1-----
& \psi_{n}(z)\Big( f \big(\psi_{n_2}(z)^{-1} m\big)\Big)\\
%----------line 2---
=&
\psi_n(z) (\eta^{st})_* (p^{st})^* (f\boxtimes (\psi_{n_2}(z)^{-1}m))\\
%---------line 3----
=&\lambda_{-1/z}
(\calF(\calV_{n}, \calE_{\vec{r}}))
(\eta^{st})_* (p^{st})^* (f\boxtimes \lambda_{-1/z}(\calF(\calV_{n_2}, 
\calE_{\vec{r}}))^{-1}m)\\
%-------------line 4
=&
(\eta^{st})_*\Big(
\lambda_{-1/z}
((\eta^{st})^*\calF(\calV_{n}, \calE_{\vec{r}}))\cdot
 (p^{st})^* (f  \boxtimes \lambda_{-1/z}
(\calF(\calV_{n_2}, \calE_{\vec{r}})^{-1}m)
\Big)
\\
%-------------line 5
=&
(\eta^{st})_*\Big(
 (p^{st})^* (f\boxtimes m)  \cdot (\lambda_{-1/z}(\eta^{st})^* \calF(\calV_{n}, 
 \calE_{\vec{r}}))\cdot (p^{st})^*
(1\boxtimes \lambda_{-1/z}(\calF(\calV_{n_2}, \calE_{\vec{r}})^{-1}))
\Big)
\\
%-------------line 6
=&
(\eta^{st})_*\Big(
 (p^{st})^* (f\boxtimes m)  \cdot 
(p^{st})^* ((q_1-q_1^{-1}+q_2-q_2^{-1}+q_3-q_3^{-1})  \lambda_{-1/z} (V_{n_1} )\boxtimes 1 )
)
\Big)\\
%-------------line 7
%-------------
=&\Big(f \prod_{t=1}^{n_1}
\frac{
(z-\nu_t-  \hbar_1)
(z-\nu_t-  \hbar_2)
(z- \nu_t-\hbar_3)}
{(z-\nu_t+  \hbar_1)
(z-\nu_t+  \hbar_2)
(z- \nu_t+\hbar_3)}\Big)(m),
\end{align*} 
where the second last equality follows from Lemma \ref{Lem:F_L}. 
This completes the proof. 
$\blacksquare$

\subsection{The Drinfeld coproduct}
\label{subsec:DrinfeldCoprodDef}
Let 
$\mathcal S \mathcal H^{\geq 0}:= \calH^0 \ltimes {\mathcal S \mathcal H}^{(Q_3,W_3)}$ be the extended spherical COHA.
In this section, we describe the Drinfeld coproduct on $\mathcal S \mathcal H^{\geq 0}$
\[
\Delta^{\Dr}:\mathcal S \mathcal H^{\geq 0} \to (\mathcal S \mathcal H^{\geq 0})\widehat{\otimes} (\mathcal S \mathcal H^{\geq 0})
\]
following \cite{YZ2}. 
The  difference with \cite{YZ2} is that 
the Cohomological Hall algebra $\H^{(Q_3, W_3)}$ in the current paper is associated to a quiver with potential. For completeness we add the detailed formulas.

%Similarly  to \S~\ref{sub:taut bundle}, we have the $\GL_n$-equivariant vector bundle $V_n:=\mathfrak{gl}_n^3 \times \C^n$ on $\mathfrak{gl}_n^3$, where $\C^n$ is the standard representation of $\GL_n$, and $\GL_n$ acts on $\mathfrak{gl}_n^3 \times \C^n$ diagonally. 
Let $\{\lambda_1, \lambda_2, \cdots, \lambda_n\}$ be the  $\GL_n$-equivariant Chern roots of $V_n$. 
Let $A=\{1, \cdots, n\}, B=\{n+1, \cdots n+m\}$, and let $\hbar_1+\hbar_2+\hbar_3=0$. Set \begin{align*}
\fac(\lambda_A\mid \lambda_{B}):=&
\prod_{s\in A}
\prod_{t\in B}
\frac{
(\lambda_s-\lambda_t-  \hbar_1)
(\lambda_s-\lambda_t-  \hbar_2)
(\lambda_s- \lambda_t-\hbar_3)}{\lambda_s- \lambda_t}. 
\end{align*}
%Consider 
%\[\widehat{\Phi}( \lambda_{B}\mid \lambda_A):=\frac{\fac(\lambda_A\mid 
%\lambda_{B})}{\fac(\lambda_B\mid \lambda_{A})}=
%\prod_{s\in A}
%\prod_{t\in B}
%\frac{
%(\lambda_s-\lambda_t-  \hbar_1)
%(\lambda_s-\lambda_t-  \hbar_2)
%(\lambda_s- \lambda_t-\hbar_3)}
%{(\lambda_s-\lambda_t+  \hbar_1)
%(\lambda_s-\lambda_t+  \hbar_2)
%(\lambda_s- \lambda_t+\hbar_3) }\]

The coproduct $\Delta^{\Dr}$ on $\mathcal S \mathcal H^{\geq 0}$ is determined as follows (see \cite{YZ2} for details). 
For $\psi(z)\in \calH^0[[z^{-1}]]$, and $P(\lambda_1, \cdots \lambda_n)\in {\mathcal S \mathcal H}^{(Q_3,W_3)}(n)$, we define
\begin{align*}
\Delta^{\Dr}( \psi(z))=\psi(z)&\otimes \psi(z), 
\\
\Delta^{\Dr}(P(\lambda_1, \cdots \lambda_n))
&=\sum_{a+b=n, A=[1, a], B=[a+1, n]} \frac{\prod_{t\in B}\psi(\lambda_t) 
P(\lambda_A\otimes \lambda_B)}{ \fac(\lambda_B\mid \lambda_A)}%\\
%&=
%\sum_{a+b=n}  \frac{\prod_{t\in [a+1, n]}\psi(\lambda_t) P(\lambda_{[1, a]} \otimes 
%\lambda_{[a+1, n]}) \prod_{s\in [1, a], t\in [a+1, n]}(\lambda_t- \lambda_s)}{ 
%\prod_{s\in [1, a], t\in [a+1, n]}(
%\lambda_t-\lambda_s-  \hbar_1)
%(\lambda_t-\lambda_s-  \hbar_2)
%(\lambda_t- \lambda_s-\hbar_3)
%}
\end{align*}
By the same reason as in \cite{YZ2}, this formula defines a coproduct on $\mathcal S \mathcal H^{\geq 0}$. 

In particular, when $n=1$, we have
\[
\calH^{(Q_3,W_3), {\bf T}_2}(1)=H_{c, \GL_1\times  
{{\bf T}_2}}^*(\mathfrak{gl}_1^3, \C)^\vee \cong 
\C[\hbar_1, \hbar_2] \otimes \C[\lambda]. 
\]
In this case, for any $P(\lambda)\in  \calH^{(Q_3,W_3), {\bf T}_2}(1)$, the Drinfeld coproduct is given by
\[
\Delta^{\Dr}(P(\lambda))=\psi(\lambda) \otimes P(\lambda)+P(\lambda)\otimes 1. 
\]

 \subsection{The action of the extended 3d COHA}
\label{subsec:3dCOHA}

Let $\vec{r}'=(r_1', r_2', r_3')\in \Z^3_{\geq 0}, \vec{r}''=(r_1'', r_2'', r_3'')\in \Z^3_{\geq 0}$. 
We denote the sum $\vec{r}'+\vec{r}''$ by $\vec{r}=(r_1, r_2, r_3)$, 
and let $|\vec{r}|=r_1+r_2+r_3$. 
We write $n=n'+n'' \in  \Z_{\geq 0}$, for any $n', n''\in \Z_{\geq 0}$.
\begin{prp}\label{prop:coprod}
%\leavevmode
Let $\mathcal S \mathcal H^{\geq 0}= \calH^0 \ltimes {\mathcal S \mathcal H}^{(Q_3,W_3)}$ be the extended spherical COHA. 
\begin{enumerate}
\item There is a map of vector spaces $l: V_{\vec{r}}\to V_{\vec{r'}}\otimes V_{\vec{r''}}$,  which is an isomorphism up to localization, such that 
\[
l(\alpha\bullet x)=\Delta^{\Dr}(\alpha)\bullet l(x), \] where $\alpha\in  \mathcal S \mathcal H^{\geq 0}$, $x\in V_{\vec{r}}$ and 
$
\Delta^{\Dr}:\mathcal S \mathcal H^{\geq 0} \to (\mathcal S \mathcal H^{\geq 0})\widehat{\otimes} (\mathcal S \mathcal H^{\geq 0})
$ is the Drinfeld coproduct of $\mathcal S \mathcal H^{\geq 0}$.  
%Similar statement holds for $Y^{\leq0}$.
\item The two maps
\[
(\Delta^{\Dr})^{\otimes|\vec{r}|}:\mathcal S \mathcal H^{\geq 0} \to (\mathcal S \mathcal H^{\geq 0})^{\widehat{\otimes} |\vec{r}|}, \,\ \text{ and} \,\  l^{\otimes |\vec{r}|}:V_{\vec{r}}\to 
(V_{1, 0, 0})^{\otimes r_1} \otimes (V_{0, 1, 0})^{\otimes r_2}\otimes (V_{0, 0, 1})^{\otimes r_3}
\] intertwine the actions of  $\mathcal S \mathcal H^{\geq 0}$ on $V_{\vec{r}}$ and $(\mathcal S \mathcal H^{\geq 0})^{\widehat{\otimes} |\vec{r}|}$ on $(V_{1, 0, 0})^{\otimes r_1} \otimes (V_{0, 1, 0})^{\otimes r_2}\otimes (V_{0, 0, 1})^{\otimes r_3}$.  %Similar statement holds for $Y^{\leq 0}$.
\end{enumerate}
\end{prp}
Let $T\subseteq \GL_{\vec{r}}$ be the maximal torus.
Here localization means base change from modules of $H_{{\bf T}_2\times\GL_{\vec{r}}}(\pt)$ to $K_{\vec{r}}$, the field of fractions of $H_{{\bf T}_2\times T}(\pt)$.

Proposition \ref{prop:coprod} can be summarized into the following two commutative diagrams. 
\[
\xymatrix{
\mathcal S \mathcal H^{\geq0} \otimes  V_{\vec{r}} \ar[r] \ar[d]_{\Delta^{\Dr} \times l}\ar[d]& V_{\vec{r}} \ar[d]^{l}\\
(\mathcal S \mathcal H^{\geq0} \widehat{\otimes} \mathcal S \mathcal H^{\geq0})  \otimes  (V_{\vec{r'}}\otimes V_{\vec{r''}}) \ar[r]& V_{\vec{r'}}\otimes V_{\vec{r''}}
}
\] and 
\[
\xymatrix@C=1em{
\mathcal S \mathcal H^{\geq0} \otimes  V_{\vec{r}} \ar[r] \ar[d]_{(\Delta^{\Dr})^{\otimes|\vec{r}|} \times l^{\otimes|\vec{r}|}}& 
V_{\vec{r}} \ar[d]^{ l^{\otimes|\vec{r}|}}\\
{\begin{matrix}
(\mathcal S \mathcal H^{\geq0})^{\widehat{\otimes}|\vec{r}|}  \otimes ((V_{1, 0, 0})^{\otimes r_1} \\ \otimes (V_{0, 1, 0})^{\otimes r_2}\otimes (V_{0, 0, 1})^{\otimes r_3}) \end{matrix}}
\ar[r]& 
((V_{1, 0, 0})^{\otimes r_1} \otimes (V_{0, 1, 0})^{\otimes r_2}\otimes (V_{0, 0, 1})^{\otimes r_3})}
\] 
The proof of Proposition~\ref{prop:coprod} will be given in Appendix \ref{App2}. 

\subsection{Nakajima-type operators}
\label{subNak}
 Let $\fM_{\vec{r}}(n)$ be the quotient $\calM_{\vec{r}}(n)^{st}/\GL_n$. 
Consider the subvariety $C^+_{\vec{r}}(n-1, n)\subset  \fM_{\vec{r}}(n-1)\times \fM_{\vec{r}}(n)$ consisting of the pairs $(V_2, V)\in \fM_{\vec{r}}(n-1)\times \fM_{\vec{r}}(n)$, such that $V_2$ is a quotient representation of $V$. There is a tautological line bundle $\calL$ on $C^+_{\vec{r}}(n-1, n)$ which classifies the one-dimensional sub-representation $V_1\subset V$. 
In the notations of \S~\ref{subsec:Action3dCOHA}, the subvariety $C^+_{\vec{r}}(n-1, n)$ is obtained by the quotient $Z_{\mathfrak{p}}^{st}/P$.

We have the following diagram \[
\xymatrix@R=1.5em @C=1em{
C^+_{\vec{r}}(n-1, n) \ar@{^{(}->}[r]&\fM_{\vec{r}}(n-1)\times \fM_{\vec{r}}(n) \ar[dl]\ar[dr] &\\
\fM_{\vec{r}}(n-1)&&\fM_{\vec{r}}(n).
}
\]
Denote by $p_1: C^+_{\vec{r}}(n-1, n) \to \fM_{\vec{r}}(n-1), p_2: 
C^+_{\vec{r}}(n-1, n) \to \fM_{\vec{r}}(n)$ the composition of the inclusion with the natural projections.

We also have the following commutative diagram
\[
\xymatrix@R=1.5em{
% line 1---
G\times_P(\mathfrak{gl}_{1}^3 \times \calM_{\vec{r}}(n-1)^{st}) 
&G\times_P Z_{\mathfrak{p}}^{st}  \ar[l]_(0.3){p^{st}}\ar[r]^{\eta^{st}} \ar[d]& \calM_{\vec{r}}(n)^{st}\ar[d]\\
%line 2---
\fM_{\vec{r}}(n-1) & C^+_{\vec{r}}(n-1, n) \ar[l]_{p_1} \ar[r]^{p_2}& \fM_{\vec{r}}(n)
}
\]
where the vertical maps are obtained by quotient by $G$. 
Let $\calL$ be the tautological line bundle on $C^+_{\vec{r}}(n-1, n)$. Let $\Psi(f(c_1(\calL)))\in \End (V_{r_1, r_2, r_3})$ be the raising operation given by convolution with $f(c_1(\calL))$ for any polynomial $f(u)\in\C[\hbar_1,\hbar_2,u]$. In other words, let $\alpha \in V_{r_1, r_2, r_3}(n-1)\cong H_{c, \GL_{\vec{r}}\times {{\bf T}_2}}( \fM_{\vec{r}}(n-1), \varphi_{(W^{fr})_{n-1}})^{\vee}$, 
\[
\Psi(f(c_1(\calL)))(\alpha):=p_{2*}\big( f(c_1(\calL)  \cdot p_1^*(\alpha) )\big) \in V_{r_1, r_2, r_3}(n). 
\]

The action of the spherical equivariant COHA $ \mathcal{SH}^{(Q_3, W_3)}$ on $V_{r_1, r_2, r_3}$ can be described in terms of the above Nakajima operators. 
\begin{prp}\label{prop:NakajimaOperators}
For any $f(\lambda)\in  \C[\hbar_1,\hbar_2,\lambda]$, view $f(\lambda)\in \mathcal{SH}^{(Q_3, W_3)}(1)
%=H_{c,\Gg_1}(\M_1, \C)^{\vee} 
\cong \C[\hbar_1,\hbar_2,\lambda]$, 
we have the equality 
\[
\Psi(f(c_1(\calL)))=\Phi(f(\lambda))\] in $\End(V_{r_1, r_2, r_3})$, where $\Phi$ is the action of the 3d spherical equivariant COHA $\mathcal{SH}^{(Q_3, W_3)}$.

\end{prp}
{\it Proof:} By the same proof as in \cite[Theorem 5.6]{YZ1}.
$\blacksquare$

\section{The Drinfeld double of COHA and its action}
\label{sec:double}

We show the extended spherical COHA  $\mathcal S \mathcal H^{\geq 0}$ is isomorphic to the Borel subalgebra of the affine Yangian $\affY$. When the framed quiver variety has a 2d description, the action of  $\affY$ on its cohomology has been known in the literature (see, e.g., \cite{SV1,T}). 
In this special case, the action of $\affY$ is compatible with the action of $\mathcal S \mathcal H^{\geq 0}$ constructed in Section \ref{sec:CartanDrinfeld}. In the current section, we prove that for  framed quiver variety with general framing $(r_1,r_2,r_3)$, the entire Yangian $\affY$ acts on its cohomology.

\subsection{The affine Yangian}
We review some facts of the  affine Yangian of $\widehat{\mathfrak{gl}(1)}$ (see \cite{T} for details). 
Let $\hbar_1, \hbar_2, \hbar_3$ be formal parameters satisfying $\hbar_1+\hbar_2+\hbar_3=0$. 
The affine Yangian 
$\affY$ is an associative 
algebra, generated by the variables $\{e_j, f_j, \psi_j\mid j\in \bbN\}$ with 
the following defining relations. 
%\footnote{Here $\hbar_1=-h_1$, $\hbar_2=-h_2$, and 
%$\hbar_3=-h_3$ compare to the parameters $h_1, h_2, h_3$ in \cite{T}. }
\begin{align}
&[\psi_{i}, \psi_j]=0,  \tag{Y0} \label{Y0}\\
%-------
&[e_{i+3}, e_j]-3[e_{i+2}, e_{j+1}]+3[e_{i+1}, e_{j+2}]-[e_i, 
e_{j+3}]\notag\\
&\phantom{1234567891012345} +\sigma_2([e_{i+1}, e_j]-[e_i, e_{j+1}])=-\sigma_3\{e_i, e_j\}. 
\tag{Y1} \label{Y1}\\
%-------
&[f_{i+3}, f_j]-3[f_{i+2}, f_{j+1}]+3[f_{i+1}, f_{j+2}]-[f_i, 
f_{j+3}]\notag\\
&\phantom{1234567891012345}+\sigma_2([f_{i+1}, f_j]-[f_i, f_{j+1}])=\sigma_3\{f_i, f_j\}. \tag{Y2} 
\label{Y2}\\
%-------
&[e_i, f_j]=\psi_{i+j} \tag{Y3} \label{Y3}\\
%-------
&[\psi_{i+3}, e_j]-3[\psi_{i+2}, e_{j+1}]+3[\psi_{i+1}, e_{j+2}]-[\psi_i, 
e_{j+3}]\notag\\
&\phantom{1234567891012345}+\sigma_2([\psi_{i+1}, e_j]-[\psi_i, e_{j+1}])=-\sigma_3\{\psi_i, e_j\}. 
\tag{Y4} 
\label{Y4}\\
%-------
& [\psi_0, e_j]=0, [\psi_1, e_j]=0, [\psi_{2}, e_j]=2e_j, \tag{Y4$'$} 
\label{Y4'}\\
%-------
&[\psi_{i+3}, f_j]-3[\psi_{i+2}, f_{j+1}]+3[\psi_{i+1}, f_{j+2}]-[\psi_i, 
f_{j+3}]\notag\\
&\phantom{1234567891012345}+\sigma_2([\psi_{i+1}, f_j]-[\psi_i, f_{j+1}])=\sigma_3\{\psi_i, f_j\}  \tag{Y5} \label{Y5}\\
%-------
& [\psi_0, f_j]=0, [\psi_1, f_j]=0, [\psi_{2}, f_j]=-2f_j, \tag{Y5$'$} 
\label{Y5'}\\
%-------
& \Sym_{\mathfrak{S}_3}[e_{i_1}, [e_{i_2}, e_{i_3+1}]]=0, \,\   
\Sym_{\mathfrak{S}_3}[f_{i_1}, [f_{i_2}, f_{i_3+1}]]=0. \tag{Y6}\label{Y6}
\end{align}
where $
\sigma_2=\hbar_1\hbar_2+\hbar_2\hbar_3+\hbar_1\hbar_3$, and 
$\sigma_3=\hbar_1\hbar_2\hbar_3$. 

Let $Y^-, Y^0, Y^+$ be the subalgebra generated by $\{f_j\}$, $\{\psi_j\}$ and 
$\{e_j\}$ respectively. 
Let $Y^{\ge 0}, Y^{\le 0}$ be the subalgebras generated by $Y^0, Y^+$, and $Y^-, 
Y^0$ respectively. The following properties can be found in \cite[Proposition 1.4]{T}. 
\begin{enumerate}
\item $Y^0$ is a polynomial algebra in the generators $\{\psi_j\}$. 
\item $Y^-$ and $Y^+$ are the algebras generated by $\{f_j\}$ and $\{e_j\}$ with 
the defining relations \eqref{Y2}\eqref{Y6}, and \eqref{Y1}\eqref{Y6} respectively. 
\item $Y^{\leq 0}$ and $Y^{\geq 0}$ are the algebras generated by $\{\psi_j, f_j\}$ 
and $\{\psi_j, e_j\}$ with the defining relations 
\eqref{Y0}\eqref{Y2}\eqref{Y5}\eqref{Y5'}\eqref{Y6}, and 
\eqref{Y0}\eqref{Y1}\eqref{Y4}\eqref{Y4'}\eqref{Y6} respectively. 
\item Multiplication induces an isomorphism of vector spaces
\[
m: Y^-\otimes Y^0\otimes Y^+ \to \affY. 
\]
\end{enumerate}

%As has been  mentioned in \cite[Page 5]{T}, the first reference in literature goes back to 
%\cite[Section 3.2]{AS}, with the particular choice $\{\hbar_1, \hbar_2, 
%\hbar_3\} = \{1, -\kappa, \kappa-1\}$, where it was shown to be isomorphic to 
%the algebra $\SH^c$ from \cite{SV2}. 
%However,  the relation $(Y1)$ in \cite{T} and the relation (3.3) in \cite{AS} are 
%different. In \S~\ref{subsec:app_Yangian} we prove Theorem~\ref{thm:positiveYangian}
%(1) below. 

We compare the Borel subalgebra $Y^{\geq 0}$ of the affine Yangian $\affY$ with the extended spherical COHA  
$\mathcal S \mathcal H^{\geq 0}$. Thanks to \cite{SV2,T},  when two of the three coordinates of $\vec{r}=(r_1,r_2,r_3)$ are zero, the entire affine Yangian  $\affY$ has an action on 
$V_{r_1, r_2, r_3}$. 
For instance, let $r_1=r_2=0, r_3=1$. 
Using the dimension reduction in \cite{D}, we have an isomorphism $V_{0, 0, 1} \cong \bigoplus_{n\in \N}H^*(\Hilb^n(\C^2), \C)$, which carries an affine Yangian action \cite{AS, SV2, T}.

\begin{thm}\label{thm:positiveYangian}
\leavevmode
\begin{enumerate}
\item 
There is an algebra homomorphism 
\[
\Psi: Y^+\to \H^{(Q_3, W_3)},  \,\  \text{by the assignment $e_{r} \mapsto \lambda^{r} \in 
\H^{(Q_3, W_3)}(1)\cong \C[\lambda]$}
\] Furthermore, it induces an isomorphism $Y^+ \cong  \mathcal{SH}^{(Q_3, W_3)}$. 
\item The isomorphism in (1) extends to an isomorphism $Y^{\geq 0} \cong \mathcal{SH}^{\geq 0}$. In particular, $Y^{\geq 0} $ acts on $V_{r_1, r_2, r_3}$, for any $(r_1, r_2, r_3)\in \N^3$. 
\item The isomorphism from (2) intertwines the Drinfeld coproduct on $Y^{\ge0}$ and that on $\mathcal{SH}^{\geq 0}$ from \S~\ref{subsec:DrinfeldCoprodDef}.
\item When $r_a=r_b=0$ and $r_c=r$, such that $\{a, b, c\}=\{1, 2, 3\}$, 
the action of $Y^{\geq 0}$ on $V_{r}$ from (2) is compatible with the actions constructed in \cite{SV2, T}.
\end{enumerate}
\end{thm}
%In Section \ref{sec:CartanDrinfeld}, we have constructed an action of $\mathcal{SH}^{\geq 0}$ on $V_{r_1, r_2, r_3}$, for any $(r_1, r_2, r_3)\in \N^3$. The isomorphism in Theorem~\ref{thm:positiveYangian}(2) gives an action of $Y^{\geq 0}$ on $V_{r_1, r_2, r_3}$, for any $(r_1, r_2, r_3)\in \N^3$. Theorem~\ref{thm:positiveYangian}(4) shows, in the case when two of the three coordinates in $\vec{r}=(r_1,r_2,r_3)$ are zero, that this is compatible with the actions in \cite{SV2, T}.
 Using this result, we will show in Section \ref{subsec:wholeY} a stronger statement that the whole Yangian $\affY$ acts on  $V_{r_1, r_2, r_3}$, for any $(r_1, r_2, r_3)\in \N^3$.

We prove Theorem~\ref{thm:positiveYangian}(1) and (4) in \S\ref{subsec:2dCOHA}, and Appendix \ref{subsec:app_Yangian}. 
The statement (3) follows directly by comparing the Drinfeld coproduct formulas on $Y^{\ge0}$ and $\mathcal{SH}^{\geq 0}$. 
We now assume Theorem~\ref{thm:positiveYangian}(1), and prove  Theorem~\ref{thm:positiveYangian}(2). 

{\it Proof of Theorem~\ref{thm:positiveYangian}(2).} 
We need to show that the action of $\calH^0$ on the spherical COHA \eqref{eq:action alg}  is compatible with the relation \eqref{Y4} of $\affY$.

It is shown in  \cite[Proposition 1.5]{T} that the defining relation \eqref{Y4} of $\affY$ is equivalent to 
\[
 P(z, \sigma^+) \psi(z) e_j+ P(\sigma^+, z) e_j \otimes \psi(z)=0, \text{for $j\in \bbN$}, 
\]
where $\sigma^+$ is the shifting operators
\begin{align*}
&\sigma^+: Y^{\geq} \to Y^{\geq}  \,\ \text{determined by} \,\ \psi_{j}\mapsto \psi_j, e_j\mapsto e_{j+1}. 
\end{align*} and $
P(z, w):=(z-w+\hbar_1)(z-w+\hbar_2)(z-w+\hbar_3)$. 
In our case, we take $f$ in \eqref{eq:action alg} to be $e_j=\lambda^j\in  
\H^{(Q_3, W_3)}(1)\cong \C[\lambda]$. Then, the relation \eqref{eq:action alg} becomes
\begin{align*}
\psi(z) \lambda^{j} \psi(z)^{-1}
=&
\frac{
(z-\lambda-  \hbar_1)
(z-\lambda-  \hbar_2)
(z- \lambda-\hbar_3)}
{(z-\lambda+  \hbar_1)
(z-\lambda+  \hbar_2)
(z- \lambda+\hbar_3) }\lambda^{j}\\
=&\frac{
(z-\sigma^+-  \hbar_1)
(z-\sigma^+-  \hbar_2)
(z- \sigma^+-\hbar_3)}
{(z-\sigma^++  \hbar_1)
(z-\sigma^++  \hbar_2)
(z- \sigma^++\hbar_3) }\lambda^{j}\\
=&-\frac{P(\sigma^+, z)}{P(z, \sigma^+)}\lambda^{j}, 
\end{align*}
which is equivalent to 
$
P(z, \sigma^+) \psi(z) \lambda^{j} =-P(\sigma^+, z)\lambda^{j} \psi(z).
$
This completes the proof of Theorem~\ref{thm:positiveYangian}(2). 
$\blacksquare$

\subsection{Review of the $2d$ COHA}\label{subsec:2dCOHA}
Let $J$ be the Jordan quiver. Recall the details of the $2d$ COHA of  $J$ a.k.a {\it preprojective COHA} (see 
\cite{SV1, SV2, YZ1}). In this section, we recall the definition and a few facts of the $2d$ COHA which has been discussed briefly in \S\ref{sec:2.2}. 
Consider the commuting variety
\[
C_n=\{(B_1, B_2) \in \mathfrak{gl}_n^2\mid [B_1, B_2]=0\}. 
\]
The two dimensional torus ${{\bf T}_2}$ acts on $C_n$ by 
\begin{equation}\label{T2action}
(t_1, t_2)\cdot ( B_1, B_2)=(t_1B_1, t_2B_2). 
\end{equation}
This action gives rise to the quotient stack. The $2d$ COHA is the cohomology of 
this stack endowed 
with the natural  Hall multiplication (see \cite{SV2, YZ1}). 
By definition,  the $2d$ COHA of $J$ as a graded vector space is isomorphic to
\[
\H^{(Q_3,W_3),{\bf T}_2}_{B_3=0}:=\bigoplus_{n\in \bbN}H_{\BM}^{\GL_n\times 
{{\bf T}_2}} (C_n), 
\]

\begin{thm}(\cite[Appendix]{ReSo} \cite[Theorem A]{YZ16})\label{thm:CritPrep}
Assume the ${{\bf T}_2}$-action on $C_n$ is given by \eqref{T2action}, and on 
$\mathfrak{gl}_n^3$ by \eqref{T3action}. There is an isomorphism of 
$\bbN$-graded associative algebras $\Xi: \H^{(Q_3,W_3), {\bf T}_2 }_{B_3=0}\to 
\H^{(Q_3, W_3), {\bf T}_2}$ 
 whose restriction to the degree-$n$ piece is
\[
\Xi_v: \H^{(Q_3,W_3), {\bf T}_2}_{B_3=0}(n)\to  \H^{(Q_3, W_3), {\bf T}_2}(n), 
\,\ \,\
f \mapsto (-1)^{{n}\choose{2}} f.
\]\end{thm}
As a consequence, Theorem \ref{thm:critaction} implies the following. 
\begin{cor}
The $2d$  COHA  $\H^{(Q_3,W_3), {\bf T}_2}_{B_3=0}$ acts on $V_{r_1, r_2, r_3}$, for any $(r_1, r_2, r_3)\in \bbN^3$. 
\end{cor}

\begin{rmk} 
Here the notation $H^{(Q_3,W_3),{\bf T}_2}_{B_3=0}$ stresses two facts. First, the $2d$ COHA is the dimensional 
reduction to the plane $B_3=0$
of the $3d$ COHA $\H^{(Q_3,W_3)}$ and second, equivariancy with respect to the 
torus
${\bf T}_2$.

There is no special notation for a $2d$ COHA. In the case of the 
Jordan quiver $J$
 the $2d$ COHA  was denoted by $\textbf{C}'$ in \cite{SV1, SV2}. The notation 
$\calP(H_{\BM}, Q)$ was used in \cite{YZ1} for a general quiver $Q$. Since all 
of those are just special cases of a more general notion of $3d$ COHA 
$\H^{(Q,W)}$),
 we use the notation which also explains how the two are related in the case of 
the quiver $J$.
\end{rmk}

%Let $R=\C[x, y]=H_{\BM}^{{{\bf T}_2}}(\pt)$, and $R_G:=H_{\BM}^{G}(\pt)$. 
The obvious inclusion $C_n\subset \mathfrak{gl}_n\times \mathfrak{gl}_n$ gives 
an $H_{\BM}^{\GL_n\times {{\bf T}_2}}(\pt)$-module homomorphism
\[
\H^{(Q_3,W_3), {{\bf T}_2}}_{B_3=0}(n) \to H_{\BM}^{\GL_n\times {{\bf T}_2}} 
(\mathfrak{gl}_n\times \mathfrak{gl}_n), 
\]
 which gives an isomorphism of the algebra $\bigoplus_{n\ge 0}\H^{(Q_3,W_3), 
{\bf T}_2}_{B_3=0}(n)$ with its image. More precisely, we have a surjective 
algebra homomorphism of $\H^{(Q_3,W_3), {\bf T}_2}_{B_3=0}(n)$  onto its 
image. This map was conjectured \cite[Conjecture 
4.4]{SV2} and later proven \cite[Proposition~4.6]{SV3} to be an isomorphism. We 
will not distinguish these two algebras.
The following description \cite[\S 4.4]{SV2} \cite{YZ1} is used to compare the affine Yangian with the 3d COHA in Appendix \ref{subsec:app_Yangian}. 

The shuffle algebra $\bSh$ is an $\bbN$-graded $\C[\hbar_1, \hbar_2]$-algebra. As a $\C[\hbar_1, \hbar_2]$-module, we 
have 
\[
\bSh=\bigoplus_{n\in\bbN}\bSh_n=\bigoplus_{n\in\bbN} \C[\hbar_1, \hbar_2]
[\lambda_1, \cdots, \lambda_n]^{\fS_n}.\]
For any $n$ and $m\in \bbN$, we consider 
$\bSh_{n}\otimes \bSh_{m}$ as a subspace of 
$\C[\hbar_1, \hbar_2][\lambda_1,\lambda_2, \cdots, \lambda_{n+m}]$ by sending 
$\lambda_s'$ to $\lambda_s$, and $\lambda''_t$ to $\lambda_{t+n}$. Set:
\begin{align}\label{fac}
\fac:
=
\prod_{s=1}^{n}
\prod_{t=1}^{m}
\frac{
(\lambda_s{'}-\lambda_t{'' }-  \hbar_1)
(\lambda_s{'}-\lambda_t{''}-  \hbar_2)
(\lambda'_s- \lambda_t''-\hbar_3)}{\lambda'_s- \lambda_t''}
\end{align}
where $\hbar_1+\hbar_2+\hbar_3=0$. 
Let $\Sh(n, m)$ be the set of shuffles of $n, m$. 
The multiplication of $f_1(\lambda_1, \cdots, \lambda_n)\in \bSh_{n}$ and 
$f_2(\lambda_1, \cdots, \lambda_m)\in \bSh_{m}$ is defined to be
\begin{equation}\label{shuffle formula}
\sum_{\sigma\in\Sh(n, m)}
\sigma(f_1\cdot f_2\cdot \fac)\in \C[\hbar_1, \hbar_2][\lambda_1, \lambda_2, \cdots, 
\lambda_{n+m}]^{\fS_{n+m}}.
\end{equation}
In \cite{SV2} the algebra $S\H^{(Q_3,W_3), {\bf T}_2}_{B_3=0}$ is 
denote by $\textbf{SC}$.
\begin{thm}\cite[Theorem 4.7]{SV2}
\label{thm:shuffle}
There is exists a unique $\C[\hbar_1, \hbar_2]$-algebra embedding $\mathcal{SH}^{(Q_3,W_3), {\bf T}_2}_{B_3=0}\to 
\bSh$, where $S\H^{(Q_3,W_3), {\bf T}_2}_{B_3=0}$ is the spherical subalgebra of 
$\H^{(Q_3,W_3), {\bf T}_2}_{B_3=0}$ generated by the first graded component 
$\H^{(Q_3,W_3), {\bf T}_2}_{B_3=0}(1)$. 
\end{thm}
In Appendix  \ref{subsec:app_Yangian}, we prove 
Theorem~\ref{thm:positiveYangian}(1) by constructing an algebra homomorphism from $\affY$ to $\bSh$. 
This will induce an algebra isomorphism $\affY \cong \mathcal{SH}^{(Q_3,W_3), {\bf T}_2}_{B_3=0}$. The latter is isomorphic to $\mathcal{SH}^{(Q_3, W_3), {\bf T}_2}$ by Theorem \ref{thm:CritPrep}. This proves Theorem~\ref{thm:positiveYangian}(1). 

When $r_a=r_b=0$ and $r_c=r$, such that $\{a, b, c\}=\{1, 2, 3\}$. The action of the positive subalgebra $Y^+$ on $V_r$ is defined via Nakajima raising operators. Theorem~\ref{thm:positiveYangian}(4) follows from Proposition \ref{prop:NakajimaOperators}.

\subsection{Action of the whole algebra $\affY$}
\label{subsec:wholeY}

In this subsection, we show the ``creation operators'' and ``annihilation operators'' of the 3d COHA 
can be glued to give rise to an action of the double COHA 
$D(\mathcal{SH}^{(Q_3,W_3)})$ on $V_{r_1,r_2,r_3}$, for any $r_1, r_2, r_3\in \N$. 
This is shown using the affine Yangian action on the 2d framed quiver varieties and the compatibility result in Proposition 
\ref{prop:coprod}.

On each free boson space $V_{1, 0, 0}, V_{0, 1, 0}, V_{0, 0, 1}$, 
we have constructed an action of the 3d COHA $\mathcal{SH}^{\geq 0} \cong Y^{\geq 0}$. 
By Proposition \ref{prop:NakajimaOperators}, this action of $\mathcal{SH}^{(Q_3, W_3)}$ can be described in terms of Nakajima raising operators. 
Similarly, we could define an action of $\mathcal{SH}^{\leq 0}\cong Y^{\leq 0}$ using the lowering operators. 
More explicitly, as notations in \S\ref{subNak},  let $\beta \in V_{r_1, r_2, r_3}(n)\cong H_{c, GL_{\vec{r}}\times {{\bf T}_2}}( \fM_{\vec{r}}(n), \varphi_{(W^{fr})_{n}})^{\vee}$, the convolution $p_{1*}\big( f(c_1(\calL)  \cdot p_2^*(\beta) )\big) \in V_{r_1, r_2, r_3}(n-1)$ gives an action of $f(\lambda)\in \mathcal{SH}^{(Q_3, W_3)}(1)
\cong \C[\hbar_1,\hbar_2,\lambda] $ on $V_{r_1, r_2, r_3}$, for any $r_1, r_2, r_3\in \N$.

Theorem~\ref{thm:positiveYangian}, together with Proposition \ref{prop:coprod}  implies the following.
Let $Y^{\geq 0}\otimes Y^{\leq 0}$ be the free non-commutative tensor product of the algebras $Y^{\geq 0}$ and $Y^{\leq 0}$. 
\begin{cor}
The action of $Y^{\geq 0}\otimes Y^{\leq0}$ on $V_{r_1,r_2,r_3}$ factors through $Y_{\hbar_1, \hbar_2, \hbar_3}(\widehat{\fg\fl_1})$. 
\[\xymatrix{
Y^{\geq 0}\otimes Y^{\leq0}\ar[r] \ar[d] &\End(V_{r_1,r_2,r_3})\\
Y_{\hbar_1, \hbar_2, \hbar_3}(\widehat{\fg\fl_1}) \ar[ru]&
}\]
\end{cor}
{\it Proof.}
It is known that (see, e.g., \cite[Theorem 2.2]{T}) the action of $Y^{\geq 0}\otimes Y^{\leq0}$ on each free boson space $V_{1, 0, 0}, V_{0, 1, 0}, V_{0, 0, 1}$ factors through $\affY$, and the coproduct $\Delta^{\Dr}: Y^{\geq 0}\otimes Y^{\leq0}\to \affY \otimes \affY$ factors through $\affY$.
Up to localization, the map (see Definition \ref{mapl})
\[
l^{\otimes |\vec{r}|}:V_{r_1, r_2, r_3} \to (V_{1, 0, 0})^{\otimes r_1} \otimes (V_{0, 1, 0})^{\otimes r_2}\otimes (V_{0, 0, 1})^{\otimes r_3}
\]
 is an isomorphism. Therefore, the action of $Y^{\geq 0}\otimes Y^{\leq0}$ on $V_{r_1, r_2, r_3}$ also factors through $\affY$.
 This completes the proof. 
$\blacksquare$

As a consequence, the actions of $\mathcal{SH}^{\geq 0}$ and $\mathcal{SH}^{\leq 0}$ on $V_{r_1, r_2, r_3}$ 
satisfy the commutation relations of the Drinfeld double $D(\mathcal{SH}^{(Q_3,W_3)}) \cong \affY$.
Thus, we have an action of 
$D(\mathcal{SH}^{(Q_3,W_3)})$ on  $V_{r_1,r_2,r_3}$, for any $r_1, r_2, r_3\in \N$.

\begin{rmk}
Closely related to the affine Yangian, 
there is an algebra called the deformed double current algebra, 
which deforms the enveloping algebra of the double current algebra $\mathfrak{gl}_N\otimes \C[u, v]$.  
Let $D_{\epsilon, \delta}(\mathfrak{gl}_N)$ denote the deformed double current algebra.
It has been conjectured by K. Costello  \cite[Section 1.8]{Costello} that
$D_{\epsilon, \delta}(\mathfrak{gl}_N)$ should 
act on the equivariant Donaldson-Thomas theory of $\C\times \widetilde{\C^2/\Z_N}$, where $\widetilde{\C^2/\Z_N}$ is the minimal resolution of ${\C^2/\Z_N}$. 
Furthermore, it is conjectured that the 3d COHA associated to the $3$-fold $\C\times \widetilde{\C^2/\Z_N}$ embeds into $D_{\epsilon, \delta}(\mathfrak{gl}_N)$, such that the action of $D_{\epsilon, \delta}(\mathfrak{gl}_N)$ restricts to the natural action of the 3d COHA on the 
same space. From the physics perspective, \cite{Costello} probes the configuration of a single $D6$-brane wrapping the whole $\C\times \widetilde{\C^2/\Z_N}$. 
This configuration is expected to lead to a "vacuum" module of the involved affine Yangian (or the deformed double current algebra).

In particular, when $N=1$, the expectation of  \cite[Section 1.8]{Costello} is 
that $D_{\epsilon, \delta}(\mathfrak{gl}_1)$ acts on $\bigoplus_{n\in \bbN} H^*_c( 
\Hilb_n(\C^3), \varphi_{W_3^{fr}})^{\vee}$. Note that $\Hilb_n(\C^3)$ is obtained from the representations $\calM_{0, 0, 1}(n)$ of the quiver with potential $(Q_3^{fr}, W_3^{fr})$ using the 
stability condition \eqref{eq:stabN} without imposing the relations
\[
B_3 I_{12}=0,  B_1 I_{23}=0,  B_2 I_{13}=0, 
\]
i.e. considering the moduli space associated to a single $D6$-brane.  This configuration is expected to lead to the MacMahon module of the affine Yangian of $\mathfrak{gl}(1)$
or the corresponding deformed double current algebra. We remark that the space $\bigoplus_{n\in \bbN} H^*_c( 
\Hilb_n(\C^3), \varphi_{W_3^{fr}})^{\vee}$ is different from the space $V_{r_1, r_2, r_3}$ considered in the present paper. 
\end{rmk}

%As  $D_{\epsilon, \delta}(\mathfrak{gl}_1)$ and $\affY$ are conjecturally two 
%different doubles of  $\affY^+$, therefore, we do not expect $\affY$ acts on 
%$\bigoplus_{n\in \bbN} H^*( \Hilb_n(\C^3), \varphi_{W_3^{fr}})^{\vee}$. 
%\footnote{Using restriction to open subset, there is a surjective map 
%$\bigoplus_{n\in \bbN} H^*( \Hilb_n(\C^3), \varphi_{W_3^{fr}})^{\vee}\to V_{1, 0, 0}$. The latter 
%$V_{1, 0, 0}$ is known to be Verma module. 
% Were this map of representations, it would have a splitting, which sounds 
%unlikely. }

\section{Free field realization}
\label{sec:free-field}
In this section, we first introduce a central extension 
$\SH^{\vec{c}}$ of the affine Yangian by a polynomial ring 
$\C[c_l^{(1)}, c_l^{(2)}, c_l^{(3)}: l\geq 0]$. 
This algebra $\SH^{\vec{c}}$ also acts on $V_{r_1, r_2, r_3}$
 such that the action of the central elements $\vec{c}=\{c_l^{(1)}, c_l^{(2)}, c_l^{(3)}: l\geq 0\}$  depends on the framing vector $\vec{r}=(r_1, r_2, r_3)$. We construct a coproduct $\Delta^{\vec{c}}$ 
 (different from the Drinfeld coproduct) on  $\SH^{\vec{c}}$.  When $c_l^{(2)}=c_l^{(3)}=0$, the algebra $\SH^{(\bold{c^{(1)}}, \bold{0}, \bold{0})}$ and the coproduct $\Delta^{(\bold{c^{(1)}}, \bold{0}, \bold{0})}$ coincide with the algebra $\SH^{\bold{c}}$ and coproduct $\bold{\Delta}$ introduced by Schiffmann--Vasserot \cite{SV2}. 

Recall the quotient $\calM_{\vec{r}}(n)^{st}/\GL_{n}$ is denoted by $\fM_{\vec{r}}(n)$. 
Geometrically, for any 1-dimensional subgroup $A$ of the maximal torus of 
$\GL_{r_1}\times\GL_{r_2}\times\GL_{r_3}$, the fixed point locus $(\fM_{\vec{r}})^A$ is a product of framed quiver varieties of the form $\prod_{k}\fM_{\vec{r^k}}$ with $\sum_k \vec{r^k}=\vec{r}$. 

For any decomposition $\vec{r}=\vec{r'}+\vec{r''}$ and a choice of an appropriate torus $A$ in $\GL_{r_1}\times\GL_{r_2}\times\GL_{r_3}$, the fixed point locus $(\fM_{\vec{r}})^A$  is isomorphic to $\fM_{\vec{r'}}\times \fM_{\vec{r''}}$.
In \S~\ref{subsec:hypLoc} we will recall the hyperbolic restriction functor and show that it commutes with the vanishing cycle functor (cf. with \cite[Proposition~5.4.1.2]{NakPCMI}). Therefore, we have the hyperbolic restriction  
$
h: V_{\vec{r}}  \to V_{\vec{r'}}\otimes 
V_{\vec{r''}}.
$

Taking sufficiently general subtorus $A\subset \GL_{r_1}\times\GL_{r_2}\times\GL_{r_3}$  we get a hyperbolic restriction
\[
h: V_{r_1,r_2,r_3}\to (V_{1, 0, 0})^{\otimes r_1} \otimes (V_{0, 1, 0})^{\otimes r_2}\otimes (V_{0, 0, 1})^{\otimes r_3}, \] 
where the target is the Fock space of $(r_1+r_2+r_3)$-free bosons. 

We show the actions of $\SH^{\vec{c}}$ on $V_{r_1,r_2,r_3}$ and on $(V_{1, 0, 0})^{\otimes r_1} \otimes (V_{0, 1, 0})^{\otimes r_2}\otimes (V_{0, 0, 1})^{\otimes r_3}$ (via the coproduct $\Delta^{\vec{c}}$) are compatible with the hyperbolic restriction map.  

%However, although $V_{1,0,0}$, $V_{0,1,0}$ and $V_{0,0,1}$
%are abstractly isomorphic to $L^{(1)}$, they are different as representations of 
%$\H^{(Q_3, W_3)}$. 
%Indeed, the representation $V_{0,0,1}$ is obtained from the 
%composition of the isomorphism in Theorem~\ref{thm:CritPrep} and the 
%representation from \cite{SV2}. 
%On the other hand, in order to get the 
%representation $V_{1,0,0}$, one needs to modify the isomorphism in 
%Theorem~\ref{thm:CritPrep} by replacing $B_1, B_2$ by $B_2,B_3$. In particular, 
%the deformation parameters $\hbar_1,\hbar_2$ are replaced by $\hbar_2,\hbar_3$.

\subsection{A central extension}\label{subsec:central}

Now we introduce the  algebra $\SH^{\vec{c}}$. This is similar to its 2d counterpart from \cite{SV2}, so we refer to the loc. cit. for the terminology and some details.

Denote by $\Lambda_{F}$ the ring of symmetric polynomials in infinitely many variables with coefficients over the field of rational functions 
$F=\C(\hbar_1, \hbar_2)$. It has a basis given by $\{J_{\lambda} \mid \lambda \,\ \text{is a partition}\}$, where $J_\lambda$ is the integral form of the Jack polynomial associated to $\lambda$, and to the parameter $-\hbar_1/\hbar_2$. 

Let  $D_{l, 0}\in \End(\Lambda_F)$ be the operator of multiplication by the standard symmetric power-sum function $p_l$. 
Let $D_{0, l}\in \End(\Lambda_F) $ be the Sekiguchi operators (see loc.cit.). The Sekiguchi operators form a collection of commuting differential operators whose joint spectrum is given by the Jack polynomials $\{J_{\lambda}\mid  \lambda \,\ \text{is a partition}\}$.

Let $\SH^{\geq 0}$ be the unital subalgebra of $\End(\Lambda_F)$ generated by 
$\{D_{0, l}, D_{l, 0} \mid l\geq 1\}$. For $l\geq 0$, we set $D_{1, l}=[D_{0, l+1}, D_{1, 0}]$.
Let $\SH^{>}$ be the unital subalgebra of $\SH^{\geq 0}$ generated by 
$\{D_{1, l}\mid l\geq 0\}$, and $\SH^0$ the subalgebra of  $\SH^{\geq 0}$ generated by $\{D_{0, l}\mid l\geq 1\}$. 
It is known that $\SH^0=F[D_{0, 1}, D_{0, 2}, \cdots]$ is the polynomial algebra and $\SH^{\geq 0}=\SH^{>}\otimes \SH^0$ (see \cite[Section 1.7]{SV2}). 

Let $\SH^{<}$ be the opposite algebra of $\SH^{>}$. Its generator corresponding to $D_{1, l}$ is denoted by $D_{-1, l}$.  
The algebra $\SH$ is generated by $\SH^{>}, \SH^{0}, \SH^<$ 
modulo a certain set of relations involving the commutators of elements from $\SH^{<}$ and $\SH^{>}$. In particular, $\SH$ has a triangular decomposition $\SH=\SH^{>}\otimes \SH^{0}\otimes \SH^{<}$. 

In the setup of the present paper, we use the following normalization of generators of $\SH^{0}, \SH^{>}, \SH^{<}$ by scaling the generators $D_{k, l}$ via appropriate powers of $\hbar_1, \hbar_2$ according to their geometric actions on the quiver varieties: 
%\begin{equation}
%D_{0, l+1}=\hbar_1^{-l}f_{0, l}, \,\ D_{1, l}=\hbar_1^{1-l}\hbar_2 f_{1, l} \text{ and } \,\ D_{-1, l}=\hbar_1^{-l} f_{-1, l}. 
%\end{equation}
\begin{equation}
f_{0, l}:=\hbar_1^{l} D_{0, l+1}, \,\  f_{1, l}:=\hbar_1^{l-1}\hbar_2^{-1}D_{1, l} \text{ and } \,\ f_{-1, l}:=\hbar_1^{l} D_{-1, l}. 
\end{equation}

Let $\GL_{r_1}\times \GL_{r_2} \times \GL_{r_3}$ be the general linear group at the framing vertices of $Q_3^{fr}$. 
Denote by $T_{r_1,r_2,r_3}\subset \GL_{r_1}\times \GL_{r_2} \times \GL_{r_3}$ its maximal torus.  
Then, 
\[
H_{T_{r_1,r_2,r_3} \times {{\bf T}_2}}(\pt)=
\C[\hbar_1, \hbar_2, \mu_1^{(1)}, \cdots, \mu_{r_1}^{(1)}, \mu_1^{(2)}, \cdots, \mu_{r_2}^{(2)}, \mu_1^{(3)}, \cdots, \mu_{r_3}^{(3)}
].
\] 
Let $K_{r_1, r_2. r_3}=\C(\hbar_1, \hbar_2, \mu_1^{(1)}, \cdots, \mu_{r_1}^{(1)}, \mu_1^{(2)}, \cdots, \mu_{r_2}^{(2)}, \mu_1^{(3)}, \cdots, \mu_{r_3}^{(3)})$ be the fractional field of $H_{T_{r_1,r_2,r_3} \times {{\bf T}_2}}(\pt)$.

We introduce a new family $\vec{c}:=(\bold{c}^{(1)}, \bold{c}^{(2)}, \bold{c}^{(3)}) =(c_0^{(1)}, c_0^{(2)}, c_0^{(3)}, c_1^{(1)},  c_1^{(2)}, c_1^{(3)}, \cdots )$ of formal parameters, which roughly speaking comes from 
$
\lim_{\{r_1, r_2, r_3\to \infty\}} H_{T_{r_1,r_2,r_3}}(\pt).  
$
Now set
\[
F^{\vec{c}}=\C(\hbar_1, \hbar_2)[c_l^{(1)}, c_l^{(2)}, c_l^{(3)}: l\geq 0], \,\
\SH^{\vec{c}, 0}=F^{\vec{c}}[f_{0, l}: l \geq 0]. 
\]

\begin{defn}
Let $\SH^{\vec{c}}$ be the $F$-algebra generated by $\SH^{>}$, $\SH^{\vec{c}, 0}$, and $\SH^{<}$, modulo 
the following set the relations
\footnote{The central elements $\{\bold{c^{(1)}_i}, i\geq 1\}$ of $\SH^{\vec{c}}$ correspond to the central elements 
$\{\bold{c_i}\hbar_1^i, i\geq 1\}$ of $\SH^{\bold{c}}$ in \cite{SV2}, and  $\bold{c^{(1)}_0}$ of 
$\SH^{\vec{c}}$ corresponds to  $\frac{\bold{c_0}}{\hbar_1}$ of $\SH^{\bold{c}}$ in \cite{SV2}.}
\begin{align*}
& \text{$c_l^{(1)}, c_l^{(2)}, c_l^{(3)}$ are central,} \,\ l\geq 0 \\
& [f_{0, l}, f_{1, k}]=f_{1, l+k},  l\geq 1, \\
&[f_{0, l}, f_{-1, k}]=-f_{-1, l+k}, l\geq 1, \\
&  [f_{1, l}, f_{-1, k}]=G_{l+k},\,\  l, k\geq 0, 
\end{align*}
where $f_{0, 0}=0$, and the elements $G_{k+l}$ are determined through the formula
\begin{align*}
1-\hbar_1\hbar_2\hbar_3 \sum_{l\geq 0} G_{l} t^{l+1}
=&
\exp(\sum_{l\geq 0} (-1)^{l+1}c_l^{(1)} \phi_{l}(t, \hbar_1)) 
\exp(\sum_{l\geq 0} (-1)^{l+1}c_l^{(2)} \phi_{l}(t, \hbar_2))\\&
\exp(\sum_{l\geq 0} (-1)^{l+1}c_l^{(3)} \phi_{l}(t, \hbar_3))  
\exp(\sum_{l\geq 0} f_{0, l} \varPhi_l(t)). 
\end{align*}
Here the functions $\phi_{l}(t, \hbar_i)$ and $\varPhi_l(t)$ are formal power series in $t$ depending on $\hbar_1, \hbar_2, \hbar_3$. 
They are given by the following formulas. 
\begin{align*}
&\exp( \sum_{l\geq 0}(-1)^{l+1} a^{l}  \phi_{l}(t, \hbar_i))
=\frac{1+t(a-\hbar_i)}{1+ta}, \\
&\exp( \sum_{l\geq 0}(-1)^{l+1} a^{l} \varPhi_l(t))
=\frac
{(1+t(a-\hbar_1))
(1+t(a-\hbar_2))
(1+t(a-\hbar_3)
)}{(1+t(a+\hbar_1))(1+t(a+\hbar_2))(1+t(a+\hbar_3))}, 
\end{align*}
where $a$ is any element in $\C$, and we impose the condition that $\hbar_1+\hbar_2+\hbar_3=0$. \end{defn}

The algebra $\SH^{\bold{c}}$ from \cite{SV2}, which  is generated by $\SH^{>}, \SH^{0}, \SH^{<}$ together with a family of central elements $\bold{c}=\{c_0, c_1, c_2\cdots\}$, modulo a certain set of relations involving the commutators $[D_{-1,k},D_{1,l}]$ \cite{AS, SV2}, is a specialization of the algebra $\SH^{\vec{c}}$ above with $c^{(2)}_l=c^{(3)}_l=0$

Similarly to \cite{SV2}, for any $\vec{r}=(r_1, r_2, r_3)$, 
there is an algebra $\SH^{(r_1, r_2, r_3)}$ obtained by the specialization of 
$\SH^{\vec{c}}\otimes_{F^{\vec{c}}} K_{r_1, r_2, r_3}$ with
\begin{align*}
&c_0^{(1)}=r_1, \,\ c_0^{(2)}=r_2, \,\ c_0^{(3)}=r_3, \\
&c_i^{(k)}=p_i(\mu_1^{(k)}, \mu_2^{(k)}, \cdots, \mu_{r_k}^{(k)}), k=1, 2, 3. 
\end{align*}
This specialization is determined using the following formula
\begin{align*}
&
\exp( \sum_{l\geq 0} (-1)^{l+1} p_l(\mu_a^{(1)}) \phi_{l}(t, \hbar_1))
\exp( \sum_{l\geq 0} (-1)^{l+1} p_l(\mu_a^{(2)}) \phi_{l}(t, \hbar_2))\\
&\cdot\exp( \sum_{l\geq 0} (-1)^{l+1} p_l(\mu_a^{(3)}) \phi_{l}(t, \hbar_3))\\
&=
\prod_{a=1}^{r_1}\frac{z-\mu_a^{(1)}+\hbar_1}{z-\mu_a^{(1)}}
\prod_{a=1}^{r_2}\frac{z-\mu_a^{(2)}+\hbar_2}{z-\mu_a^{(2)}}
\prod_{a=1}^{r_3}\frac{z-\mu_a^{(3)}+\hbar_3}{z-\mu_a^{(3)}}, 
\,\ \text{where $z=-1/t$, }
\end{align*}
and $\mu_1^{(k)}, \cdots, \mu_{r_k}^{(k)}$ are free variables. In the action in the following proposition, these variables are specialized to be the Chern roots of $\calE_{r_k}$, $k=1, 2, 3$.

\begin{prp}\label{prop:shuffle specialization}
\leavevmode
\begin{enumerate}
\item
There is an algebra homomorphism $\SH^{\vec{c}} \to \affY$, given by
\[
f_{1, l} \mapsto e_l,\,\  f_{-1, l}\mapsto f_l, \,\  G_{l} \mapsto \psi_{l}, \,\  \bold{c_i}^{(k)} \mapsto 0, \,\ k=1, 2, 3. 
\]
\item
The algebra $\SH^{\vec{c}}$ acts on $V_{r_1, r_2, r_3}$ through the action of the specialization $\SH^{(r_1, r_2, r_3)}$ on $V_{r_1, r_2, r_3}$, for any $\vec{r}=(r_1, r_2, r_3)\in \Z^3_{\geq 0}$. 
\end{enumerate}
\end{prp}
{\it Proof.}
For (1): By \cite[Corollary 6.4]{SV2}, there is an algebra isomorphism of $\SH^{>}$ and the spherical 2d COHA 
$\mathcal{SH}^{(Q_3,W_3),{\bf T}_2}_{B_3=0}$. The latter is in turn isomorphic to $Y^+$ by Theorem \ref{thm:positiveYangian}. 
Thus, we have the isomorphisms of algebras
\[
\SH^{>}\cong Y^{+}, f_{1, l} \mapsto e_l \,\ \text{and} \,\  \SH^{<}\cong Y^{-}, f_{-1, l} \mapsto f_l. 
\]
Define a natural map $\SH^{\vec{c}, 0}\to Y^0$ by $G_{l} \mapsto \psi_{l}, \,\  \bold{c_i}^{(k)} \mapsto 0$. 
It is straightforward to check that the above assignments define an algebra homomorphism $\SH^{\vec{c}} \to \affY$. 

For (2): Using the isomorphisms $\SH^{>}\cong Y^{+}, \SH^{<}\cong Y^{+}$, we deduce that the algebras $\SH^{>}$ and  $\SH^{<}$ act on $V_{r_1, r_2, r_3}$. We require the central elements $\{c_i^{(k)}\mid i\in \Z_{\geq 0}, k=1, 2, 3\}$ acts by the formula 
$p_i(\mu_1^{(k)}, \mu_2^{(k)}, \cdots, \mu_{r_k}^{(k)})$. 
The series $1-\hbar_1\hbar_2\hbar_3 \sum_{l\geq 0} G_{l} t^{l+1}$ acts by $\lambda_{t}(\calF(\calV_{n}, \calE_{\vec{r}}))$, which is given by the formula \eqref{action of Y0}. In particular, the action of $f_{0, l}\in \SH^{\vec{c}, 0}$ is determined by the formula 
\begin{align*}
\exp(\sum_{l\geq 0} f_{0, l} \varPhi_l(t))\cdot m=\left(\prod_{d=1}^{n}\frac{z-\lambda_d-\hbar_1}{z-\lambda_d+\hbar_1}
\frac{z-\lambda_d-\hbar_2}{z-\lambda_d+\hbar_2}
\frac{z-\lambda_d-\hbar_3}{z-\lambda_d+\hbar_3}\right)^-  \cdot m, 
\end{align*} for $m\in V_{r_1, r_2, r_3}$, $z=-1/t$, and $\lambda_1, \cdots, \lambda_d$ are the Chern roots of the tautological bundle $\calV$. 
It is straightforward to check that the above assignments give an action of  $\SH^{\vec{c}}$ on $V_{r_1, r_2, r_3}$. 
The action obviously factors through the algebra $\SH^{(r_1, r_2, r_3)}$. 
$\blacksquare$

\begin{rmk}
\begin{enumerate}
\item
The action of $\SH^{\vec{c}}$ on $V_{r_1, r_2, r_3}$ does not factor through the quotient map $\SH^{\vec{c}} \to \affY$, as the actions of the central elements $\bold{c_i}^{(k)}$ are not trivial. 
Instead, the action of $\SH^{\vec{c}}$ factors through the specialization $\SH^{(r_1, r_2, r_3)}$ via $c_i^{(k)}=p_i(\mu_1^{(k)}, \mu_2^{(k)}, \cdots, \mu_{r_k}^{(k)})$. 
\item One could send the central elements $\bold{c_i}^{(k)}$ to nonzero central elements in $\affY$. Then, using the identity
\begin{align*}
1-\hbar_1\hbar_2\hbar_3 \sum_{l\geq 0} \psi_{l} t^{l+1}=&
\exp(\sum_{l\geq 0} (-1)^{l+1}c_l^{(1)} \phi_{l}(t, \hbar_1)) 
\exp(\sum_{l\geq 0} (-1)^{l+1}c_l^{(2)} \phi_{l}(t, \hbar_2))\\&
\exp(\sum_{l\geq 0} (-1)^{l+1}c_l^{(3)} \phi_{l}(t, \hbar_3))  
\exp(\sum_{l\geq 0} f_{0, l} \varPhi_l(t)), 
\end{align*}
we could solve the generators $\{f_{0, l} \mid l\geq 0\}$ of $\SH^{\vec{c}}$ 
in terms of $\{\psi_{l} \mid l\geq 0 \}\subset Y^0$, which would give the image of $f_{0, l}$ in $Y^0\subset \affY$. The algebra homomorphism $\SH^{\vec{c}}\to \affY$ defined this way is compatible with their actions on $V_{\vec{r}}$.
\end{enumerate}
\end{rmk}

\subsection{The coproduct of $\SH^{\vec{c}}$}
\label{sec:twisted coproduct}

Introduce the following elements of $\SH^{\vec{c}}$
\begin{align*}
&B_{-l}:=\frac{\ad^{l-1} }{(l-1)!}( f_{1, 1})( f_{1, 0})=\frac{1}{(l-1)!}[f_{1, 1}, [f_{1, 1}, \cdots, [f_{1, 1}, f_{1, 0}]]], 
\\
&B_l:=\frac{(  f_{-1, 0})\ad^{l-1}( f_{-1, 1})}{(l-1)!}=\frac{1}{(l-1)!} [[[[f_{-1, 0}, f_{-1, 1}],  f_{-1, 1}],  \cdots, f_{-1, 1}], f_{-1, 1}], \\
&B_0:=G_1= [f_{1, 1}, f_{-1, 0}], \,\ G_0= [f_{1, 0}, f_{-1, 0}]. 
\end{align*}
The following lemma is analogous to \cite[Proposition 1.40]{SV2}, and the proof is similar. 
\footnote{
When $\bold{c^{(1)}, c^{(2)}}=0$, in \cite{SV2}, the Heisenberg subalgebra of $\SH^{\bold{c}}$ is generated by 
$\{b_{-l}, b_l, b_0, E_0\mid l\geq 0\}$. To compare with the notation in the current paper, we have 
\[
B_{-l}=\frac{b_{-l}}{\hbar_1}, \,\ B_l=\frac{b_l}{\hbar_1},  \,\ G_0=\frac{E_0}{\hbar_1 \hbar_2}=\frac{\bold{c}_0^{(3)}}{\hbar_1 \hbar_2}, \,\
B_0:=G_1=\frac{E_1}{\hbar_2}=-\frac{b_0}{\hbar_1}. 
\]}
\begin{lmm}\label{Heisen}
The elements $\{B_{-l}, B_l, B_0, G_0\mid l\geq 0\}$ form a Heisenberg subalgebra of $\SH^{\vec{c}}$. 
\end{lmm}
By a direct computation (see \cite[Remark 1.41]{SV2}), we have $G_0=\frac{\hbar_1 \bold{c}_0^{(1)}+\hbar_2 \bold{c}_0^{(2)}+\hbar_3 \bold{c}_0^{(3)}}{\hbar_1\hbar_2\hbar_3}$.

\begin{prp}\label{prop:coprod_formula}
The algebra $\SH^{\vec{c}}$ is equipped with a Hopf algebra structure. The coproduct is uniquely determined by the following 
formulas. 
\begin{align*}
&\bold{\Delta}^{\vec{c}}(\bold{c}_l^{(k)})=\bold{c}_l^{(k)}\otimes 1+1\otimes \bold{c}_l^{(k)},   \,\  l\geq 0, k=1, 2, 3\\
&\bold{\Delta}^{\vec{c}}(B_l)=B_l\otimes 1+1\otimes B_l, \,\  l\neq 0\\
&\bold{\Delta}^{\vec{c}}(B_0)=B_0\otimes 1+1\otimes B_0
+\hbar_1 \bold{c}_0^{(1)}\otimes \bold{c}_0^{(1)}
+\hbar_2 \bold{c}_0^{(2)}\otimes \bold{c}_0^{(2)}
+\hbar_3 \bold{c}_0^{(3)}\otimes \bold{c}_0^{(3)}, \\
&\bold{\Delta}^{\vec{c}} (f_{0, 1})= f_{0, 1} \otimes 1+1\otimes f_{0, 1}+(\hbar_1\hbar_2\hbar_3) \sum_{l\geq 1} l 
B_{l}\otimes B_{-l}
\end{align*}
\end{prp}
The proof of Proposition \ref{prop:coprod_formula} will be given in Appendix~\ref{subsec:proof_hyperLoc}.

When $\bold{c}^{2}=\bold{c}^{3}=\bold{0}$, $\SH^{\vec{c}}$ becomes the algebra $\SH^{\bold{c}}$ in \cite{SV2}, 
and the above coproduct specializes to the coproduct in \cite[Theorem 7.9]{SV2}. 

In \S~\ref{subsec:hypLoc} we construct a map of vector spaces 
$h: V_{\vec{r}}\to V_{\vec{r'}}\otimes V_{\vec{r''}}$ via the hyperbolic localization, 
which is compatible with the coproduct $\Delta^{\vec{c}}$. (see Proposition \ref{prp:coprod_hyperbolic}).

\subsection{Hyperbolic localization}\label{subsec:hypLoc}
Fix $\vec{r}=\vec{r'}+\vec{r''}$, and write $\C^{\vec{r}}=\C^{\vec{r'}}\oplus\C^{\vec{r''}}$.
Let $A\subseteq \GL_{\vec{r}}$ be a 1-dimensional subtorus with coordinate $t$ so that the eigenvalues on $\C^{\vec{r'}}$ are $t$ and the eigenvalues on $\C^{\vec{r''}}$ are 1.
Then, taking direct sum of framed representations defines an isomorphism $\fM_{\vec{r'}}\times\fM_{\vec{r''}}\cong \fM_{\vec{r}}^A$, with their tautological bundles satisfying $\calV_1\oplus\calV_2\cong\calV$ and $\calE_{\vec{r'}}\oplus\calE_{\vec{r''}}\cong \calE_{\vec{r}}$. Moreover, the potential function $\tr W$ on  $\fM_{\vec{r}}$, when restricted to the fixed point set $\fM_{\vec{r}}^A$, agrees with the potential function $\tr W\boxplus \tr W$ on the product $\fM_{\vec{r'}}\times\fM_{\vec{r''}}$ under this isomorphism. 

The attracting set $\calA_{\fM_{\vec{r}}}$, consisting of the subset of  $\fM_{\vec{r}}$ that has a limit as $t\to 0$, can be identified as a quotient of the total space of the affine bundle $\calH om(\calV_1,\calV_2)^3\oplus \calH om(\calE_{\vec{r'}},\calV_2)\oplus \calH om(\calV_1,\calE_{\vec{r''}})$. In particular, the potential function $\tr W$ of $\fM_{\vec{r}}$, when restricted to $\calA_{\fM_{\vec{r}}}$, is equal to the pull-back of $\tr W\boxplus \tr W$ on $\fM_{\vec{r}}^A$, as the trace of an endomorphism only depends on the diagonal.

We have the natural correspondence \[\xymatrix{\fM_{\vec{r}}^A&\calA_{\fM_{\vec{r}}}\ar[l]_{p}\ar[r]^{\eta}&\fM_{\vec{r}}}.\]

Recall that the variety $\fM_{\vec{r}}(n)$ is smooth.

\begin{lmm}\label{lmm:smooth}
The variety $\calA_{\fM_{\vec{r}}(n)}$ is an affine bundle on $\fM_{\vec{r}}(n)^A$, whose rank is equal to the codimension of $\calA_{\fM_{\vec{r}}(n)}$ in $\fM_{\vec{r}}(n)$. 
\end{lmm}
{\it Proof.}
Let $g:\C^n\to \C^n$ be an automorphism of  $(B_i,I_{jk},J_{jk})$, hence commutes with $B_i$ and $I_{jk}$. In particular, it is closed on the subspaces $V_l$ for $l=1,2,3$. By definition of $V_1$, it is also closed under $B_2, B_3$, hence, we have a representation $(B_2,B_3,I_{23})$ on the vector space $V_1$ of the standard framed double loop quiver which is furthermore stable in the usual sense (e.g., in \cite[\S~3.2]{Ginz}). Therefore, by the standard results about  stability conditions (see e.g., \cite[Lemma 3.2.3.]{Ginz}) $g|_{V_1}$ is trivial. Similarly, $g|_{V_2}$ and $g|_{V_3}$ are both trivial. As $\C^n=V_1+V_2+V_3$, this implies that $g$ is the identity on $\C^n$.
$\blacksquare$

\begin{rmk}
Note that the $\GL_n$ action on a stable framed representation of a quiver does not factor through the quotient
$\GL_n/\C^*$ by the scalar matrices. Indeed, an automorphism of a framed representation becomes the identity map when restricting on the framing, which rules out scalars as automorphisms. 
\end{rmk}

\begin{prp}\label{prop:hyper_def}
We have a natural isomorphism $p_*\circ\eta^!\circ\varphi_{W}\C_{\fM_{\vec{r}}}\cong \varphi_{W}\C_{\fM_{\vec{r}}^A}[2d]$, where $d=\codim \fM_{\vec{r}}^A$.
\end{prp}
The proof essentially follows the same argument as in \cite[pp.29-31]{NakPCMI}. For completeness, we include the details.
\begin{lmm}\label{lem:hyper_phi}
The map $p_*\eta^!$ intertwines the vanishing cycle functors on $\fM_{\vec{r}}(n)$ and $\fM_{\vec{r}}(n)^A$
\[\varphi_{W} p_*\eta^!\cong p_*\eta^!\varphi_{W}.\]
Moreover, $p_*\eta^!\C_{\fM_{\vec{r}}(n)}[\dim \fM_{\vec{r}}(n)]\cong \C_{\fM_{\vec{r}}(n)^A}[\dim \fM_{\vec{r}}(n)^A]$.
\end{lmm}
{\it Proof.}
For simplicity, here in the proof we denote $\fM_{\vec{r}}(n)$ by $\fX$, with the potential function $f=\tr W:\fX\to \C$. It restricts to functions on $\calA_{\fX}$ and $\fX^A$ denoted by $f_{\calA}$ and $f_{A}$, with the zero locus denoted by $X,\calA_X$ and $X^A$ respectively. We have the diagram 

\begin{equation}\label{eqn:nearby}
\xymatrix{
\fX^A & \calA_{\fX} \ar[l]^{p_{\fX}} \ar[r]_{j_{\fX}} & \fX \\
X^A\ar[u]^{r_A}&\calA_X\ar[u]^{r_\calA}\ar[l]^{p_X}\ar[r]_{j_X}&X\ar[u]^{r_X}\\
}
\end{equation}

We first show that the hyperbolic localization commutes with the nearby cycle functor \begin{equation}\label{eqn:Psi_hyp}
p_{X*}j_{X}^!\Psi_{f}=\Psi_{f_A}p_{\fX^*}j_{\fX}^!.
\end{equation}
Let $c:\tilde{\C}^*\to \C$ be the universal cover of $\C^*\subseteq \C$, and we have the fiber diagram \[\xymatrix{
\tilde{\fX}\ar[r]\ar[d]_{\tilde{c_X}}&\tilde{\C}^*\ar[d]_c\\
\fX\ar[r]&\C.
}\]
Similarly we have $\tilde{c_{\calA}}:\tilde{\calA_{\fX}}\to \calA_{\fX}$ and $\tilde{c_{A}}:\tilde{\fX^A}\to \fX^A$.
Recall that $\Psi_f=r_X^*\tilde{c_X}_{*}\tilde{c_X}^{*}$, and similar for $\Psi_{f_A}$.
We need to show that 
\begin{enumerate}
\item $p_{\tilde{\fX}*}j_{\tilde{\fX}}^!\tilde{c_X}^*=\tilde{c_A}^*p_{\fX*}j_{\fX}^!$;
\item $\tilde{c_A}_*p_{\tilde{\fX}*}j_{\tilde{\fX}}^!=p_{\fX*}j_{\fX}^!\tilde{c_X}_*$;
\item $p_{X*}j_X^!r_X^*=r_{A}^*p_{\fX*}j_{\fX}^!$.
\end{enumerate}
By Lemma~\ref{lmm:smooth},  $p_{\tilde{\fX}*}j_{\tilde{\fX}}^!\tilde{c_X}^*$ and $p_{\tilde{\fX}!}j_{\tilde{\fX}}^*\tilde{c_X}^*$ differ only by a homological degree shifting. We also have 
\begin{align*}
p_{\tilde{\fX}!}j_{\tilde{\fX}}^*\tilde{c_X}^*=p_{\tilde{\fX}!}\tilde{c_{\calA}}^*j_{\fX}^!=\tilde{c_{A}}^*p_{\fX!}j_{\fX}^!,
\end{align*}
which in turn differs from $\tilde{c_A}^*p_{\fX*}j_{\fX}^!$ by the same shifting as above. This proves (1). The other two are proved in a similar way: (3) is proven using diagram \eqref{eqn:nearby}, and replacing $p_{\tilde{\fX}*}j_{\tilde{\fX}}^!$ by $p_{\tilde{\fX}!}j_{\tilde{\fX}}^*$ and a homological shifting; (2) is straightforward. 

It follows from the same proof that this natural isomorphism \eqref{eqn:Psi_hyp} moreover is compatible with the natural transform $r^*_A\to \Psi_{f_A}$ and $r^*_X\to \Psi_{f}$.

Now for the vanishing cycle functor, we apply $r_{A*}$ to the isomorphism of natural transforms 
\[
p_{X*}j_{X}^!\circ (r^*_X\to \Psi_{f})=(r^*_A\to \Psi_{f_A})\circ p_{\fX*}j_{\fX}^!,\]
the left hand side becomes $r_{A*}p_{X*}j_X^!$ composed with $(r^*_A\to \Psi_{f_A})$. Note that the right square of diagram~\eqref{eqn:nearby} is cartesian, hence $r_{A*}p_{X*}j_X^!= p_{\fX*}r_{\calA*}j^!_{X}=p_{\fX*}j_{\fX}^!r_{X*}$. Therefore, we get $p_{\fX*}j_{\fX}^!r_{X*}\circ(r^*_X\to \Psi_{f})=r_{A*}(r^*_A\to \Psi_{f_A})\circ p_{\fX*}j_{\fX}^!$. Taking the mapping cones on both sides, we get \[p_{\fX*}j_{\fX}^!\varphi_f=\varphi_{f_A}p_{\fX*}j_{\fX}^!.\]

According to  Lemma~\ref{lmm:smooth}, $\calA_{\fX}\to \fX$ is a regular embedding, and $\calA_{\fX}\to \fX^A$ is an affine bundle. Hence, $p_*\eta^!\C_{\fM_{\vec{r}}(n)}$ is $\C_{\fM_{\vec{r}}(n)^A}$ up to a homological shifting. Moreover, the codimension of $\calA_{\fX}\to \fX$ is equal to the dimension of the fibers of $\calA_{\fX}\to \fX^A$, hence the shifting is the codimension of $\fX^A$ in $\fX$.
Therefore,   $p_*\eta^!\C_{\fM_{\vec{r}}(n)}[\dim \fX]\cong \C_{\fM_{\vec{r}}(n)^A}[\dim \fX^A]$.
$\blacksquare$

{\it Proof of Proposition~\ref{prop:hyper_def}.}  
Proposition~\ref{prop:hyper_def} follows directly by combining the two parts in this lemma. 
$\blacksquare$

Now we define the map 
\[
h: H^*_{c,T}(\fM_{\vec{r}},\varphi)^\vee\to H^{*+d}_{c,T}(\fM_{\vec{r'}}\times\fM_{\vec{r''}},\varphi)^\vee
\]
as follows. 

 Recall here that $T$ is the maximal torus of $\GL_{\vec{r}}$.
We have the natural adjunction $\eta_!\eta^!\to \id$. Applying to $\varphi_{W}\C_{\fM_{\vec{r}}}$, we get $\eta_!\eta^!\varphi_{W}\C_{\fM_{\vec{r}}}\to \varphi_{W}\C_{\fM_{\vec{r}}}$. Apply to it the functor $\bbD\circ p_{\fM_{\vec{r}}!}$, we obtain 
\begin{align*}
H^*_{c,T}(\fM_{\vec{r}},\varphi)^\vee&=\bbD p_{\fM_{\vec{r}}!} \varphi_{W}\C_{\fM_{\vec{r}}} && \text{(by definition)}\\
&\to\bbD p_{\fM_{\vec{r}}!} \eta_!\eta^!\varphi_{W}\C_{\fM_{\vec{r}}} && \text{(by the adjunction $\eta_!\eta^!\to \id$)}\\
&=\bbD  p_{\calA_{\fM_{\vec{r}}}!}\eta^!\varphi_{W}\C_{\fM_{\vec{r}}} && \text{(by $p_{\fM_{\vec{r}}} \circ \eta= p_{\fM_{\vec{r}}^A}$)}\\
&\cong\bbD  p_{\fM_{\vec{r}}^A!}p_*\eta^!\varphi_{W}\C_{\fM_{\vec{r}}}[-d] 
&& \text{(by $p$ is an affine bundle)}\\
&\cong\bbD  p_{\fM_{\vec{r}}^A!}\varphi_{W}\C_{\fM_{\vec{r}}^A}[d]  &&\text{(by Proposition \ref{prop:hyper_def})}\\
&\cong H^{*+d}_{c,T}(\fM_{\vec{r}}^A,\varphi)^\vee && \text{(by definition)}. 
\end{align*}
Keeping in mind that $\fM_{\vec{r}}^A \cong \fM_{\vec{r'}}\times\fM_{\vec{r''}}$, therefore we have obtained
\[
h: H^*_{c,T}(\fM_{\vec{r}},\varphi)^\vee\to H^{*+d}_{c,T}(\fM_{\vec{r'}}\times\fM_{\vec{r''}},\varphi)^\vee.
\]
The functor $p_0\circ \eta^!$ is the hyperbolic restriction functor \cite{Bra}.
Note that this map $h$ is induced by the natural isomorphism in Proposition \ref{prop:hyper_def},  hence will be referred to as the hyperbolic localization map.

\begin{rmk}\label{rmk:isom up to localization}
The map $h$ is an isomorphism up to localization. 
More explicitly,  both $H^*_{c,T}(\fM_{\vec{r}},\varphi)^\vee$ and $H^{*+d}_{c,T}(\fM_{\vec{r'}}\times\fM_{\vec{r''}},\varphi)^\vee$ 
are modules over 
\[
H^*_{c,T}(\pt)^\vee=\mathbb{C}[\Lie T].\]  
The Lie algebra of $T$ has coordinates $(\mu_1^1,\dots, \mu_1^{r_1},\mu_2^1,\dots, \mu_2^{r_2},\mu_3^1,\dots, \mu_3^{r_3})$.
 Let $\mathbb{C}(\Lie T)$ be the quotient field of $\mathbb{C}[\Lie T]$. 
By the algebraic equivariant localization \cite[Theorem~6.2(2)]{GKM}, the map $h$ induces an isomorphism between 
$H^*_{c,T}(\fM_{\vec{r}},\varphi)^\vee\otimes_{\mathbb{C}[\Lie T]} \mathbb{C}(\Lie T)$ and $H^{*+d}_{c,T}(\fM_{\vec{r'}}\times\fM_{\vec{r''}},\varphi)^\vee\otimes_{\mathbb{C}[\Lie T]} \mathbb{C}(\Lie T)$ over $\mathbb{C}(\Lie T)$.
\end{rmk}

The following is the geometric interpretation of the formulas in Proposition~\ref{prop:coprod_formula}.
In order to keep the flow of the paper, we present the proof in Appendix~\ref{subsec:proof_hyperLoc}.
\begin{prp}\label{prp:coprod_hyperbolic}
The coproduct $\Delta^{\vec{c}}$ from Proposition~\ref{prop:coprod_formula} is compatible with the hyperbolic localization. That is, for $\vec{r}=\vec{r}'+\vec{r}''$, we have
\[
h(\alpha\bullet x)=\Delta^{\vec{c}}(\alpha)\bullet h(x)\] for $\alpha\in  \SH^{\vec{c}}$ and $x\in V_{\vec{r}}$.  
\end{prp}

As a special case, take the subtorus $A\subset \GL_{r_1}\times\GL_{r_2}\times\GL_{r_3}$ to be generic enough, so that we have 
 $ \fM_{\vec{r}}^A= (\fM_{100})^{\times r_1}\times   (\fM_{010})^{\times r_2} \times  (\fM_{001})^{\times r_3}$. 
In this case, the hyperbolic localization gives a map 
\[
h: V_{r_1,r_2,r_3} \to (V_{1, 0, 0})^{\otimes r_1} \otimes (V_{0, 1, 0})^{\otimes r_2}\otimes (V_{0, 0, 1})^{\otimes r_3}, \] 
where the target is the Fock space of $(r_1+r_2+r_3)$-free Bosons. Furthermore, by Remark \ref{rmk:isom up to localization}, 
after tensoring the quotient field $\mathbb{C}(\Lie (T_{r_1}\times T_{r_2}\times T_{r_3}))$, $h$ becomes an isomorphism. 
In particular, we have a free field realization of $\SH^{\vec{c}}$ coming from the coproduct $\Delta^{\vec{c}}$, or equivalently, hyperbolic localization. 
\begin{equation}\label{diag:hyperbolic restriction}
\xymatrix@C=4.5em@R=0.5em{
\SH^{\vec{c}}\ar[r]^{(\Delta^{\vec{c}})^{r_1+r_2+r_3}}&(\SH^{\vec{c}})^{\otimes 
(r_1+r_2+r_3)}\\
&\\
\ar@(ul,ur)[]&\ar@(ul,ur)[]\\
V_{r_1,r_2,r_3}\ar[r]^(0.3){\text{hyperbolic}}_(0.3){\text{localization}}&(V_{1, 0, 0})^{\otimes r_1} \otimes (V_{0, 1, 0})^{\otimes r_2}\otimes (V_{0, 0, 1})^{\otimes r_3}
}
\end{equation}

\section{Vertex operator algebra $\W_{r_1, r_2, r_3}$}\label{sec:VOA_W}
In this section, we review the free field representation and some properties of a class of one-parameter families of vertex operator algebras
$\W_{r_1, r_2, r_3}$ for $r_i\in \Z_{\geq 0}, i=1,2,3$, introduced in \cite{GR}. Generic modules of $\W_{r_1, r_2, r_3}$
can be identified with the modules $V_{r_1,r_2,r_3}$ of the algebra $\SH^{\vec{c}}$ as shown by comparing the free field 
representation with the one coming from the coproduct of $\SH^{\vec{c}}$. 

\subsection{$\W_{r_1, r_2, r_3}$ as extensions of the Virasoro algebra}

The algebra $\W_{r_1, r_2, r_3}$ to be defined in the next section is generated by the modes of
fields $W_1,W_2,\dots, W_n$ of conformal weight $1, 2,\dots, n$ for some integer $n$
\[
W_i(z)=\sum_{l}W_{i,l}z^{-l-i}. 
\]
It is an extension of the Virasoro algebra with the modes of $W_2$ satisfy the Virasoro relations 
\[
\left [W_{2,m},W_{2,n}\right ]=(m-n)W_{2,m+n}+\frac{c}{12}(m^3-m)\delta_{m+n,0}
\]
for $c$ a central element called the central charge and $W_i$ for $i\neq 2$ are primary fields of conformal weight $i$, i.e. their modes have
the following commutation relations with the Virasoro generators
\[
[W_{2,n},W_{i,m}]=((i-1)n-m)W_{i,m+n}.
\]
The generators $W_{i,l}$ are endomorphisms on a vector space spanned by
\[
W_{i_1,-l_1}\dots W_{i_t,-l_1}|0\rangle,\qquad l_i\geq i, \qquad t \geq 0
\]
for $|0\rangle$ the vacuum state, i.e. a state annihilated by all $W_{i_t,-l_1}$ with $l_i< i$. For general values of parameters $r_1,r_2,r_3$, the algebra $\W_{r_1, r_2, r_3}$ is not freely
generated and one needs to work modulo modes of some composite fields to satisfy Jacobi identities.

It is also convenient to introduce the derivative of fields and the operation of normal ordered product. 
Modes of derivatives of a field $\mathcal{O}$ are defined by
\[
(\partial^n \mathcal{O})_m=(-1)^n\prod_{i=0}^{n-1} 
(h_{\mathcal{O}}+m+i)\mathcal{O}_m
\]
where $h_{\mathcal{O}}$ is the conformal dimension ($[W_{2,0},\cdot]$ eigenvalue) of the 
field $\mathcal{O}$. On the other hand, modes of the normal ordered product 
are given by
\[
(\mathcal{O}^1\mathcal{O}^2)_m=\sum_{n\leq 
-h_{\mathcal{O}^1}}\mathcal{O}^1_{n}\mathcal{O}^2_{m-n}+\sum_{n> 
-h_{\mathcal{O}^1}}\mathcal{O}^2_{m-n}\mathcal{O}^1_{n}.
\]

\subsection{Definition of $\W_{r_1, r_2, r_3}$}
Let $\mathcal{H}$ be the Heisenberg algebra generated by modes $b_i$ $i\in \mathbb{Z}$. We define $\W_{r_1,r_2,r_3}$ as a subalgebra of a tensor product of $m=r_1+r_2+r_3$ Heisenberg algebras. The algebra $\W_{r_1,r_2,r_3}$ can be defined according to \cite{BFM,LS,PR2} as an intersection of kernels of vertex operators of the form
\[
V[\alpha](z)=T_{\alpha}\exp \left [\alpha\sum_{n<0}\frac{b_n}{n}z^{-n}\right ]\exp \left  [\alpha\sum_{n>0}\frac{b_n}{n}z^{-n}\right ],
\]
associated to these $m$ Heisenberg algebras, where $T_{\alpha}$ sends the vacuum $|0\rangle$ to a more general highest weight vector $|\alpha\rangle$ satisfying $b_0|\alpha\rangle= \alpha |\alpha\rangle$.
\begin{defn}\label{def:screen}
Let $\mathcal{H}^{\otimes m}$  be a tensor product of $m$ Heisenberg vertex operator 
algebras generated by fields $b^{(3)}_i(z), b^{(2)}_j(z), b^{(1)}_k(z)$ for $i=1,\dots, r_3$, $j=r_3+1,\dots,r_3+ r_2$, $k=r_3+ r_2+1,\dots, r_3+ r_2+r_1$ of conformal weight 1
together with commutation relations
\[
[b^{(a)}_i, b^{(b)}_j]= -\frac{\hbar_a}{\hbar_1\hbar_2\hbar_3}\delta_{a,b}\delta_{i+j,0}.
\label{boson_norm}
\]
Let us consider a set of $m-1$ screening currents given by products of vertex operators of the form 
\begin{eqnarray}\nonumber
S^{33}_i(z)&=&V^{(1)}_i\left [-\hbar_1\right ](z)\otimes V^{(1)}_{i+1}\left [\hbar_1\right ](z),\qquad i\in [1,r_3-1]\\ \nonumber
S^{32}_{r_3}(z)&=&V^{(1)}_{r_3}\left [-\hbar_2\right ](z)\otimes V^{(2)}_{1}\left [\hbar_3\right ](z)\\ \nonumber
S^{22}_i(z)&=&V^{(2)}_i\left [-\hbar_3\right ](z)\otimes V^{(2)}_{i+1}\left [\hbar_3\right ](z),\qquad i\in [r_3+1,r_3+r_2-1]\\ \nonumber
S^{21}_{r_3+r_2}(z)&=&V^{(2)}_{r_2}\left [-\hbar_1\right ](z)\otimes V^{(3)}_{1}\left [\hbar_2\right ](z)\\ \nonumber
S^{11}_i(z)&=&V^{(3)}_i\left [-\hbar_2\right ](z)\otimes V^{(3)}_{i+1}\left [\hbar_2\right ](z),\qquad i\in [r_3+r_2+1,r_3+r_2+r_1-1]
\end{eqnarray}
where $V_{i}^{(\kappa)}[\alpha]$ is a vertex operator for $b^{(\kappa)}_i$. The algebra $\W_{r_1,r_2,r_3}$ is then defined as 
the subalgebra of $\mathcal{H}^{\otimes m}$ given by an intersection 
of kernels of the zero modes of the above screening currents $\oint dz S$.
\end{defn}

Let us now give a few important examples.

\begin{exa}
The algebra $\W_{0, 0, 1}$ can be simply identified with the Heisenberg algebra $\mathcal{H}$.
\end{exa}
\begin{exa}
The algebra $\W_{0, 0, 2}$ can be identified with the Virasoro algebra tensored with the Heisenberg 
algebra $Vir \times \mathcal{H}$.
\end{exa}
\begin{exa}
More generally, the algebra  $\W_{0, 0, r}$ is the well known $\mathcal{W}_r$-algebra \cite{FF2,DF}. This algebra 
is closely related to the specialization $\SH^{(r)}$ introduced by Schiffmann-Vasserot 
in \cite{SV2}. 
\end{exa}
\begin{exa}
$\W_{0,1,1}$ can be identified with a quotient of the well studied algebra $\mathcal{W}_3\times 
\mathcal{H}$ at central charge $c=-2$ from \cite{Wang}. We will return to this example later.
\end{exa}

Apart from the free field realization, there exist two more (conjecturally equivalent) definitions 
of the algebra $\W_{r_1, r_2, r_3}$. 
It was argued in \cite{PR} that $\W_{r_1, r_2, r_3}$ can be viewed as a quotient of a
family of algebras $\mathcal{W}_{1+\infty}[\lambda_1,\lambda_2,\lambda_3]$ 
depending on three complex parameters $\lambda_i\in \C$ subject to the 
constraint\footnote{The parameters $\lambda_i$ can be expressed in terms of the parameters above as 
$\lambda_i=(r_1\hbar_1+r_2\hbar_2+r_3\hbar_3)/\hbar_i$. Note also that the commutation relations of the VOA depend only
on the scaling invariant combination $\Psi=-\hbar_1/\hbar_2$.}
\[
\frac{1}{\lambda_1}+\frac{1}{\lambda_2}+\frac{1}{\lambda_3}=0
\]
 in the same way as $\SH^{(r_1, r_2, r_3)}$ are specializations of $\SH^{\vec{c}}$.  The algebra 
$\mathcal{W}_{1+\infty}[\lambda_1,\lambda_2,\lambda_3]$ was recently constructed 
in \cite{L} based on the previous work of \cite{H,GG,P1}. It was 
argued\footnote{See also \cite{L} and references therein for discussion of 
special classes of such truncations.} in \cite{P2,PR} that the algebra 
$\mathcal{W}_{1+\infty}[\lambda_1,\lambda_2,\lambda_3]$ contains an ideal
whenever the parameters are specialized to
\[
\frac{r_1}{\lambda_1}+\frac{r_2}{\lambda_2}+\frac{r_3}{\lambda_3}=1.
\label{trunc}
\]
Conjecturally, the algebra $\W_{r_1,r_2,r_3}$ can be defined as a quotient of  
$\mathcal{W}_{1+\infty}[\lambda_1,\lambda_2,\lambda_3]$ by such an ideal and depends 
on the continuous parameter $\Psi=-\frac{\lambda_1}{\lambda_2}=-\frac{\hbar_2}{\hbar_1}$.

The original definition of $\W_{r_1,r_2,r_3}[\Psi]$ from \cite{GR} was given in terms of 
a combination of the quantum  Drinfeld-Sokolov reduction \cite{FF1} and the BRST coset \cite{HR}
of super Kac-Moody algebras. We refer reader to \cite{GR,PR} for detailed discussion.

The definition of the algebra is invariant under the mutual permutation of the parameters $r_i$
and $\hbar_i$. In \cite{GR}, this triality symmetry was a non-trivial consequence of the S-duality
of boundary conditions. From the spiked instanton configuration point of view, this triality is 
a trivial consequence of relabeling of the coordinates $z_i$ for $i=1,2,3$ associated to
the three $\C\subset \C^3$.

Note also that one might consider different ordering of the free bosons in the free field 
realizations above. Each ordering gives a different realization of the same algebra. For the purpose
of this work, we use the ordering corresponding to the screening charges above.

\subsection{$\W_{0,0,2}$ example in detail}

The simplest non-trivial example $\W_{0,0,2}$ is freely generated by the Heisenberg field $W_1$ and 
the Virasoro field $W_2$. For completeness, let us write down explicitly their free field realization in terms
of a pair of free bosons $b_1^{(3)}$ and $b_2^{(3)}$:
\begin{eqnarray} \nonumber
W_1&=&b_1^{(3)}+b_2^{(3)},\\ \nonumber
W_2&=&-\frac{\hbar_1\hbar_2}{4}\left (\tilde{b},\tilde{b}\right )-\frac{\hbar_3}{2}\partial \tilde{b}
\end{eqnarray}
where we have introduced\footnote{Note that modes of $W_1$ commute with the modes of $W_2$. One can recover the standard commutation relations by adding a multiple of normal ordered product $(W_1W_1)$ to $W_2$.}
\[
 \tilde{b}=b_1^{(3)}-b_2^{(3)}.
\]

\subsection{$\W_{0,1,1}$ example in detail}
\label{subsec:W3-2}
Let us further discuss $\W_{0,1,1}$ that plays an important role in the 
comparison of the geometric action on the equivariant cohomology on the moduli 
space of spiked instantons with generic modules of $\W_{r_1,r_2,r_3}$. The 
algebra $\W_{0,1,1}$ is generated by modes of the Heisenberg algebra field $W_1$, 
the stress-energy tensor $W_2$ and a primary field of spin three $W_3$ with the 
following commutation relation
\begin{eqnarray}
\left [W_{3,m},W_{3,n}\right 
]&=&\frac{4}{3}(m-n)\sum_{k=-\infty}
^{\infty}W_{2,m+n-k}W_{2,k}-\frac{6}{5}m(m^2-1)(m^2-4)\delta_{m,-n} \notag\\
&&-\sqrt{6}(m-n)(m+2)(2m+3)W_{2,m+n}. \notag
\end{eqnarray}
The algebra $\W_{0,1,1}$ is not generated freely and there exists a null field
\[
X=(W_3W_3)-128(W_2(W_2W_2))-76(\partial W_2 \partial 
W_2)-112(\partial^2W_2 W_2)+\frac{32}{3}\partial^4 W_2
\]
whose modes need to be factored out. This combination will be automatically zero in 
the free field realization of the VOA as one can easily check.

The full vertex operator algebra decomposes into the subalgebra of positive, negative and zero modes. 
Let us now describe how to induce the modules of interest from the one 
dimensional modules for the algebra of positive and zero modes. Let us consider a 
vector $|u,h,w\rangle$ annihilated by all the positive modes and being an 
eigenstate of the zero modes
\begin{eqnarray}\nonumber
W_{i,n}|w_1,w_2,w_3\rangle=0 \ \mbox{for}\ n>0,i=1,2,3,\\ \nonumber
W_{i,0}|w_1,w_2,w_3\rangle=w_i|w_1,w_2,w_3\rangle .
\end{eqnarray}
Note also that the algebra of zero modes is commutative.
The existence of the null state $X$ puts a constraint on the weights $u,h,w$ 
since $X_0|w_1,w_2,w_3\rangle=0$. Acting by $X_0$ on the state $|w_1,w_2,w_3\rangle$ leads to 
the constraint \cite{Wang,PR2}
\[
w_3^2=16w_2^2(8w_2+1).
\]
Having a triple of numbers $w_1,w_2,w_3$ satisfying the above condition, one 
can consider a module of the full algebra generated by an action of negative 
modes\footnote{Highest weight representations of this form (and similarly for a general $\mathcal{W}_{r_1,r_2,r_3}$) 
are paremetrized by the spectrum of the Zhu algebra \cite{BN,Z,L} that turn out to be commutative in the case of 
$\mathcal{W}_{r_1,r_2,r_3}$. We expect the representation theory associated to more general toric Calabi-Yau
$3$-folds to be more complicated. \cite{PR} conjectured appearance of modules induced from generic Gelfand-Tsetlin modules of \cite{FGR} and various irregular modules  of \cite{GP,GT}.} 
on the vector $|w_1,w_2,w_3\rangle$ and quotient by the negative modes of $X$.

Using the free field realization in terms of a pair of free bosons, one can 
induce the above modules from the vertex operators (Fock modules) of the free fields. Let us 
consider two Heisenberg algebras generated by fields $b^{(2)}_1$ and $b^{(3)}_2$
normalized as above.
The fields $W_1,W_2,W_3$ can then be then realized as\footnote{For simplicity, we again decouple 
$W_1$ as in the Virasoro case above.}
\begin{eqnarray}\nonumber
W_1&=&b^{(2)}_1+b^{(3)}_2\\ \label{W011}
W_2&=&\frac{1}{2}\left [ (\tilde{b}\tilde{b})+\partial \tilde{b} \right ]\\ \nonumber
W_3&=& 
4(\tilde{b}(\tilde{b}\tilde{b}))+6(\tilde{b}\partial \tilde{b})+\partial^2\tilde{b}
\end{eqnarray}
where now we have introduced
\[
\tilde{b}=\hbar_3b^{(2)}_1-\hbar_2b^{(3)}_2
\]
having trivial commutation relations with $W_1$.

Generic modules can be induced using the free field realization by an action of the 
positive modes of $W_1,W_2,W_3$ on the highest weight state
\[
|q_1,q_2\rangle =V_1^{(3)}[q_1](0)|0\rangle\otimes V_1^{(2)}(0)[q_2](0)|0\rangle.
\]
One can in particular compute the action of the zero modes to establish the 
conection with the parameters $w_1,w_2,w_3$ above. One finds
\begin{eqnarray}\nonumber
w_1&=&-\frac{q_1}{\hbar_1 \hbar_3}-\frac{q_2}{\hbar_1\hbar_2}\\ \nonumber
w_2&=&\frac{q (q+1)}{2}\\ \nonumber
w_3&=&-2q(q+1)(2q+1)
\end{eqnarray}
where we have introduced $q=(q_1-q_2)/\hbar_1$.

\section{The action of $\W_{r_1, r_2, r_3}$ on framed quiver varieties}
\label{sec:Y110}

In this section, we show the following.
\begin{thm}\label{thm in sec9}
The action of  $\SH^{\vec{c}}$ on  $V_{r_1, r_2, r_3}$, constructed in \S~\ref{subsec:central},  factors through the vertex operator algebra $\calW_{r_1,r_2,r_3}$.
\end{thm}
In \S~\ref{sec10_1}, we define a map $\SH^{\vec{c}}\to \calH^{\otimes |\vec{r}|}$. Using the calculations of \S~\ref{sec10_2} and \ref{sec10_3}, in \S~\ref{sec10_4} we prove that the image of the map $\SH^{\vec{c}}\to \calH^{\otimes |\vec{r}|}$ lands in $\calW_{r_1,r_2,r_3}\subseteq \calH^{\otimes |\vec{r}|}$.
In particular, this finishes Step 4 in proof of Theorem \ref{thm:main}, hence complete the proof of the main theorem.

\subsection{Action on $V_{0,0,1}$}\label{sec10_1}
The Heisenberg elements $B_i$ of the algebra $\SH^{\vec{c}}$ acting on $V_{0,0,1}$ can be identified with the Heisenberg VOA generators $b_i^{(3)}$. Moreover,
the action of $f_{0,1} \in \SH^{\vec{c}}$ 
on such module can be identified according to \cite[(8.46)]{SV2} as
\begin{eqnarray}\nonumber
f_{0,1}&\mapsto &\frac{\hbar_1^2\hbar_2^2}{2}\sum_{k,l\geq 1}\left ( b^{(3)}_{-l-k}b^{(3)}_{l}b^{(3)}_{k}+b^{(3)}_{-l}b^{(3)}_{-k}b^{(3)}_{l+k} \right ) \label{(56)}\\
& +&\frac{\hbar_1 \hbar_2 \hbar_3}{2}\sum_{l\geq 1}(l-1)b^{(3)}_{-l}b^{(3)}_{l}+\mu_1^{(3)}\hbar_1\hbar_2\sum_{l\geq 1}b^{(3)}_{-l}b^{(3)}_{l}, 
\label{f_one} \label{(57)}
\end{eqnarray}
where $\mu_1^{(3)}$ is the Chern root of the line bundle $\calE_{0, 0, 1}$ on $\calM_{(0, 0, 1)}$ (see \S\ref{sub:taut bundle}). 

Similarly for the free bosons in other directions, one get analogous expression with the parameters $\hbar_i$ permuted. 

The algebra $\SH^{\vec{c}}$ is generated by the Heisenberg elements $\{B_{i}\}$ and $f_{0, 1}$ over $F^{\vec{c}}$.
Indeed, when $\bold{c}^{(1)}=\bold{c}^{(2)}=0$, it is shown in \cite[pg. 217, Proposition. (c)]{SV2} that the resulting algebra $\SH^{\bold{c}^{(3)}}$ is generated by such elements. The claim now follows from the embedding $\SH^{\vec{c}}\subset \bigoplus_{k=1, 2, 3}\SH^{\bold{c}^{(k)}}$ (see \S~\ref{subsec:coprod_form_conclud}). Note that in Lemma \ref{Heisen}, the element $G_0=\frac{\hbar_1 \bold{c}_0^{(1)}+\hbar_2 \bold{c}_0^{(2)}+\hbar_3 \bold{c}_0^{(3)}}{\hbar_1\hbar_2\hbar_3}$ lies in $F^{\vec{c}}$.

By the formula \eqref{(56)}\eqref{(57)}, the image of $f_{0, 1}$ in $\End(V_{0, 0, 1})$ under the action $\SH^{\vec{c}} \to \End(V_{0,0,1})$ 
lies in the image of $\calH$ in $\End(V_{0, 0, 1})$. Therefore, the action of  $\SH^{\vec{c}}$ on $V_{0,0,1}$ factors through the Heisenberg algebra. Similarly for the actions on $V_{0, 1, 0}, V_{1, 0, 0}$. 
Let $m=r_1+r_2+r_3 \in \Z_{\geq 0}$. One can compose the coproduct 
$(\Delta^{\vec{c}})^m:\SH^{\vec{c}} \rightarrow (\SH^{\vec{c}})^{\otimes m}$ with the maps $\SH^{\vec{c}}$ to the three Heisenberg algebras 
with $r_i$ $\SH^{\vec{c}}$ factors mapped in the $i$'th way. Thus, the action of $\SH^{\vec{c}}$ on 
$(V_{1, 0, 0})^{\otimes r_1}\otimes (V_{0,1, 0})^{\otimes r_2} \otimes (V_{0, 0, 1})^{\otimes r_3}$ factors through the tensor $\calH^{\otimes m}$. 

\subsection{Comparison for $V_{0,0,2}$}\label{sec10_2}
We have the hyperbolic restriction map $h: V_{0, 0, 2}\to V_{0,0,1}\otimes V_{0, 0, 1}$. Keeping in mind also Proposition~\ref{prop:shuffle specialization}, 
by the computation in \cite[\S 8.8]{SV2}, we have the following. 
\begin{prp}
The action of $\SH^{(0,0,2)}$ on $V_{0,0,1}\otimes V_{0, 0, 1}$ factors through 
$\W_{0,0,2}$.
\end{prp}
The algebra $\SH^{\vec{c}}$ is generated by $\{B_{i}\}$ and $f_{0, 1}$. 
It suffices to show the actions of $\{B_{i}\}$ and $f_{0, 1}$ on $V_{0,0,1}\otimes V_{0, 0, 1}$ factor through $\W_{0,0,2}$. 
This is done in \cite[Lemma 8.20]{SV2}.

\subsection{Comparison for $V_{0,1,1}$}\label{sec10_3}
We have the hyperbolic restriction map $h: V_{0, 1, 1}\to V_{0,1,0}\otimes V_{0, 0, 1}$. 
\begin{prp}
The action of $\SH^{(0,1,1)}$ on $V_{0,1,0}\otimes V_{0, 0, 1}$ factors through 
$\W_{0,1,1}$.
\end{prp}
{\it Proof.}
It suffices to consider the generators $\{B_{i}\}$ and $f_{0, 1}$ of $\SH^{\vec{c}}$. 
The action of $B_i$ is obvious since we can identify
\[
W_1=b^{(2)}_1+b^{(3)}_2.
\]
The action of $f_{0,1}$ can be checked by expressing the free field realization of $f_{0,1}$ in terms of modes of the freely realized $W_i$.
This can be seen along the lines of \cite[\S 8]{SV2}. From the coproduct and the expression (\ref{f_one}), the free field realization of $f_{0,1}$ is
\begin{eqnarray}\nonumber
f_{0,1}&=&\frac{\hbar_1 \hbar_2 \hbar_3}{2}\sum_{l\geq 1}(l-1)b^{(2)}_{-l}b^{(2)}_{l}+\frac{\hbar_1^2\hbar_3^2}{2}\sum_{k,l\geq 1}\left ( b^{(2)}_{-l-k}b^{(2)}_{l}b^{(2)}_{k}+b^{(2)}_{-l}b^{(2)}_{-k}b^{(2)}_{l+k} \right )\\ \nonumber
&&+\frac{\hbar_1 \hbar_2 \hbar_3}{2}\sum_{l\geq 1}(l-1)b^{(3)}_{-l}b^{(3)}_{l}+\frac{\hbar_1^2\hbar_2^2}{2}\sum_{k,l\geq 1}\left ( b^{(3)}_{-l-k}b^{(3)}_{l}b^{(3)}_{k}+b^{(3)}_{-l}b^{(3)}_{-k}b^{(3)}_{l+k} \right )\\ \nonumber
&&+\mu_1^{(2)}\hbar_1\hbar_3\sum_{l\geq 1}b^{(2)}_{-l}b^{(2)}_{l}+\mu_1^{(3)}\hbar_1\hbar_2\sum_{l\geq 1}b^{(3)}_{-l}b^{(3)}_{l}+\hbar_1\hbar_2\hbar_3\sum_{l\geq 1}l b^{(2)}_{l}b^{(3)}_{-l}. 
\end{eqnarray}
Using $\mu_1^{(2)}=\hbar_1 \hbar_3 b^{(2)}_0$ and $\mu_1^{(3)}=\hbar_1 \hbar_2 b^{(3)}_0+\hbar_3$, we can rewrite the expression as
\begin{align*}\nonumber 
f_{0,1}&=\frac{\hbar_1 \hbar_2 \hbar_3}{4}\sum_{l=-\infty}^{\infty}|l|:\left ( b^{(2)}_{-l}b^{(2)}_{l}+b^{(3)}_{-l}b^{(3)}_{l}\right ):+\frac{\hbar_1 \hbar_2 \hbar_3}{4}\sum_{l=-\infty}^{\infty}:\left ( b^{(3)}_{-l}b^{(3)}_{l}-b^{(2)}_{-l}b^{(2)}_{l}\right ):\\ 
&+\sum_{k,l=-\infty}^{\infty}:\left (\frac{\hbar_1^2\hbar_3^2}{6} b^{(2)}_{-l-k}b^{(2)}_{l}b^{(2)}_{k}+\frac{\hbar_1^2\hbar_2^2}{6}b^{(3)}_{-l-k}b^{(3)}_{l}b^{(3)}_{k} \right ):+\hbar_1\hbar_2\hbar_3\sum_{l\geq 1}l b^{(2)}_{l}b^{(3)}_{-l}+c
\end{align*}

Let us show that this expression can be written in terms of modes of $W_1,W_2,W_3$ from (\ref{W011}). Similarly as in the Virasoro case, one can easily get rid of the term containing  the absolute value by subtracting a $\hbar_1 \hbar_2 \hbar_3/4$ multiple of
\[
\sum_{l=-\infty}^{\infty}|l|:W_{1,-l}W_{1,l}:=2\sum_{l=-\infty}^{\infty}|l|:b_{-l}^{(2)}b_{l}^{(3)}:+\sum_{l=-\infty}^{\infty}|l|:b^{(2)}_{-l}b^{(2)}_l+b^{(3)}_{-l}b^{(3)}_l:
\]
The term containing absolute value cancels and the tail of the twisted coproduct combine with the first term on the right hand side into an expression of the form of a zero mode of some vertex operator. In total, one gets
\begin{eqnarray}\nonumber 
&&\frac{\hbar_1 \hbar_2 \hbar_3}{4}\sum_{l\in \Z}:\left ( b^{(2)}_{-l}b^{(2)}_{l}-b^{(3)}_{-l}b^{(3)}_{l}\right ):-\frac{\hbar_1\hbar_2\hbar_3}{2}\sum_{l=-\infty}^{\infty}l b^{(2)}_{-l}b^{(3)}_l\\ \nonumber 
&&+\sum_{k,l\in \Z}:\left (\frac{\hbar_1^2\hbar_3^2}{6} b^{(2)}_{-l-k}b^{(2)}_{l}b^{(2)}_{k}+\frac{\hbar_1^2\hbar_2^2}{6}b^{(3)}_{-l-k}b^{(3)}_{l}b^{(3)}_{k} \right ):+c.
\end{eqnarray}
Note now that subtracting the zero mode of the combination
\begin{eqnarray}\nonumber 
\frac{1}{24}(\hbar_3-\hbar_2)W_3+\hbar_2\hbar_3 (W_1,W_2)+\frac{\hbar_2^2\hbar_3^2}{6}(W_1,(W_1,W_1))\\ \nonumber 
-\frac{\hbar_2^2\hbar_3^2}{2(\hbar_2-\hbar_3)} (W_1,\partial W_1)+\frac{\hbar_1^2+2\hbar_1 \hbar_2+2\hbar_2^2}{4(\hbar_2-\hbar_3)}\partial W_2
\end{eqnarray}
cancels the cubic term together with the sum containing $l$. One is left with
\begin{eqnarray}\nonumber 
&&\frac{\hbar_1 \hbar_2 \hbar_3}{4}\sum_{l\in \Z}:\left ( b^{(2)}_{-l}b^{(2)}_{l}-b^{(3)}_{-l}b^{(3)}_{l}\right ):+\frac{\hbar_1 \hbar_2 \hbar_3}{2}\frac{\hbar_2+\hbar_3}{\hbar_2-\hbar_3}\sum_{l=\infty}^{\infty}b^{(3)}_{-l}b^{(2)}_{l}+c
\end{eqnarray}
Finally, subtracting
\begin{eqnarray}\nonumber 
\frac{1}{\hbar_1(\hbar_3-\hbar_2)}W_{2,0}-\frac{\hbar_1^2+2\hbar_1\hbar_2+2\hbar_2^2}{4\hbar_1(\hbar_2-\hbar_3)}(W_1W_1),
\end{eqnarray}
one can get rid of all the terms containing $b^{(j)}$. 
This completes the proof. 
$\blacksquare$

\subsection{Comparison for $V_{r_1,r_2,r_3}$}
\label{sec10_4}
To summarize, we have proven that the action of $\SH^{(0,0,2)}$ on $V_{0,0,1}\otimes V_{0,0,1}$ factors through $\W_{0,0,2}$ and the action of $\SH^{(0,1,1)}$ on $V_{0,1,0}\otimes V_{0,0,1}$ factors through $\W_{0,1,1}$. Analogous proposition holds also for the other combinations $\W_{0,2,0}$, $\W_{2,0,0}$, $\W_{1,1,0}$ simply by permuting the parameters $\hbar_i$.
Keeping in mind the definition of $\W_{r_1,r_2,r_3}$ in terms of an intersection of screening-charges kernels, the consequence of the above is the action of $\SH^{(r_1,r_2,r_3)}$ on $(V_{1, 0, 0})^{\otimes r_1}\otimes (V_{ 0,1, 0})^{\otimes r_2}\otimes (V_{0, 0, 1})^{\otimes r_3}
$ factors through $\W_{r_1,r_2,r_3}$. 

Recall by \eqref{diag:hyperbolic restriction}, we have the hyperbolic restriction map 
\[
h: V_{r_1, r_2, r_3}\to (V_{1, 0, 0})^{\otimes r_1}\otimes (V_{0,1, 0})^{\otimes r_2} \otimes (V_{0, 0, 1})^{\otimes r_3},\] which is compatible with the coproduct $\Delta^{\vec{c}}$ by 
Proposition \ref{prp:coprod_hyperbolic}. 
Furthermore, by Remark \ref{rmk:isom up to localization}, 
$h$ becomes an isomorphism after tensoring the quotient field $\mathbb{C}(\Lie (T_{r_1}\times T_{r_2}\times T_{r_3}))$. 
This implies the action of $\calW_{r_1,r_2,r_3}$ on $V_{r_1, r_2, r_3}$. 
This finishes the proof of Theorem \ref{thm in sec9}. 

\appendix

\section{The proof of Theorem~\ref{thm:positiveYangian}}
\label{subsec:app_Yangian}
In this section, we prove Theorem~\ref{thm:positiveYangian} (1). 
Let $Y^+$ be the positive part of the affine Yangian $\affY$, and $\bSh$ be the shuffle algebra associated to the 2d COHA $\H^{(Q_3,W_3),{\bf T}_2}_{B_3=0}$ in \S\ref{subsec:2dCOHA}. By Theorem \ref{thm:CritPrep} and Theorem \ref{thm:shuffle}, it suffices to show there is an algebra homomorphism from $Y^+$ to the shuffle algebra $\bSh$, the multiplication of which we denote by $\star$. We now check the assignment 
\[
\Psi: Y^+\to \bSh,  \,\  \text{by  $e_{r} \mapsto \lambda^{r} \in 
\bSh(1)\cong \C[\hbar_1,\hbar_2,\lambda]$}
\]
preserves the relations \eqref{Y1} and 
\eqref{Y6}. 

Let $\lambda_{12}=\lambda_1-\lambda_2$, and let 
\begin{align*}
&\fac(\lambda_{12}):= 
\frac{(\lambda_{12}-\hbar_1)(\lambda_{12}-\hbar_2)(\lambda_{12}-\hbar_3)}{
\lambda_{12}}.\\&
\fac'(\lambda_{12}):=\fac(\lambda_{21})= 
\frac{(\lambda_{12}+\hbar_1)(\lambda_{12}+\hbar_2)(\lambda_{12}+\hbar_3)}{
\lambda_{12}}.
\end{align*}
Let 
$\sigma_{2}:=\hbar_1\hbar_2+\hbar_2\hbar_3+\hbar_1\hbar_3$. 
Under the condition $\hbar_1+\hbar_2+\hbar_3=0$, we have the following identities. 
\begin{align}
&(X + \hbar_1) (X + \hbar_2) (X + \hbar_3) - (X - \hbar_1) (X - \hbar_2) (X - 
\hbar_3)=
2 \hbar_1 \hbar_2 \hbar_3. \label{eq:commutator}\\
&(X + \hbar_1) (X + \hbar_2) (X + \hbar_3) + (X - \hbar_1) (X - \hbar_2) (X - 
\hbar_3)=
2 X^3+2\sigma_2X. \label{eq:symmetric}
\end{align}
Therefore, 
\begin{align*}
&\fac(\lambda_{12})-\fac'(\lambda_{12})=\fac(\lambda_{12})-\fac(\lambda_{21}
)=-2\frac{\hbar_1\hbar_2\hbar_3}{\lambda_{12}}, \\
&\fac(\lambda_{12})+\fac'(\lambda_{12})=\fac(\lambda_{12})+\fac'(\lambda_{21}
)=2\frac{ \lambda_{12}^3+\sigma_2\lambda_{12}}{\lambda_{12}}. 
\end{align*}

Let $R=\C[\hbar_1, \hbar_2]$. By the shuffle formula \eqref{shuffle formula}, the multiplication of $\bSh$ is given by
\begin{align*}
&R[\lambda_1]\otimes R[\lambda_2] \to R[\lambda_1, \lambda_2]^{\mathfrak{S}_2}, 
\\
&(f(\lambda_1), g(\lambda_2))\mapsto f(\lambda_1)g(\lambda_2)
\fac(\lambda_{12})
+f(\lambda_2)g(\lambda_1)
\fac'(\lambda_{12})
\end{align*}
Therefore, for any $a, b\in \bbN$, 
$
\lambda^a\star \lambda^b=\lambda_1^a\lambda_2^b \fac(\lambda_{12})
+\lambda_1^b\lambda_2^a \fac'(\lambda_{12})
$. 
This gives that
\begin{align}
&\lambda^a\star \lambda^b-\lambda^b\star 
\lambda^a=-2\frac{\hbar_1\hbar_2\hbar_3}{\lambda_{12}}\Big(\lambda^a_1 
\lambda^b_2-\lambda^b_1 \lambda^a_2\Big) \notag\\
&\lambda^a\star \lambda^b+\lambda^b\star \lambda^a=\Big(\lambda_1^a\lambda_2^b 
+\lambda^b_1\lambda^a_2  \Big)
2 \frac{\lambda_{12}^3+\sigma_2\lambda_{12}}{\lambda_{12}} \label{eq:lambda 
symm}
\end{align}
Using \eqref{eq:commutator}, we now compute
\begin{align*}
&\Psi([e_{i+3}, e_j]-3[e_{i+2}, e_{j+1}]+3[e_{i+1}, e_{j+2}]-[e_i, 
e_{j+3}]+\sigma_2([e_{i+1}, e_j]-[e_i, e_{j+1}]))\\
=&-2\frac{\hbar_1\hbar_2\hbar_3}{\lambda_{12}}
\Big(
\lambda^{i+3}_1 \lambda^{j}_2-\lambda^{j}_1 \lambda^{i+3}_2
-3(\lambda^{i+2}_1 \lambda^{j+1}_2-\lambda^{j+1}_1 \lambda^{i+2}_2)
+3(\lambda^{i+1}_1 \lambda^{j+2}_2-\lambda^{j+2}_1 \lambda^{i+1}_2)\\&
-(\lambda^{i}_1 \lambda^{j+3}_2-\lambda^{j+3}_1 \lambda^{i}_2)
+\sigma_2(
(\lambda^{i+1}_1 \lambda^{j}_2-\lambda^{j}_1 \lambda^{i+1}_2)
-(\lambda^{i}_1 \lambda^{j+1}_2-\lambda^{j+1}_1 \lambda^{i}_2)
)
\Big)\\
%----------------------
=&-2\frac{\hbar_1\hbar_2\hbar_3}{\lambda_{12}}
\Big( \lambda_1^{i}\lambda_2^j 
(\lambda_{12}^{3}+\sigma_2\lambda_{12})-\lambda_1^{j}\lambda^i_2
(\lambda_{21}^{3}+\sigma_2\lambda_{21})
\Big)\\
=&-2\frac{\hbar_1\hbar_2\hbar_3}{\lambda_{12}}
( \lambda_1^{i}\lambda_2^j+\lambda_1^{j}\lambda^i_2) 
(\lambda_{12}^{3}+\sigma_2\lambda_{12})
\end{align*}
By \eqref{eq:lambda symm}, the above is the same as
\[
-\hbar_1\hbar_2\hbar_3 (\lambda^i\star \lambda^j+\lambda^j\star 
\lambda^i)=\Psi(-\sigma_3 \{e_i, e_j\}). 
\]
Therefore, the assignment $\Psi$ preserves the relation \eqref{Y1}. 

By the shuffle formula \eqref{shuffle formula}, the multiplication of $\bSh$ is given by
\begin{align*}
R[\lambda_1]\otimes R[\lambda_2, \lambda_3]^{\mathfrak{S}_2}& \to R[\lambda_1, 
\lambda_2, \lambda_3]^{\mathfrak{S}_3}, \\
(f(\lambda_1), g(\lambda_2, \lambda_3))
\mapsto 
%---------------
&f(\lambda_1)g(\lambda_2, \lambda_3)
\fac(\lambda_{12})
\fac(\lambda_{13})\\&
%---------------
+f(\lambda_2)g(\lambda_1, \lambda_3)
\fac(\lambda_{21})
\fac(\lambda_{23})
%---------------
+f(\lambda_3)g(\lambda_2, \lambda_1)
\fac(\lambda_{32})
\fac(\lambda_{31}),
\end{align*}
and
\begin{align*}
R[\lambda_1, \lambda_2]^{\mathfrak{S}_2}\otimes R[\lambda_3]& \to R[\lambda_1, 
\lambda_2, \lambda_3]^{\mathfrak{S}_3}, \\
(f(\lambda_1, \lambda_2), g(\lambda_3))
\mapsto 
%---------------
&f(\lambda_1, \lambda_2)g(\lambda_3))
\fac(\lambda_{13})
\fac(\lambda_{23})\\
%---------------
+&f(\lambda_3, \lambda_2)g(\lambda_1))
\fac(\lambda_{31})
\fac(\lambda_{21})
%---------------
+f(\lambda_1, \lambda_3)g(\lambda_2))
\fac(\lambda_{12})
\fac(\lambda_{32}). 
\end{align*}

Therefore, we have
\begin{align*}
&e_{c}\star [e_{a}, e_{b}] =
-2\sigma_3 \lambda_1^c\star \frac{\Big
(\lambda^a_2 \lambda^b_3-\lambda^b_2 \lambda^a_3\Big)}{\lambda_{23}}\\
=&-2\sigma_3 
\frac{(\lambda^a_2 \lambda^b_3-\lambda^b_2 \lambda^a_3)\lambda_1^c 
}{\lambda_{23}}
\fac(\lambda_{12})
\fac(\lambda_{13})
%--------------
-2\sigma_3 \frac{(\lambda^a_1 \lambda^b_3-\lambda^b_1 \lambda^a_3)\lambda_2^c 
}{\lambda_{13}}
\fac(\lambda_{21})
\fac(\lambda_{23})\\&
%-----------
-2\sigma_3 \frac{(\lambda^a_2 \lambda^b_1-\lambda^b_2 \lambda^a_1)\lambda_3^c 
}{\lambda_{21}}
\fac(\lambda_{32})
\fac(\lambda_{31}), 
\end{align*}
and
\begin{align*}
&[e_{a}, e_{b}] \star e_{c}
=-2\sigma_3\frac{\Big(\lambda^a_1 \lambda^b_2-\lambda^b_1 
\lambda^a_2\Big)}{\lambda_{12}}\star \lambda_3^c\\
=
&
-2\sigma_3\frac{(\lambda^a_1 \lambda^b_2-\lambda^b_1 \lambda^a_2)\lambda_3^c 
}{\lambda_{12}}
\fac(\lambda_{13})
\fac(\lambda_{23})
%---------------
-2\sigma_3\frac{(\lambda^a_3 \lambda^b_2-\lambda^b_3 \lambda^a_2)\lambda_1^c 
}{\lambda_{32}}
\fac(\lambda_{31})
\fac(\lambda_{21})\\&
%---------------
-2\sigma_3\frac{(\lambda^a_1 \lambda^b_3-\lambda^b_1 \lambda^a_3)\lambda_2^c 
}{\lambda_{13}}
\fac(\lambda_{12})
\fac(\lambda_{32}). 
\end{align*}

We compute
\begin{align*}
&[e_c, [e_a, e_b]]
=e_c\star [e_a, e_b]-[e_a, e_b]\star e_c\\
=&-2\sigma_3 
\frac{(\lambda^a_2 \lambda^b_3-\lambda^b_2 \lambda^a_3)\lambda_1^c 
}{\lambda_{23}}
\fac(\lambda_{12})
\fac(\lambda_{13})
%--------------
-2\sigma_3 \frac{(\lambda^a_1 \lambda^b_3-\lambda^b_1 \lambda^a_3)\lambda_2^c 
}{\lambda_{13}}
\fac(\lambda_{21})
\fac(\lambda_{23})\\&
%-----------
-2\sigma_3 \frac{(\lambda^a_2 \lambda^b_1-\lambda^b_2 \lambda^a_1)\lambda_3^c 
}{\lambda_{21}}
\fac(\lambda_{32})
\fac(\lambda_{31})\\&
%-----
+2\sigma_3\frac{(\lambda^a_1 \lambda^b_2-\lambda^b_1 \lambda^a_2)\lambda_3^c 
}{\lambda_{12}}
\fac(\lambda_{13})
\fac(\lambda_{23})
%---------------
+2\sigma_3\frac{(\lambda^a_3 \lambda^b_2-\lambda^b_3 \lambda^a_2)\lambda_1^c 
}{\lambda_{32}}
\fac(\lambda_{31})
\fac(\lambda_{21})\\&
%---------------
+2\sigma_3\frac{(\lambda^a_1 \lambda^b_3-\lambda^b_1 \lambda^a_3)\lambda_2^c 
}{\lambda_{13}}
\fac(\lambda_{12})
\fac(\lambda_{32})\\
%-------------------------------------------------
%-------------------------------------------------
=&-2\sigma_3 
\lambda_1^c(\lambda^a_2 \lambda^b_3-\lambda^b_2 \lambda^a_3) 
\Big(\frac{
\fac(\lambda_{12})
\fac(\lambda_{13})}{\lambda_{23}}
+\frac{
\fac(\lambda_{21})
\fac(\lambda_{31})}{\lambda_{32}}
\Big)\\
%------
&-2 \sigma_3 \lambda_2^c (\lambda^a_1 \lambda^b_3-\lambda^b_1 \lambda^a_3)
\Big(
\frac{\fac(\lambda_{21})
\fac(\lambda_{23})}{\lambda_{13} }
+\frac{\fac(\lambda_{12})
\fac(\lambda_{32})}{\lambda_{31} }
\Big)
\\&
%-----------
-2\sigma_3 \lambda_3^c (\lambda^a_2 \lambda^b_1-\lambda^b_2 \lambda^a_1)
\Big(
\frac{\fac(\lambda_{32})
\fac(\lambda_{31})}{\lambda_{21}}
+
\frac{\fac(\lambda_{23})
\fac(\lambda_{13})}{\lambda_{12}}
\Big). 
\end{align*}
Plug the above formula into the following
\begin{align}
 \Sym_{\mathfrak{S}_3}[e_{i_1}, [e_{i_2}, e_{i_3+1}]]
 =&[e_{i_1}, [e_{i_2}, e_{i_3+1}]]
 +[e_{i_2}, [e_{i_1}, e_{i_3+1}]]
+ [e_{i_3}, [e_{i_2}, e_{i_1+1}]] \notag\\&
+ [e_{i_1}, [e_{i_3}, e_{i_2+1}]]
 +[e_{i_2}, [e_{i_3}, e_{i_1+1}]]
+ [e_{i_3}, [e_{i_1}, e_{i_2+1}]] \label{eq:serre check}
\end{align}
The term in \eqref{eq:serre check} involving 
$\lambda_1^{i_1}\lambda_2^{i_2}\lambda_3^{i_3}$ is
\begin{align*}
-2\sigma_{3}\lambda_1^{i_1}\lambda_2^{i_2}\lambda_3^{i_3}
&\Big(
-\fac(\lambda_{12})\fac(\lambda_{13})
+\fac(\lambda_{21})\fac(\lambda_{31})
-\fac(\lambda_{21})\fac(\lambda_{23})\\&
+\fac(\lambda_{12})\fac(\lambda_{32})
-\fac(\lambda_{32})\fac(\lambda_{31})
+\fac(\lambda_{23})\fac(\lambda_{13})
\Big)=0
\end{align*}
By symmetry, all other terms involving 
$\lambda_1^{i_a}\lambda_2^{i_b}\lambda_3^{i_c}$, for $\{a, b, c\}=\{1, 2, 3\}$, 
in  \eqref{eq:serre check} are zero. 
Therefore, 
$\Psi(  \Sym_{\mathfrak{S}_3}[e_{i_1}, [e_{i_2}, e_{i_3+1}]])=0$.  
This completes the proof. 
$\blacksquare$

\section{The proof of Proposition~\ref{prop:coprod}}
\label{App2}
In this section, we prove Proposition~\ref{prop:coprod}. We follow the notations from \eqref{eq:corresp stab}.
\begin{defn}\label{mapl}
Passing to the localization of $V_{\vec{r}}$ and $ V_{\vec{r'}}\otimes V_{\vec{r''}}$, we define the map 
$l: V_{\vec{r}}\to V_{\vec{r'}}\otimes V_{\vec{r''}}$ in Proposition \ref{prop:coprod} as the inverse of the following morphism:
\[
(-1)^{(r_1'+r_2'+r_3')n''} \eta^{st}_{*} \circ (p^{st})^*: V_{\vec{r'}}(n')\otimes V_{\vec{r''}}(n'')\to V_{\vec{r}}(n). 
\]
\end{defn}

\begin{lmm}\label{coprod:Y0}
Notations as in \S~\ref{subsec:3dCOHA}, let $\calV_n$ and $\calE_{\vec{r}}$ be the tautological bundles on $\calM_{\vec{r}}(n)$. 
We have
\begin{enumerate}
\item $p^*(\calV_{n'}\boxtimes\calV_{n''})=\eta^*\calV_{n}$; 
\item $p^*(\calE_{\vec{r'}}\boxtimes\calE_{\vec{r''}})=\eta^*\calE_{\vec{r''}}$;
\item Consequently,  for $\psi(z)\in \calH^0$, and for all $x\in V_{r_1, r_2, r_3}$, we have
\[
\Delta^{\Dr}(\psi(z))\bullet l(x)=l(\psi(z) \bullet x)
\] 
\end{enumerate}
\end{lmm}
{\it Proof.}
The first two are clear by definition. Now we prove the last equation.
By the definition of $l$, it suffices to show that 
\[
\eta^{st}_{*} \circ (p^{st})^* (\Delta^{\Dr}(\psi(z))\bullet \tilde{l} x)=\psi(z)\bullet x, \]
where $\tilde{l}:=(-1)^{(r_1'+r_2'+r_3')n''} l$.  
By definition $\Delta^{\Dr}(\psi(z))=\psi(z)\otimes \psi(z)$ and $\psi(z) \bullet x=\lambda_{-1/z}(\calF_{n,\vec{r}})\cdot x$. 
Therefore, 
\begin{align*}
\eta^{st}_{*} \circ (p^{st})^* (\Delta^{\Dr}(\psi(z))\bullet \tilde{l} x)
=&\eta^{st}_{*} \circ (p^{st})^* (\psi(z)\otimes \psi(z)\bullet \tilde{l} x)\\
=&\eta^{st}_{*} \circ (p^{st})^* \Big(\lambda_{-1/z}(\calF_{n',\vec{r}'} \boxtimes \calF_{n',\vec{r}'})\cdot  \tilde{l} x\Big)\\
=&\eta^{st}_{*} \Big(  \lambda_{-1/z}((\eta^{st})^* \calF_{n,\vec{r}}) \cdot  (p^{st})^*(\tilde{l} x)\Big)\\
=&\lambda_{-1/z}(\calF_{n,\vec{r}}) \cdot \eta^{st}_{*}  (p^{st})^*(\tilde{l} x) \\
=& \psi(z)\bullet \tilde{l}^{-1} \tilde{l}(x)
= \psi(z)\bullet x. \end{align*}
This completes the proof. 
$\blacksquare$

We now prove Proposition~\ref{prop:coprod}. 
Lemma \ref{coprod:Y0} implies Proposition~\ref{prop:coprod} when $\alpha$ is an element in $\calH^0$. 
We will now focus on the case when $\alpha\in \mathcal{SH}^{(Q_3, W_3)}$. 
The proof below is similar to the proof of associativity of the Hall multiplication.

%Similar to the correspondence $Z_{(\vec{r}', \vec{r}'' )}(n', n'')$, we also have $Z_{(\vec{r}', \vec{r}'')}(n'+1, n'')$ and $  Z_{(\vec{r}', \vec{r}'')}(n', n''+1)$.
The diagram \eqref{eq:asso} and correspondence \eqref{eq:corresp stab}  induce the following diagram
\[
\xymatrix@C=1em{
%------- line 1--------
{\begin{matrix}
(\mathfrak{gl}_1)^3 \times \calM_{ \vec{r}'}^{st}(n')\\
 \times \calM_{\vec{r}''}^{st}(n'') 
\end{matrix}}
&
{\begin{matrix}
(Z_{(\vec{0}, \vec{r}')}(1, n')^{st} \times \calM_{\vec{r}''}^{st}(n''))  \\ 
\sqcup  (\calM_{\vec{r}'}^{st}(n')\times Z_{(\vec{0}, \vec{r}'')}(1, n'')^{st})  
\end{matrix}}\ar[l]\ar[r]&  
{\begin{matrix}
(\calM_{\vec{r}'}^{st}(n'+1)  \times \calM_{\vec{r}''}^{st}(n'')) \\
\sqcup (\calM_{\vec{r}'}^{st}(n')  \times \calM_{\vec{r}''}^{st}(n''+1))
\end{matrix}}
\\
%------- line 2--------
{\begin{matrix}
(\mathfrak{gl}_1)^3 \times \\
Z_{(\vec{r}', \vec{r}'')}(n',n'')^{st}
\end{matrix}}
\ar[u]\ar[d]
&
{
\begin{matrix}
F_{(\vec{0}, \vec{r}', \vec{r}'')}(1, n', n'')^{st}\\
\sqcup F_{(\vec{r}', \vec{0},\vec{r}'')}(n',1, n'')^{st}
\end{matrix}
} \ar[u]\ar[d]\ar[l]\ar[r]&
{\begin{matrix}
Z_{(\vec{r}', \vec{r}'')}(n'+1, n'')^{st}  \\ 
\sqcup  Z_{(\vec{r}', \vec{r}'')}(n', n''+1)^{st}  
\end{matrix}}
\ar[u]\ar[d]\\
%------- line 3--------
(\mathfrak{gl}_1)^3 \times
\calM_{\vec{r}}^{st}(n) 
& Z_{(\vec{0}, \vec{r})}(1, n)^{st}\ar[l]\ar[r]& \calM_{\vec{r}}^{st}(n+1) }
\]
The disjoint union in the above breaks the diagram into two diagrams, one for each component of the disjiont union. The maps in the two diagrams are schematically represented as follows. 
Here the numbers $1, n', n''$ in the rectangle are the sizes of the corresponding matrices, which also schematically represents the subquotients in \eqref{eq:asso}.
\begin{equation}\label{diag1}
\xymatrix{
%------- line 1--------
\begin{tikzpicture}
 \draw (0,0.4) rectangle  (0.4,0.8) node[pos=.5] {1} ;
 \draw (0.4,0) rectangle  (0.8, 0.4) node[pos=.5] {$n'$};
\draw(0.8, 0) rectangle (1.2, -0.4) node[pos=.5] {$n''$}; 
\end{tikzpicture}
&
\begin{tikzpicture}
  \draw (0,0.4) rectangle  (0.4,0.8) node[pos=.5] {1} ;
 \draw (0.4,0) rectangle  (0.8, 0.4) node[pos=.5] {$n'$};
\draw(0.8, 0) rectangle (1.2, -0.4) node[pos=.5] {$n''$}; 
\draw (0.4,0.4) rectangle (0.8, 0.8) node[pos=.5] {$*$}; 
\end{tikzpicture}
\ar[l]\ar[r]&  
\begin{tikzpicture}
  \draw (0,0.4) rectangle  (0.4,0.8) node[pos=.5] {1} ;
 \draw (0.4,0) rectangle  (0.8, 0.4) node[pos=.5] {$n'$};
\draw(0.8, 0) rectangle (1.2, -0.4) node[pos=.5] {$n''$}; 
\draw (0.4,0.4) rectangle (0.8, 0.8) node[pos=.5] {$*$}; 
\draw (0,0) rectangle (0.4, 0.4) node[pos=.5] {$*$}; 
\end{tikzpicture}
\\
%------- line 2--------
\begin{tikzpicture}
 \draw (0,0.4) rectangle  (0.4,0.8) node[pos=.5] {1} ;
 \draw (0.4,0) rectangle  (0.8, 0.4) node[pos=.5] {$n'$};
\draw(0.8, 0) rectangle (1.2, -0.4) node[pos=.5] {$n''$}; 
\draw (0.8,0) rectangle  (1.2, 0.4) node[pos=.5] {$*$};
\end{tikzpicture}
\ar[u]\ar[d]
&
\begin{tikzpicture}
  \draw (0,0.4) rectangle  (0.4,0.8) node[pos=.5] {1} ;
 \draw (0.4,0) rectangle  (0.8, 0.4) node[pos=.5] {$n'$};
\draw(0.8, 0) rectangle (1.2, -0.4) node[pos=.5] {$n''$}; 
\draw (0.4,0.4) rectangle (0.8, 0.8) node[pos=.5] {$*$}; 
\draw (0.8,0) rectangle  (1.2, 0.4) node[pos=.5] {$*$};
\draw (0.8,0.4) rectangle  (1.2, 0.8) node[pos=.5] {$*$};
\end{tikzpicture}
\ar[l]\ar[r]
\ar[u]\ar[d]&  
\begin{tikzpicture}
  \draw (0,0.4) rectangle  (0.4,0.8) node[pos=.5] {1} ;
 \draw (0.4,0) rectangle  (0.8, 0.4) node[pos=.5] {$n'$};
\draw(0.8, 0) rectangle (1.2, -0.4) node[pos=.5] {$n''$}; 
\draw (0.4,0.4) rectangle (0.8, 0.8) node[pos=.5] {$*$}; 
\draw (0,0) rectangle (0.4, 0.4) node[pos=.5] {$*$}; 
\draw (0.8,0) rectangle  (1.2, 0.4) node[pos=.5] {$*$};
\draw (0.8,0.4) rectangle  (1.2, 0.8) node[pos=.5] {$*$};
\end{tikzpicture}\ar[u]\ar[d]
\\
%------- line 3--------
\begin{tikzpicture}
 \draw (0,0.4) rectangle  (0.4,0.8) node[pos=.5] {1} ;
 \draw (0.4,0) rectangle  (0.8, 0.4) node[pos=.5] {$n'$};
\draw(0.8, 0) rectangle (1.2, -0.4) node[pos=.5] {$n''$}; 
\draw (0.8,0) rectangle  (1.2, 0.4) node[pos=.5] {$*$};
\draw (0.4, -0.4) rectangle  (0.8, 0) node[pos=.5] {$*$};
\end{tikzpicture}
&
\begin{tikzpicture}
  \draw (0,0.4) rectangle  (0.4,0.8) node[pos=.5] {1} ;
 \draw (0.4,0) rectangle  (0.8, 0.4) node[pos=.5] {$n'$};
\draw(0.8, 0) rectangle (1.2, -0.4) node[pos=.5] {$n''$}; 
\draw (0.4,0.4) rectangle (0.8, 0.8) node[pos=.5] {$*$}; 
\draw (0.8,0) rectangle  (1.2, 0.4) node[pos=.5] {$*$};
\draw (0.8,0.4) rectangle  (1.2, 0.8) node[pos=.5] {$*$};
\draw (0.4, -0.4) rectangle  (0.8, 0) node[pos=.5] {$*$};
\end{tikzpicture}
\ar[l]\ar[r]&  
{\begin{tikzpicture}
  \draw (0,0.4) rectangle  (0.4,0.8) node[pos=.5] {1} ;
 \draw (0.4,0) rectangle  (0.8, 0.4) node[pos=.5] {$n'$};
\draw(0.8, 0) rectangle (1.2, -0.4) node[pos=.5] {$n''$}; 
\draw (0.4,0.4) rectangle (0.8, 0.8) node[pos=.5] {$*$}; 
\draw (0,0) rectangle (0.4, 0.4) node[pos=.5] {$*$}; 
\draw (0.8,0) rectangle  (1.2, 0.4) node[pos=.5] {$*$};
\draw (0.8,0.4) rectangle  (1.2, 0.8) node[pos=.5] {$*$};
\draw (0.4, -0.4) rectangle  (0.8, 0) node[pos=.5] {$*$};
\draw (0, -0.4) rectangle  (0.4, 0) node[pos=.5] {$*$};
\end{tikzpicture}}
}
\end{equation}
and
\begin{equation}\label{diag2}
\xymatrix@C=1em{
%------- line 1--------
%====1.1
\begin{tikzpicture}
 \draw (0,0.4) rectangle  (0.4,0.8) node[pos=.5] {1} ;
 \draw (0.4,0) rectangle  (0.8, 0.4) node[pos=.5] {$n'$};
\draw(0.8, 0) rectangle (1.2, -0.4) node[pos=.5] {$n''$}; 
\end{tikzpicture} \cong
\begin{tikzpicture}
 \draw (0,0.4) rectangle  (0.4,0.8) node[pos=.5] {$n'$} ;
 \draw (0.4,0) rectangle  (0.8, 0.4) node[pos=.5] {1};
\draw(0.8, 0) rectangle (1.2, -0.4) node[pos=.5] {$n''$}; 
\end{tikzpicture}
%====1.2
&
\begin{tikzpicture}
 \draw (0,0.4) rectangle  (0.4,0.8) node[pos=.5] {$n'$} ;
 \draw (0.4,0) rectangle  (0.8, 0.4) node[pos=.5] {$1$};
\draw(0.8, 0) rectangle (1.2, -0.4) node[pos=.5] {$n''$}; 
\draw (0.8,0) rectangle  (1.2, 0.4) node[pos=.5] {$*$};
\end{tikzpicture}
\ar[l]\ar[r]&  
%====1.3
\begin{tikzpicture}
 \draw (0,0.4) rectangle  (0.4,0.8) node[pos=.5] {$n'$} ;
 \draw (0.4,0) rectangle  (0.8, 0.4) node[pos=.5] {$1$};
\draw(0.8, 0) rectangle (1.2, -0.4) node[pos=.5] {$n''$}; 
\draw (0.8,0) rectangle  (1.2, 0.4) node[pos=.5] {$*$};
\draw (0.4, -0.4) rectangle  (0.8, 0) node[pos=.5] {$*$};
\end{tikzpicture}
\\
%------- line 2--------
%===2.1
\begin{tikzpicture}
 \draw (0,0.4) rectangle  (0.4,0.8) node[pos=.5] {1} ;
 \draw (0.4,0) rectangle  (0.8, 0.4) node[pos=.5] {$n'$};
\draw(0.8, 0) rectangle (1.2, -0.4) node[pos=.5] {$n''$}; 
\draw (0.8,0) rectangle  (1.2, 0.4) node[pos=.5] {$*$};
\end{tikzpicture}
\ar@<3.5ex>[u]\ar[d]
%===2.2
&
\begin{tikzpicture}
  \draw (0,0.4) rectangle  (0.4,0.8) node[pos=.5] {$1$} ;
 \draw (0.4,0) rectangle  (0.8, 0.4) node[pos=.5] {$n'$};
\draw(0.8, 0) rectangle (1.2, -0.4) node[pos=.5] {$n''$}; 
\draw (0.4,0.4) rectangle (0.8, 0.8) node[pos=.5] {$*$}; 
\draw (0.8,0) rectangle  (1.2, 0.4) node[pos=.5] {$*$};
\draw (0.8,0.4) rectangle  (1.2, 0.8) node[pos=.5] {$*$};
\end{tikzpicture}
\,\ \,\ \,\ \,\ 
\begin{tikzpicture}
  \draw (0,0.4) rectangle  (0.4,0.8) node[pos=.5] {$n'$} ;
 \draw (0.4,0) rectangle  (0.8, 0.4) node[pos=.5] {$1$};
\draw(0.8, 0) rectangle (1.2, -0.4) node[pos=.5] {$n''$}; 
\draw (0.4,0.4) rectangle (0.8, 0.8) node[pos=.5] {$*$}; 
\draw (0.8,0) rectangle  (1.2, 0.4) node[pos=.5] {$*$};
\draw (0.8,0.4) rectangle  (1.2, 0.8) node[pos=.5] {$*$};
\end{tikzpicture}
\ar[l] \ar[r]
\ar@<-3.5ex>[u]\ar@<-3.5ex>[d]
%===2.3
&  
\begin{tikzpicture}
  \draw (0,0.4) rectangle  (0.4,0.8) node[pos=.5] {$n'$} ;
 \draw (0.4,0) rectangle  (0.8, 0.4) node[pos=.5] {$1$};
\draw(0.8, 0) rectangle (1.2, -0.4) node[pos=.5] {$n''$}; 
\draw (0.4,0.4) rectangle (0.8, 0.8) node[pos=.5] {$*$}; 
\draw(0.4, -0.4) rectangle (0.8, 0) node[pos=.5] {$*$}; 
\draw (0.8,0) rectangle  (1.2, 0.4) node[pos=.5] {$*$};
\draw (0.8,0.4) rectangle  (1.2, 0.8) node[pos=.5] {$*$};
\end{tikzpicture}\ar[u]\ar@<3.5ex>[d]
\\
%------- line 3--------
%====3.1
\begin{tikzpicture}
 \draw (0,0.4) rectangle  (0.4,0.8) node[pos=.5] {1} ;
 \draw (0.4,0) rectangle  (0.8, 0.4) node[pos=.5] {$n'$};
\draw(0.8, 0) rectangle (1.2, -0.4) node[pos=.5] {$n''$}; 
\draw (0.8,0) rectangle  (1.2, 0.4) node[pos=.5] {$*$};
\draw (0.4, -0.4) rectangle  (0.8, 0) node[pos=.5] {$*$};
\end{tikzpicture}
&
%====3.2
\begin{tikzpicture}
  \draw (0,0.4) rectangle  (0.4,0.8) node[pos=.5] {1} ;
 \draw (0.4,0) rectangle  (0.8, 0.4) node[pos=.5] {$n'$};
\draw(0.8, 0) rectangle (1.2, -0.4) node[pos=.5] {$n''$}; 
\draw (0.4,0.4) rectangle (0.8, 0.8) node[pos=.5] {$*$}; 
\draw (0.8,0) rectangle  (1.2, 0.4) node[pos=.5] {$*$};
\draw (0.8,0.4) rectangle  (1.2, 0.8) node[pos=.5] {$*$};
\draw (0.4, -0.4) rectangle  (0.8, 0) node[pos=.5] {$*$};
\end{tikzpicture}
\ar[l]\ar[r]&  
%====3.3
{\begin{tikzpicture}
  \draw (0,0.4) rectangle  (0.4,0.8) node[pos=.5] {1} ;
 \draw (0.4,0) rectangle  (0.8, 0.4) node[pos=.5] {$n'$};
\draw(0.8, 0) rectangle (1.2, -0.4) node[pos=.5] {$n''$}; 
\draw (0.4,0.4) rectangle (0.8, 0.8) node[pos=.5] {$*$}; 
\draw (0,0) rectangle (0.4, 0.4) node[pos=.5] {$*$}; 
\draw (0.8,0) rectangle  (1.2, 0.4) node[pos=.5] {$*$};
\draw (0.8,0.4) rectangle  (1.2, 0.8) node[pos=.5] {$*$};
\draw (0.4, -0.4) rectangle  (0.8, 0) node[pos=.5] {$*$};
\draw (0, -0.4) rectangle  (0.4, 0) node[pos=.5] {$*$};
\end{tikzpicture}}
\cong
{\begin{tikzpicture}
  \draw (0,0.4) rectangle  (0.4,0.8) node[pos=.5] {$n'$} ;
 \draw (0.4,0) rectangle  (0.8, 0.4) node[pos=.5] {$1$};
\draw(0.8, 0) rectangle (1.2, -0.4) node[pos=.5] {$n''$}; 
\draw (0.4,0.4) rectangle (0.8, 0.8) node[pos=.5] {$*$}; 
\draw (0,0) rectangle (0.4, 0.4) node[pos=.5] {$*$}; 
\draw (0.8,0) rectangle  (1.2, 0.4) node[pos=.5] {$*$};
\draw (0.8,0.4) rectangle  (1.2, 0.8) node[pos=.5] {$*$};
\draw (0.4, -0.4) rectangle  (0.8, 0) node[pos=.5] {$*$};
\draw (0, -0.4) rectangle  (0.4, 0) node[pos=.5] {$*$};
\end{tikzpicture}}
}
\end{equation}

We have two ways from the lower left corner $(\mathfrak{gl}_1)^3 \times \calM_{ \vec{r}}^{st}(n)$ 
to the upper right corner $
(\calM_{\vec{r}'}^{st}(n'+1)  \times \calM_{\vec{r}''}^{st}(n'')) 
\sqcup (\calM_{\vec{r}'}^{st}(n')  \times \calM_{\vec{r}''}^{st}(n''+1))$. The corresponding maps
\begin{align*}
&H^*_c((\mathfrak{gl}_1)^3 \times \calM_{ \vec{r}}^{st}(n);\varphi_{\tr})^\vee\to 
\big(V_{\vec{r}'}(n'+1)\otimes V_{\vec{r}''}(n'')\big)\bigoplus \big(V_{\vec{r}'}(n')\otimes V_{\vec{r}''}(n''+1)\big)
%\\&H^*_c((\calM_{\vec{r}'}^{st}(n'+1)  \times \calM_{\vec{r}''}^{st}(n'')) ;\varphi_{\tr})^\vee\oplus H^*_c((\calM_{\vec{r}'}^{st}(n')  \times \calM_{\vec{r}''}^{st}(n''+1));\varphi_{\tr})^\vee
\end{align*}
will simply be referred to as Way 1 and Way 2.
 
Way 1: using the bottom horizontal and right vertical correspondences of \eqref{diag1} \eqref{diag2}, and follow the standard procedure as in \S~\ref{sec:CritCOHA_act}.

Way 2: using the left vertical and top horizontal correspondences of \eqref{diag1} \eqref{diag2}, and follow the standard procedure as in \S~\ref{sec:CritCOHA_act}.

For any $\alpha\in \mathcal{SH}^{(Q_3, W_3)}$, $x \in V_{\vec{r}}$, clearly Way 1 applied to $(\alpha\otimes x)$ gives $l (\alpha \bullet x)$. 
We claim that Way 2  applied to $(\alpha\otimes x)$ gives $\Delta^{\Dr}(\alpha) \bullet l(x)$. 
Then, convolutions on the level of critical cohomology of  diagrams \eqref{diag1},  \eqref{diag2} lead to $l (\alpha \bullet x)=\Delta^{\Dr}(\alpha) \bullet l(x)$. This in turn  implies Proposition~\ref{prop:coprod}(1). As the maps $l$ and $\Delta^{\Dr}$ are co-associative in the natural sense, applying Proposition~\ref{prop:coprod}(1) iteratively gives Proposition~\ref{prop:coprod}(2).

In the rest of this section we present the proof in relative details. 

The extension 
$\begin{tikzpicture}
\draw (0,0) rectangle  (0.4,0.4) node[pos=.5] {1} ;
\draw (0.4,-0.4) rectangle  (0.8, 0) node[pos=.5] {$n$};
\draw [pattern=north west lines] (0.4,0) rectangle  (0.8, 0.4);
\end{tikzpicture}$ in \eqref{diag3} is 
\[
Z_{(0, \vec{r})}(1, n)^{st}=\{(B_i, I_{ab}, J_{ab})_{i\in \underline{3}, ab\in \overline{3}}\in \calM_{ \vec{r}}(1+n)^{st}\mid 
B_i(\calV_1) \subset \calV_1, 
J_{ab}(\calV_1 )=0
\}.\] The extension $\begin{tikzpicture}
\draw (0,0) rectangle  (0.4,0.4) node[pos=.5] {$n$} ;
\draw (0.4,-0.4) rectangle  (0.8, 0) node[pos=.5] {$1$};
\draw [pattern=dots] (0.4,0) rectangle  (0.8, 0.4);
\end{tikzpicture}$ consists of
\[
Z_{(\vec{r}, 0)}(n,1)^{op}:=\{(B_i, I_{ab}, J_{ab})_{i\in \underline{3}, ab\in \overline{3}}\in \calM_{ \vec{r}}(1+n)\mid 
B_i(\calV) \subset \calV,  I_{ab}(\calE_{r_c})\subset \calV, a\neq b\neq c\}, 
\] 
where $op$ means the opposite stability condition. 
%\[
%Z_{(\vec{r}, 0)}(n,1)^{op}:=\{(B_i, I_{ab}, J_{ab})_{i\in \underline{3}, ab\in \overline{3}}\in \calM_{ \vec{r}}(1+n)\mid 
%B_i(\calV) \subset \calV,  J_{ab}(\calV_1 )=0\}.
%\] 

We consider \[
\xymatrix@R=0.5em{
(\mathfrak{gl}_1)^3 \times \calM_{ \vec{r}}^{st}(n) & Z_{(0, \vec{r})}(1, n)^{st} \ar[l]_{p_1}\ar[r]^{\eta_1}& \calM_{\vec{r}}^{st}(1+n)\\
\calM_{ \vec{r}}^{st}(n)\times (\mathfrak{gl}_1)^3 & Z_{(\vec{r}, 0)}(n,1)^{op}\ar[l]_{p_2}\ar[r]^{\eta_2}& \calM_{\vec{r}}^{st}(1+n)
}
\]
\begin{equation}\label{diag3}
\xymatrix{
	%------- line 1--------
	\begin{tikzpicture}
	\draw (0,0) rectangle  (0.4,0.4) node[pos=.5] {1} ;
	\draw (0.4,-0.4) rectangle  (0.8, 0) node[pos=.5] {$n$};
	\end{tikzpicture} \ar[d]^{\cong}
	&
	\begin{tikzpicture}
	\draw (0,0) rectangle  (0.4,0.4) node[pos=.5] {1} ;
	\draw (0.4,-0.4) rectangle  (0.8, 0) node[pos=.5] {$n$};
	\draw [pattern=north west lines] (0.4,0) rectangle  (0.8, 0.4);
	\end{tikzpicture}
	\ar[l]_{p_1}\ar[r]^{\eta_1}&
	\begin{tikzpicture}
	\draw (0,0) rectangle  (0.4,0.4) node[pos=.5] {1} ;
	\draw (0.4,-0.4) rectangle  (0.8, 0) node[pos=.5] {$n$};
	\draw (0.4,0) rectangle  (0.8, 0.4) node[pos=.5] {$*$};
	\draw (0,-0.4) rectangle  (0.4, 0) node[pos=.5] {$*$};
	\end{tikzpicture}\ar[d]^{\cong}
	\\
	%------- line 2--------
	\begin{tikzpicture}
	\draw (0,0) rectangle  (0.4,0.4) node[pos=.5] {$n$} ;
	\draw (0.4,-0.4) rectangle  (0.8, 0) node[pos=.5] {$1$};
	\end{tikzpicture} 
	&
	\begin{tikzpicture}
	\draw (0,0) rectangle  (0.4,0.4) node[pos=.5] {$n$} ;
	\draw (0.4,-0.4) rectangle  (0.8, 0) node[pos=.5] {$1$};
	\draw [pattern=dots] (0.4,0) rectangle  (0.8, 0.4);
	\end{tikzpicture}
	\ar[l]_{p_2}\ar[r]^{\eta_2}&
	\begin{tikzpicture}
	\draw (0,0) rectangle  (0.4,0.4) node[pos=.5] {$n$} ;
	\draw (0.4,-0.4) rectangle  (0.8, 0) node[pos=.5] {$1$};
	\draw (0.4,0) rectangle  (0.8, 0.4) node[pos=.5] {$*$};
	\draw (0,-0.4) rectangle  (0.4, 0) node[pos=.5] {$*$};
	\end{tikzpicture}
}\end{equation}
\begin{lmm} \label{lem:switch}
For $\alpha(\lambda)\in \H^{Q_3, W_3}(1)=\C[\lambda]$, and $x\in V_{r_1, r_2, r_3}(n)$, 
we have the equality
\[
(\eta_1)_*(p_1)^*(\alpha(\lambda), x)=(-1)^{r_1+r_2+r_3}(\eta_2)_*\big( \psi(\lambda) \cdot (p_2)^*(x, \alpha(\lambda))\big). 
\]
\end{lmm}
{\it Proof.}
Let us compute the difference between the two maps induced by the two diagrams above. Let $\calV$ (resp. $\calE_{r_i}, i=1, 2, 3$) be the tautological dimension $n$ (resp. dimension $r_i, i=1, 2, 3$) bundle on $\calM_{r_1, r_2, r_3}(n)$
and $\calV_1$ is the tautological line bundle on $\calM_{0, 0, 0}(1)$. 

Let $\lambda$ be the Chern root of $\calV_1$, $\lambda_1, \cdots, \lambda_n$ the Chern roots of $\calV$, and 
$\mu_1, \cdots, \mu_{r_i}$ the Chern roots of $\calE_{r_i}$, $i=1, 2, 3$. 
The difference of 
$\begin{tikzpicture}
\draw (0,0) rectangle  (0.4,0.4) node[pos=.5] {$n$} ;
\draw (0.4,-0.4) rectangle  (0.8, 0) node[pos=.5] {$1$};
\draw [pattern=dots] (0.4,0) rectangle  (0.8, 0.4);
\end{tikzpicture}$
and
$\begin{tikzpicture}
\draw (0,0) rectangle  (0.4,0.4) node[pos=.5] {1} ;
\draw (0.4,-0.4) rectangle  (0.8, 0) node[pos=.5] {$n$};
\draw [pattern=north west lines] (0.4,0) rectangle  (0.8, 0.4);
\end{tikzpicture}
$ is 
\begin{align*}
&\Big(\Hom(\calV_1, \calV) \oplus \Hom(\calV_1, \calE_{r_1})
 \oplus \Hom(\calV_1, \calE_{r_2})
  \oplus \Hom(\calV_1, \calE_{r_3})
\Big)\\
&-\Big(\Hom(\calV, \calV_1) \oplus \Hom(\calE_{r_1}, \calV_1)
\oplus \Hom(\calE_{r_2}, \calV_1)\oplus \Hom(\calE_{r_3}, \calV_1)
\Big)\end{align*}
%\begin{align*}
%&\Big(\Hom(\calV_1, \calV)^3 \oplus \Hom(\calV_1, \calE_{r_1})
%\oplus \Hom(\calV_1, \calE_{r_2})
%\oplus \Hom(\calV_1, \calE_{r_3})
%\Big)\\
%&-\Big(\Hom(\calV, \calV_1)^3 \oplus \Hom(\calV_1, \calE_{r_1})
%\oplus \Hom(\calV_1, \calE_{r_2})
%\oplus \Hom(\calV_1, \calE_{r_3})
%\Big)\end{align*}
Taking into account the torus action in \S \ref{subsec:torus action}, we have
\begin{align*}
&eu\Big(
\oplus_{i=1}^3\Hom(\calV_1, q_i\calV) -\oplus_{i=1}^3\Hom(\calV, q_i\calV_1) \Big)\\
=&\frac
{eu(\calV_1^*\otimes q_1\calV)eu(\calV_1^*\otimes q_2\calV)eu(\calV_1^*\otimes q_3\calV)
}
{eu(\calV^*\otimes q_1\calV_1) eu(\calV^*\otimes q_2\calV_1) eu(\calV^*\otimes q_3\calV_1)
}
\\
=&
\prod_{i=1}^n \frac{  \lambda_i -\lambda+\hbar_1}{ \lambda-\lambda_i +\hbar_1}
\prod_{i=1}^n \frac{  \lambda_i -\lambda+\hbar_2}{ \lambda-\lambda_i +\hbar_2}
\prod_{i=1}^n \frac{  \lambda_i -\lambda+\hbar_3} { \lambda-\lambda_i +\hbar_3}\\
=&
(-1)^{3n}\prod_{i=1}^n \frac{  \lambda -\lambda_i-\hbar_1}{ \lambda-\lambda_i +\hbar_1}
\prod_{i=1}^n \frac{  \lambda -\lambda_i-\hbar_2}{ \lambda-\lambda_i +\hbar_2}
\prod_{i=1}^n \frac{  \lambda -\lambda_i-\hbar_3} { \lambda-\lambda_i +\hbar_3}.
\end{align*}
Similarly,
%\begin{align*}
%&eu\Big(\Hom(\calV_1, q_2 q_3 \calE_{r_1}) -\Hom(\calV_1,\calE_{r_1}) \Big) 
%=eu\Big(\calV_1^*\otimes q_2 q_3 \calE_{r_1} -\calE_{r_1}\otimes \calV_1^* \Big) \\
%=&\prod_{i=1}^{r_1}\frac{ \mu_i-\lambda-\hbar_1}{\mu_i-\lambda}
%= \prod_{i=1}^{r_1}\frac{\lambda- \mu_i+\hbar_1}{\lambda-\mu_i},
%\end{align*}
\begin{align*}
&eu\Big(\Hom(\calV_1, q_2 q_3 \calE_{r_1}) -\Hom(\calE_{r_1}, \calV_1) \Big) 
=eu\Big(\calV_1^*\otimes q_2 q_3 \calE_{r_1} -\calE_{r_1}^*\otimes \calV_1 \Big) \\
=&\prod_{i=1}^{r_1}\frac{ \mu_i-\lambda-\hbar_1}{\lambda-\mu_i}
=(-1)^{r_1} \prod_{i=1}^{r_1}\frac{\lambda- \mu_i+\hbar_1}{\lambda-\mu_i}. 
\end{align*}
and
\begin{align*}
&eu\Big(\Hom(\calV_1, q_1 q_3 \calE_{r_2}) -\Hom(\calE_{r_2}, \calV_1) \Big) 
=(-1)^{r_2} \prod_{i=1}^{r_2}\frac{\lambda- \mu_i+\hbar_2}{\lambda-\mu_i},
\\
&eu\Big(\Hom(\calV_1, q_1 q_2 \calE_{r_3}) -\Hom(\calE_{r_3}, \calV_1) \Big) 
=(-1)^{r_3} \prod_{i=1}^{r_3}\frac{\lambda- \mu_i+\hbar_3}{\lambda-\mu_i}. 
\end{align*}
%and
%\begin{align*}
%&eu\Big(\Hom(\calV_1, q_1 q_3 \calE_{r_2}) -\Hom(\calV_1,\calE_{r_2} ) \Big) 
%= \prod_{i=1}^{r_2}\frac{\lambda- \mu_i+\hbar_2}{\lambda-\mu_i},
%\\
%&eu\Big(\Hom(\calV_1, q_1 q_2 \calE_{r_3}) -\Hom(\calV_1,\calE_{r_3}) \Big) 
%= \prod_{i=1}^{r_3}\frac{\lambda- \mu_i+\hbar_3}{\lambda-\mu_i}. 
%\end{align*}
Therefore, 
\begin{align}
&eu\Big(\big(\Hom(\calV_1, \calV) \oplus \Hom(\calV_1, \calE_{r_1})
 \oplus \Hom(\calV_1, \calE_{r_2})
  \oplus \Hom(\calV_1, \calE_{r_3})
\big) \notag\\
&-\big(\Hom(\calV, \calV_1) \oplus \Hom(\calE_{r_1}, \calV_1)
\oplus \Hom(\calE_{r_2}, \calV_1)\oplus \Hom(\calE_{r_3}, \calV_1)
\big)\Big) \notag\\
=& (-1)^{3n+r_1+r_2+r_3}
\Big(\prod_{a=1}^{r_1}\frac{\lambda-\mu_a+\hbar_1}{\lambda-\mu_a}
\prod_{b=1}^{r_2}\frac{\lambda-\mu_b+\hbar_2}{\lambda-\mu_b}
\prod_{b=1}^{r_3}\frac{\lambda-\mu_c+\hbar_3}{\lambda-\mu_c} \notag\\
&
\phantom{123456789}
\cdot\prod_{d=1}^{n}\frac{\lambda-\lambda_d-\hbar_1}{\lambda-\lambda_d+\hbar_1}
\frac{\lambda-\lambda_d-\hbar_2}{\lambda-\lambda_d+\hbar_2}
\frac{\lambda-\lambda_d-\hbar_3}{\lambda-\lambda_d+\hbar_3}\Big). 
\label{eq:eu of trig}
\end{align}
By the formula \eqref{action of Y0} of $\calH^0$-action on $x\in V_{r_1, r_2, r_3}(n)$,  \eqref{eq:eu of trig} coincides with 
$(-1)^{3n+r_1+r_2+r_3}\psi(\lambda)\cdot x $. 
Recall that in the definition of $\eta_{i*}$ as in  \eqref{eq:action} Step~5, a change of group from a parabolic to $\GL_{n+1}$ appears, where the parabolics are opposite for $\eta_1$ and $\eta_2$, which results a sign $(-1)^n$ in the pushforwards, canceling the sign $(-1)^{3n}$ above.
This completes the proof. 
$\blacksquare$

\begin{lmm}
Way 2 gives $\Delta^{\Dr}(\alpha) \bullet l(x)$. 
\end{lmm}
{\it Proof.}
Assume $\alpha=f(\lambda)$ is an element in $\H^{(Q_3, W_3)}(1) =\C[\lambda]$, for some polynomial $f$. 
By definition of the Drinfeld coproduct, we have 
\begin{align*}
\Delta^{\Dr}(\alpha)
=&\psi(\lambda) \otimes \alpha+\alpha \otimes 1\\
=& 1\otimes \alpha-(\hbar_1\hbar_2\hbar_3) \sum_{j\geq 0} \psi_{j} \otimes \lambda^{j+1} \alpha +\alpha \otimes 1
\in Y^{\geq 0}\widehat{\otimes} Y^{\geq 0}. 
\end{align*}

For $\alpha\in \H^{(Q_3, W_3)}(1)$, and $x_1\in V_{\vec{r}'}, x_2\in V_{\vec{r}''}$, 
applying pullback and pushforward via the following correspondence, 
\[
\xymatrix{
%------- line 1--------
{\begin{matrix}
(\mathfrak{gl}_1)^3 \times 
\calM_{ \vec{r}'}^{st}(n')\\
 \times \calM_{\vec{r}''}^{st}(n'') 
 \end{matrix}}
&
Z_{(\vec{0}, \vec{r}')}(1, n')^{st} \times \calM_{\vec{r}''}^{st}(n'')
\ar[l]\ar[r]&  
(\calM_{\vec{r}'}^{st}(n'+1)  \times \calM_{\vec{r}''}^{st}(n'')) 
}\]
we obtain $(\alpha \otimes 1)\bullet (x_1\otimes x_2)=(\alpha \bullet x_1) \otimes  x_2$. 

We now use the following correspondence
\[
\xymatrix{
{\begin{matrix}
(\mathfrak{gl}_1)^3 \times \calM_{ \vec{r}'}^{st}(n')\\
 \times \calM_{\vec{r}''}^{st}(n'') 
\end{matrix}}
& \calM_{\vec{r}'}^{st}(n')\times Z_{(\vec{0}, \vec{r}'')}(1, n'')^{st}
\ar[l]\ar[r] &  
 (\calM_{\vec{r}'}^{st}(n')  \times \calM_{\vec{r}''}^{st}(n''+1))
}
\]
Notice that we need to switch $(\mathfrak{gl}_1)^3$ with $\calM_{ \vec{r}'}^{st}(n')$ in order the apply pullback and pushforward to the above correspondence. 
By Lemma \ref{lem:switch}, we obtain $(-1)^{r_1'+r_2'+r_3'}(\psi(\lambda) \otimes \alpha) \bullet (x_1\otimes x_2)$.  
Therefore, Way 2 gives $\Delta^{\Dr}(\alpha) \bullet l(x)$. 
This completes the proof. 
$\blacksquare$

\vspace{3mm}
%%%%%%%%%%%%%%%%%%%%%%%%%%%%%%%%%%%%%%%%%%%%%%%%%%%

\section{More on the coproduct $\Delta^{\vec{c}}$}
\label{subsec:proof_hyperLoc}

In this section we prove Proposition~\ref{prop:coprod_formula} and Proposition~\ref{prp:coprod_hyperbolic}.
The outline of the proofs is as follows:

1. Prove the primitivity of the Heisenberg operator $B_1$ using the definition of the hyperbolic localization. 

2. Take the opposite $A$-action on $\fM_{\vec{r}}$, and the comparison of these two $A$-actions defines an $R$-matrix. This $R$-matrix defines a subalgebra of $\End(\oplus_{\vec{r}}V_{\vec{r}})$ (a generalization of the Maulik-Okounkov Yangian in this 3d setting). Step 1, together with the geometric interpretations of the central and Cartan elements shows that $\SH^{\vec{c}}$ has a natural map to this subalgebra.

3. When two of the three coordinates of $\vec{r}=(r_1,r_2,r_3)$ are zero, the hyperbolic localization $h:V_{\vec{r}}\to V_{\vec{r'}}\otimes V_{\vec{r''}}$, under the dimensional reduction, agrees with the one from Nakajima and hence the stable envelope of Maulik-Okounkov. 

4. Observe that the map $\Delta^{\vec{c}}:\SH^{\vec{c}}\to (\SH^{\vec{c}})^{\otimes 2}$ is uniquely determined by the fact that, when applied to modules $V_{\vec{r}}$ with two of the three coordinates of $\vec{r}=(r_1,r_2,r_3)$ being zero, it agrees with \cite{MO,NakPCMI,SV2}. This in particular proves both Proposition~\ref{prop:coprod_formula}
 and Proposition~\ref{prp:coprod_hyperbolic}. 

These will be done respectively in the rest of this section. Through out this section, for simplicity we ignore homological degree shiftings. 

\subsection{The primitivity of $B_1$}\label{subsubsec:B_1}
Consider the following commutative diagram
\begin{equation}\label{diag:hyp_Corr}
\xymatrix@R=1em @C=0.5em{
%------- line 1--------
\fM_{\vec{r}}(n)
&
\calA_{\fM_{\vec{r}}(n)}\ar[l]_{\eta_1}\ar[rr]^{p_1}&  &
\fM_{\vec{r}}(n)^A \supseteq \fM_{\vec{r_1}}(n_1)\times \fM_{\vec{r_2}}(n_2)
\\
%-----Line 1.5
&&\bullet\ar[ru]^(0.4){q_3'}\ar[rd]^(0.4){p_2''}&\\
%------- line 2--------
C_{\vec{r}}(n, n+1)
\ar[uu]_{q_1}\ar[dd]^{b_1}
&
\calA_{C} \ar[uu]_{q_2}\ar[dd]^{b_2}\ar[l]_{\eta_2}\ar[rr]^(0.25){p_2} \ar[ru]^{p_2'}&&
C_{\vec{r}}(n, n+1)^A\supseteq 
{\begin{matrix}
(C_{\vec{r_1}}(n_1, n_1+1)\times \fM_{\vec{r_2}}(n_2))\\
\coprod(\fM_{\vec{r_1}}(n_1)\times C_{\vec{r_2}}(n_2, n_2+1))
\end{matrix}}
\ar[uu]_{q_3}\ar[dd]^{b_3} \ar[ld]_(0.7){b_3'}\\
%-----Line 2.5
&&\bullet\ar[rd]^(0.3){p_3''}&\\
%------- line 3--------
\fM_{\vec{r}}(n+1)
&\calA_{\fM_{\vec{r}}(n+1)}\ar[l]_{\eta_3}\ar[rr]^(0.25){p_3}\ar[ru]^{p_3'}
&  &
\fM_{\vec{r}}(n+1)^A
\supseteq 
{\begin{matrix}
(\fM_{\vec{r_1}}(n_1+1)\times \fM_{\vec{r_2}}(n_2))\\
\coprod(\fM_{\vec{r_1}}(n_1)\times \fM_{\vec{r_2}}(n_2+1))
\end{matrix}}
}
\end{equation}
The right two squares in the diagram are not Cartesian. Nevertheless, the maps $p_2$ and $p_3$ factorize with $p_2''$ and $p_3''$ affine bundles, and the squares with $p_1,q_2, p_2',q_3'$ and $p_2,b_2,p_3',b_3'$ both Cartesian. For clarity, we remark that $p_2''$, defined as a map on the union of two varieties,  is defined as the identity map when restricted to the component with $(n_2+1)$ in the index.  Similarly, $p_3''$ is defined as the identity map on the component with $(n_1+1)$ in the index.

\begin{rmk}
It is tempting at this point to generalize $q_1^*:H^*_c(\fM_{\vec{r}}(n),\varphi
)^\vee\to H^*_c(C_{\vec{r}}(n,n+1),\varphi
)^\vee$ and similarly $q_3^*$ to the correspondence where dimensions differ  by 1. However, it is crucial that the potential function on $C_{\vec{r}}(n, n+1)$ is the pullback of the potential function on $\fM_{\vec{r}}(n)$, a fact that does not generalize to higher correspondences.
\end{rmk}

For simplicity, for any variety $X$, the structure map $X\to\pt$ is denoted by $t_X.$
Composing $\bbD t_{\fM_{\vec{r}}(n+1)!}$, $(\id\to \eta_{3!}\eta_3^!)$, and $(\id\to b_{1*}b_1^*)$, we get a commutative diagram \[\xymatrix{
\bbD t_{\fM_{\vec{r}}(n+1)!}\varphi\C_{\fM_{\vec{r}}(n+1)}&\bbD t_{\fM_{\vec{r}}(n+1)!}b_{1*}b_1^*\varphi\C_{\fM_{\vec{r}}(n+1)}\ar[l]\\
\bbD t_{\fM_{\vec{r}}(n+1)!} \eta_{3!}\eta_3^!\varphi\C_{\fM_{\vec{r}}(n+1)}\ar[u]&\bbD t_{\fM_{\vec{r}}(n+1)!} \eta_{3!}\eta_3^!b_{1*}b_1^*\varphi\C_{\fM_{\vec{r}}(n+1)}\ar[u]\ar[l]
}\]
We have $t_{\fM_{\vec{r}}(n+1)!} \eta_{3!}=t_{\fM_{\vec{r}}(n+1)^A!} p_{3*}$. We compare $t_{\fM_{\vec{r}}(n+1)^A!} p_{3*}\eta_3^!(\id\to b_{1*}b_1^*)$ and $t_{\fM_{\vec{r}}(n+1)^A!} (\id\to b_{3*}b_3^*)p_{3*}\eta_3^!$. 
By commutativity and Cartesian properties of \eqref{diag:hyp_Corr}, we have  
\begin{eqnarray*}&b_{3*}b_3^*p_{3*}\eta_3^!=b_{3*}b_3^{'*}p_3^{''*}p_{3*}''p_{3*}^{'}\eta_3^!\to b_{3*}b_3^{'*}p_{3*}'\eta_3^!=b_{3*}p_{2*}b_2^*\eta_3^!\\&=p_{3*}b_{2*}b_2^*\eta_3^!=p_{3*}b_{2*}\eta_2^!b_1^*=p_{3*}\eta_3^!b_{1*}b_1^*,
\end{eqnarray*}
we get a natural transform $(t_{\fM_{\vec{r}}(n+1)^A!} (\id\to b_{3*}b_3^*)p_{3*}\eta_3^!)\to (t_{\fM_{\vec{r}}(n+1)^A!} p_{3*}\eta_3^!(\id\to b_{1*}b_1^*))$, which therefore gives a commutative diagram 
\[\xymatrix{
\bbD t_{\fM_{\vec{r}}(n+1)!} \eta_{3!}\eta_3^!\varphi\C_{\fM_{\vec{r}}(n+1)}&\bbD t_{\fM_{\vec{r}}(n+1)!} \eta_{3!}\eta_3^!b_{1*}b_1^*\varphi\C_{\fM_{\vec{r}}(n+1)}\ar[l]\\
\bbD t_{\fM_{\vec{r}}(n+1)^A!} p_3^*\eta_3^!\varphi\C_{\fM_{\vec{r}}(n+1)}\ar[u]&\bbD t_{\fM_{\vec{r}}(n+1)^A!} b_{3*}b_3^*p_3^*\eta_3^!\varphi\C_{\fM_{\vec{r}}(n+1)}\ar[l] \ar[u].
}\]
This proves the commutativity of hyperbolic restrictions and pushforwards in the diagram \eqref{diag:hyp_Corr}.

Similarly, composing $\bbD t_{\fM_{\vec{r}}(n)!}$, $\id\to \eta_{1!}\eta_1^!$, and $(q_{1!}q_1^!\to\id)$, we get a commutative square
\[\xymatrix{
\bbD t_{\fM_{\vec{r}}(n)!}q_{1!}q_1^!\varphi\C_{\fM_{\vec{r}}(n)}&\bbD t_{\fM_{\vec{r}}(n)!}\varphi\C_{\fM_{\vec{r}}(n)}\ar[l]\\
\bbD t_{\fM_{\vec{r}}(n)!}\eta_{1!}\eta_1^!q_{1!}q_1^!\varphi\C_{\fM_{\vec{r}}(n)}\ar[u]&\bbD t_{\fM_{\vec{r}}(n)!}\eta_{1!}\eta_1^!\varphi\C_{\fM_{\vec{r}}(n)}\ar[l]\ar[u].
}\] We have $t_{\fM_{\vec{r}}(n)!}\eta_{1!}=t_{\fM_{\vec{r}}(n)^A!}p_{1*}$.
Now we compare $p_{1*}\eta_1^!(q_{1!}q_1^!\to\id)$ and $(q_{3!}q_3^!\to\id)p_{1*}\eta_1^!$.
As we have \[q_{3!}q_3^!p_{1*}\eta_1^!\to q_{3!}p_{2*}''p_2^{''!}q_3^!p_{1*}\eta_1^!
=q_{3!}p_{2*}q_2^!\eta_1^!
=p_{1*}q_{2!}q_2^!\eta_1^!
=p_{1*}q_{2!}\eta_2^!q_1^!
=p_{1*}\eta_{1}^!q_{1!}q_1^!,
\]
we get a natural transform $(p_{1*}\eta_1^!(q_{1!}q_1^!\to\id))\to ((q_{3!}q_3^!\to\id)p_{1*}\eta_1^!)$, which therefore gives a commutative diagram 
\[\xymatrix{
\bbD t_{\fM_{\vec{r}}(n)^A!}p_{1*}\eta_1^!q_{1!}q_1^!\varphi\C_{\fM_{\vec{r}}(n)}&\bbD t_{\fM_{\vec{r}}(n)^A!}p_{1*}\eta_1^!\varphi\C_{\fM_{\vec{r}}(n)}\ar[l]\\
\bbD t_{\fM_{\vec{r}}(n)^A!}q_{3!}q_3^!p_{1*}\eta_1^!\varphi\C_{\fM_{\vec{r}}(n)}\ar[u]&\bbD t_{\fM_{\vec{r}}(n)^A!}p_{1*}\eta_1^!\varphi\C_{\fM_{\vec{r}}(n)}\ar[l]\ar[u].
}\]
This proves the commutativity of hyperbolic restrictions and pullbacks in the diagram \eqref{diag:hyp_Corr}.

\subsection{The $R$-matrices}\label{subsec:R_matrix}
Using the correspondence \[\xymatrix{\fM_{\vec{r}}^A&\calA_{\fM_{\vec{r}}}\ar[l]_{p}\ar[r]^{\eta}&\fM_{\vec{r}}},\]
we have constructed the hyperbolic restriction map
\[
h: V_{r_1, r_2, r_3} \to V_{r_1', r_2', r_3'} \otimes V_{r_1'', r_2'', r_3''}, 
\]
where $(r_1, r_2, r_3)=(r_1', r_2', r_3')+(r_1'', r_2'', r_3'')$. 
Note that the map $h$ depends on the torus $A$ action $t \mapsto A(t)$, $t\in \C^*$. 
Consider the opposite action of $A$, that is, the action given by $t \mapsto A(t^{-1})$,  $t\in \C^*$. 
This opposite action gives rise to another hyperbolic restriction map
\[
h^{op}: V_{r_1, r_2, r_3} \to V_{r_1', r_2', r_3'} \otimes V_{r_1'', r_2'', r_3''}, 
\]
which is an isomorphism after localization by the localization theorem in equivariant cohomology.

Similarly to \cite{MO}, we define the R-matrix to be 
\[
h^{op}\circ h^{-1}: V_{r_1', r_2', r_3'} \otimes V_{r_1'', r_2'', r_3''} \to (V_{r_1', r_2', r_3'} \otimes V_{r_1'', r_2'', r_3''}).
\]
The geometrically defined $R$-matrix gives rise to a Yangian $\mathbb{Y}_{(Q_3^{fr}, W_3^{fr})}$ in the 3d case using the RTT=TTR formalism (see \cite[Section 5.2]{MO} for details). 
In particular, $\mathbb{Y}_{(Q_3^{fr}, W_3^{fr})}$ is a subalgebra
\[
\mathbb{Y}_{(Q_3^{fr}, W_3^{fr})}\subset \prod_{i_1, \dots, i_n} \End(F_{i_1}(u_1)\otimes F_{i_2}(u_2)\cdots \otimes F_{i_n}(u_n) ), 
\]
where $F_{i}$'s form a family of representations, with each  $V_{r_1,r_2,r_3}$ for  $(r_1, r_2, r_3)\in \Z_{\geq 0}^{3}$ being so that  $V_{r_1,r_2,r_3}\cong F_{i_1}(u_{i_1})\otimes\cdots\otimes F_{i_m}(u_{i_m})$ for some $i_1,\dots,i_m$ \cite[\S~5.2.14]{MO}. Here the parameters $u_i$'s are  identified with the $T$-equivariant parameters in the definition of $V_{r_1,r_2,r_3}$.

\begin{lmm}
There is an algebra homomorphism $\SH^{\vec{c}}\to \mathbb{Y}_{(Q_3^{fr}, W_3^{fr})}$, such that the following diagram commutes. 
\[
\xymatrix@R=1em @C=1em{
\SH^{\vec{c}}\ar[r]\ar[d]&\End(\oplus_{\vec{r}}V_{\vec{r}})\\
 \mathbb{Y}_{(Q_3^{fr}, W_3^{fr})}\ar@{^{(}->}[ru]
}
\]
\end{lmm}
{\it Proof.}
The algebra  $\SH^{\vec{c}}$ is generated by $B_1, B_{-1}$ and $\SH^{\vec{c}, 0}$. 
To show $\SH^{\vec{c}}$ maps to $\mathbb{Y}_{(Q_3^{fr}, W_3^{fr})}$, 
we only need to show the actions of $B_1, B_{-1}$ and $\SH^{\vec{c}, 0}$ on $V_{\vec{r}}$ come from the Yangian, for any $\vec{r}\in \Z_{\geq 0}^{3}$. 
By Section \ref{sec:CartanDrinfeld}, the action of $\SH^{\vec{c}, 0}$ is given by the Chern classes of the tautological bundles, therefore, it lies in the Yangian. 
By step 1, the element $B_1$ is primitive. Thus, to show that it lies in the image of of the map $\mathbb{Y}_{(Q_3^{fr}, W_3^{fr})}\to \End_{r\in  \Z_{\geq 0}}(V_{\vec{r}})$, it suffices to prove this in the case when  $\vec{r}$ has two components being zero. 
In this 2d case, it is known that $B_1$ is in the image of the map $\mathbb{Y}_{(Q_3^{fr}, W_3^{fr})}\to \End_{r\in  \Z_{\geq 0}}(V_{\vec{r}})$ \cite[\S~5.2]{BFN}.  Similar argument holds for $B_{-1}$. 
$\blacksquare$

Note that in particular, the generators of $\SH^{\vec{c}}$ are closed under the comultiplication on $\mathbb{Y}_{(Q_3^{fr}, W_3^{fr})}$ which in turn is induced by $h$. Denote this comultiplication  by $\Delta^h$. In \S~\ref{subsec:coprod_form_conclud} we will show $\Delta^h$ agrees with $\bold{\Delta}^{\vec{c}}$ from Proposition~\ref{prop:coprod_formula}.

\subsection{Dimensional reduction}
\label{subsec:dim_red}
Now we work with $\vec{r}=(r_1,r_2,r_3)$ such that two of the three coordinates are zero. Without loss of generality, assume $r_1=r\ne0, r_2=r_3=0$. In this section, write 
$X=\{(B_2,B_3,I_{23},J_{23}) \in  \calM_{\vec{r}}(n) \mid \C\langle B_2,B_3\rangle I_{23}(\C^{r_{23}})=\C^n \}/\GL_n$, and hence by the definition of stability conditions we have $\fM_{\vec{r}}(n)=X\times\mathfrak{gl}_n$, endowed with a projection $\pi_X:X\times\mathfrak{gl}_n\to X$. Following the notations as in \cite{D}, here we write
$Z\subset X$ be the locus consisting of orbits of $(B_2,B_3,I_{23},J_{23})$ where $\tr_W(B_k,I_{ij},J_{ij})_{k,i,j=1,2,3}=0$ for all $B_1\in \mathfrak{gl}_n$, with the natural embedding $i_Z:Z\to X$. Note that  $Z$ is isomorphic to the Nakajima quiver variety of dimension $n$ and framing $r$. 

Let $A\subseteq \GL_r$ be such that  $\fM_{\vec{r}}(n)^A=\coprod_{n'+n''=n}\fM_{\vec{r'}}(n')\times \fM_{\vec{r''}}(n'')$. We focus on one component labeled by $n'+n''=n$, and denote the fixed points loci and the attracting loci by $(X\times\mathfrak{gl}_n)^A\cong X^A\times \mathfrak{l}$, $X^A$, $Z^A$, $\calA_{X\times\mathfrak{gl}_n}\cong \calA_X\times \mathfrak{p}$, $\calA_X$, $\calA_Z$ respectively. We have the following diagram.
\begin{equation*}
\xymatrix{
&&&X^A\times\mathfrak{p}\ar[dl]_{d}\ar[dddl]_{\pi'_{X^A}}\\
X\times\mathfrak{gl}_n\ar[dd]_{\pi_X}&\calA_X\times\mathfrak{p}\ar[urr]^{c}\ar[dd]_{\pi_{\calA_X}}\ar[l]_{j_{X\times\mathfrak{gl}_n}}\ar[r]_{q_{X\times\mathfrak{gl}_n}}&X^A\times\mathfrak{l}\ar[dd]_{\pi_{X^A}}\\
&&&Y\ar[dl]^{d_Y}\\
X&\calA_X\ar[l]_{j_X}\ar[r]_{q_X}\ar[urr]^(0.8){q'}&X^A\\
Z\ar[u]^{i_{Z}}&\calA_Z\ar[u]^{i_{\calA_Z}}\ar[l]_{j_Z}\ar[r]_{q_Z}&Z^A\ar[u]^{i_{Z^A}}\ar[uur]_{i_Y}
}
\end{equation*}
Here we introduce the variety 
\[
Y=\{(c, x, x^*)\mid c\in \mathfrak{p}, x\in \mathfrak{l}, x^*\in \mathfrak{l}\mid c_{\mathfrak{l}}=[x, x^*]\}, 
\] 
where $c_\mathfrak{l}\in\mathfrak{l}$ is the projection of $c\in \mathfrak{p}$ under the natural map $\mathfrak{p}\surj  \mathfrak{l}$. 
Under $\mathfrak{p}\cong \mathfrak{n}\oplus \mathfrak{l}$, we decompose $c$ as $c=c_\mathfrak{l}+c_\mathfrak{n}$, where $c_\mathfrak{n}$ is the corresponding element in $\mathfrak{n}$. 
Note that we have the isomorphism $Y\cong \mathfrak{n}\times\mathfrak{l}^2, (c, x, x^*)\mapsto (c_\mathfrak{n}, x, x^*)$. 
The maps are given by 
$q'_X: \mathcal{A}_X\to Y, (x, x^*)\mapsto ([x, x^*], x_{\mathfrak{l}}, x^*_{\mathfrak{l}})$, 
$d_Y: Y \to X^A, (c, x, x^*)\mapsto (x, x^*)$,  and 
$i_Y: Z^A\to Y, (x, x^*)\mapsto (0, x, x^*)$. 
It is straightforward to check that  the square with $i_Y$, $q'_X$, $q_Z$, and $i_{\calA_Z}$ is Cartesian.
The maps $c,d$ are natural projections. Therefore, the square with $c$, $\pi'_{X^A}$, $q_X$, and $\pi_{\calA_X}$ is Cartesian; The two left squares are both Cartesian. 

The following identities are easy to prove, using the commutativity of the diagrams and the Cartesian properties liste above, as well as Lemma~\ref{lmm:smooth}.
\begin{eqnarray}
\label{eqn1}
q_{Z!}j_Z^*i_Z^*=i_Y^*q'_!j_X^*\\  \label{eqn2}
q_{X\times\mathfrak{gl}_n!}j^*_{X\times\mathfrak{gl}_n}\pi_X^*=d_!d^*\pi^*_{X^A}q_{X!}j_X^*\\ \label{eqn3}
i^*_{Z^A}q_{X!}=i^*_{Z^A}d_{Y!}q'_!=i^{*}_Yd_Y^*d_{Y!}q'_!\\ \label{eqn4}
\pi_{X^A!}q_{X\times\mathfrak{gl}_n!}j^*_{X\times\mathfrak{gl}_n}=q_{X!}j_X^*\pi_{X!}\\ \label{eqn5}
i_{Z^A*}q_{Z!}j_Z^*=q_{X!}j_X^*i_{Z*}
\end{eqnarray}
We have the dimensional reduction $\pi_{X!}\pi_X^!i_{Z*}i_Z^*=\pi_{X!}\varphi_{\tr W}\pi_X^*$ and $\pi_{X^A!}\pi_{X^A}^!i_{Z^A*}i_{Z^A}^*=\pi_{X^A!}\varphi_{\tr W}\pi_{X^A}^*$.

We write $\fg=\mathfrak{gl}_n$ for short. We have the natural isomorphisms of functors 
\begin{align*}
&q_{X!}j_X^*\pi_{X!}\varphi\pi_X^*=\pi_{X^A!}q_{X\times\fg!}j_{X\times\fg}^*\varphi\pi_X^*=\pi_{X^A!}\varphi q_{X\times\mathfrak{gl}_n!}j^*_{X\times\mathfrak{gl}_n} \pi_X^*\\
=&\pi_{X^A!}\varphi d_!d^*\pi^*_{X^A}q_{X!}j_X^*\to \pi_{X^A!}\varphi \pi^*_{X^A}q_{X!}j_X^*,
\end{align*} 
with the equalities given by \eqref{eqn4}, Lemma~\ref{lem:hyper_phi}, and \eqref{eqn2} respectively. Here we apply Lemma~\ref{lem:hyper_phi} to the setting when $A$ is the opposite 1-dimensional subgroup, and then use \cite[Theorem~1]{Bra} to identify the hyperbolic localization functors of Lemma~\ref{lem:hyper_phi} in that setting with the functors $q_{X\times\mathfrak{gl}_n!}j^*_{X\times\mathfrak{gl}_n}$ and  $q_{X\times\fg!}j_{X\times\fg}^*$ used here.
We also have
\begin{align*}
&q_{X!}j_X^*\pi_{X!}\pi_{X}^*i_{Z*}i_Z^*=\pi_{X^A!}q_{X\times\mathfrak{gl}_n!}j^*_{X\times\mathfrak{gl}_n}\pi_{X}^*i_{Z*}i_Z^*=\pi_{X^A!}d_!d^*\pi^*_{X^A}q_{X!}j_X^*i_{Z*}i_Z^*\\
\to& \pi_{X^A!}\pi^*_{X^A}q_{X!}j_X^*i_{Z*}i_Z^*
=\pi_{X^A!}\pi^*_{X^A}i_{Z^A*}q_{Z!}j_Z^*i_Z^*=\pi_{X^A!}\pi^*_{X^A}i_{Z^A*}i_Y^*q'_!j_X^*
\\\to &
 \pi_{X^A!}\pi^*_{X^A}i_{Z^A*}i_{Z^A}^*d_Y^*d_{Y!}q'_!j_X^*=\pi_{X^A!}\pi^*_{X^A}i_{Z^A*}i^*_{Z^A}q_{X!}j_X^*
\end{align*}
with the equalities given by \eqref{eqn4}, \eqref{eqn2}, \eqref{eqn5}, \eqref{eqn1}, \eqref{eqn3} respectively.

\begin{lmm}\label{lem:hyp_red}
The following diagram is commutative. 
\begin{equation}
\xymatrix{
q_{X!}j_X^*\pi_{X!}\varphi\pi_X^*\ar[d]_{\cong}\ar[r]&\pi_{X^A!}\varphi \pi^*_{X^A}q_{X!}j_X^*\ar[d]_{\cong}\\
q_{X!}j_X^*\pi_{X!}\pi_{X}^*i_{Z*}i_Z^*\ar[r]&\pi_{X^A!}\pi^*_{X^A}i_{Z^A*}i^*_{Z^A}q_{X!}j_X^*
}
\end{equation}
with the vertical arrows the dimensional reductions, and horizontal arrows given as above. 
\end{lmm}
{\it Proof.} We follow a similar strategy as in \cite[\S~A1]{D}, and by the same construction as in {\it loc. cit.} we have $\tilde{i}_{Z}$, $\tilde{i}_{\calA_Z}$, and $\tilde{i}_{Z^A}$, with the natural transforms $\pi_{X!}\tilde{i}_{Z*}\tilde{i}_Z^*\varphi_{\tr W}\pi_X^*\to \pi_{X!}\tilde{i}_{Z*}\tilde{i}_Z^*\pi_X^*$, and $\pi_{X^A!}\tilde{i}_{Z^A*}\tilde{i}_{Z^A}^*\varphi_{\tr W}\pi_{X^A}^*\to \pi_{X^A!}\tilde{i}_{Z^A*}\tilde{i}_{Z^A}^*\pi_{X^A}^*$ being isomorphisms.
Starting with $q_{X!}j_X^*$ composed with $\pi_{X!}\tilde{i}_{Z*}\tilde{i}_Z^*\varphi_{\tr W}\pi_X^*\to \pi_{X!}\tilde{i}_{Z*}\tilde{i}_Z^*\pi_X^*$ , using (4) we get \[\xymatrix{
q_{X!}j_X^*\pi_{X!}\tilde{i}_{Z*}\tilde{i}_Z^*\varphi_{\tr W}\pi_X^*\ar[r]\ar@{=}[d]&q_{X!}j_X^*\pi_{X!}\tilde{i}_{Z*}\tilde{i}_Z^*\pi_X^*\ar@{=}[d]\\
\pi_{X^A!}q_{X\times\mathfrak{gl}_n!}j^*_{X\times\mathfrak{gl}_n}\tilde{i}_{Z*}\tilde{i}_Z^*\varphi_{\tr W}\pi_X^*\ar[r]& \pi_{X^A!}q_{X\times\mathfrak{gl}_n!}j^*_{X\times\mathfrak{gl}_n}\tilde{i}_{Z*}\tilde{i}_Z^*\pi_X^*
.}\] which is a commutative square of natural isomorphisms. Then by \cite[Lemma~A.4]{D} and Lemma~\ref{lem:hyper_phi}, we get a natural isomorphism 
\begin{equation}
\pi_{X^A!}\varphi_{\tr W}q_{X\times\mathfrak{gl}_n!}j^*_{X\times\mathfrak{gl}_n}\pi_X^*\to \pi_{X^A!}q_{X\times\mathfrak{gl}_n!}j^*_{X\times\mathfrak{gl}_n}\tilde{i}_{Z*}\tilde{i}_Z^*\varphi_{\tr W}\pi_X^*, 
\end{equation}
which composes to \begin{equation}\label{eqn:proof_red_par1}
\pi_{X^A!}\varphi_{\tr W}q_{X\times\mathfrak{gl}_n!}j^*_{X\times\mathfrak{gl}_n}\pi_X^*\to \pi_{X^A!}q_{X\times\mathfrak{gl}_n!}j^*_{X\times\mathfrak{gl}_n}\tilde{i}_{Z*}\tilde{i}_Z^*\pi_X^.
\end{equation}

Similarly, starting with $\pi_{X^A!}\varphi_{\tr W}\pi_{X^A}^*\to \pi_{X^A!}\tilde{i}_{Z^A*}\tilde{i}_{Z^A}^*\varphi_{\tr W}\pi_{X^A}^*\to \pi_{X^A!}\tilde{i}_{Z^A*}\tilde{i}_{Z^A}^*\pi_{X^A}^*$ composed with $q_{X!}j_X^!$, and use (2) we get a commutative diagram
\begin{equation}\label{eqn:eqn:proof_red_par2}
\xymatrix @C=0.8em{
\pi_{X^A!}\varphi_{\tr W}\pi_{X^A}^*q_{X!}j_X^!\ar[d]&\pi_{X^A!}\varphi_{\tr W}d_!d^*\pi_{X^A}^*q_{X!}j_X^!\ar[l]\ar@{=}[r]\ar[d]&\pi_{X^A!}\varphi_{\tr W}q_{X\times\mathfrak{gl}_n!}j^*_{X\times\mathfrak{gl}_n}\pi_X^*\ar[d]\\ 
\pi_{X^A!}\tilde{i}_{Z^A*}\tilde{i}_{Z^A}^*\pi_{X^A}^*q_{X!}j_X^!&\pi_{X^A!}\tilde{i}_{Z^A*}\tilde{i}_{Z^A}^*d_!d^*\pi_{X^A}^*q_{X!}j_X^!\ar[l]\ar@{=}[r]&
\pi_{X^A!}\tilde{i}_{Z^A*}\tilde{i}_{Z^A}^*q_{X\times\mathfrak{gl}_n!}j^*_{X\times\mathfrak{gl}_n}\pi_X^*
}
\end{equation}
Now using the analogues of (3), (1), and (5), we get a natural map 
\[
\pi_{X^A!}\tilde{i}_{Z^A*}\tilde{i}_{Z^A}^*q_{X\times\mathfrak{gl}_n!}j^*_{X\times\mathfrak{gl}_n}\pi_X^*\to \pi_{X^A!}q_{X\times\mathfrak{gl}_n!}j^*_{X\times\mathfrak{gl}_n}\tilde{i}_{Z*}\tilde{i}_Z^*\pi_X^*.\] Its composition with the right vertical map in \eqref{eqn:eqn:proof_red_par2} gives the same map as \eqref{eqn:proof_red_par1}, due to the naturality of the transforms $\id\to \tilde{i}_{Z*}\tilde{i}_Z^*$ and $\varphi_{\tr_W}\to \id$. 
$\blacksquare$

Composing $\bbD t_{X!}$, $\id\to j_{X!}j_X^!$, and $\pi_{X!}\pi_X^!i_{Z*}i_Z^*=\pi_{X!}\varphi_{\tr W}\pi_X^*$, we get the commutative square
\[\xymatrix{
\bbD t_{X!}\pi_{X!}\pi_X^!i_{Z*}i_Z^*\ar@{=}[r]&\bbD t_{X!}\pi_{X!}\varphi_{\tr W}\pi_X^*\\\
\bbD t_{X!}j_{X!}j_X^!\pi_{X!}\pi_X^!i_{Z*}i_Z^*\ar@{=}[r]\ar[u]&\bbD t_{X!}j_{X!}j_X^!\pi_{X!}\varphi_{\tr W}\pi_X^*\ar[u],
}\]
where the left vertical arrow gives the map of Braverman, Finkelberg, and Nakajima \cite[(5.6.4)]{BFN}, and the right vertical arrow gives $h$.
Applying the natural isomorphism of functors $t_{X!}j_{X!}j_X^!=t_{X^A!}q_{X!}j_X^*$, we get the commutative square
\[\xymatrix{
\bbD t_{X!}j_{X!}j_X^!\pi_{X!}\pi_X^!i_{Z*}i_Z^*\ar@{=}[d]\ar@{=}[r]&\bbD t_{X!}j_{X!}j_X^!\pi_{X!}\varphi_{\tr W}\pi_X^*\ar@{=}[d]\\
\bbD t_{X^A!}q_{X!}j_X^*\pi_{X!}\pi_X^!i_{Z*}i_Z^*\ar@{=}[r]&\bbD t_{X^A!}q_{X!}j_X^*\pi_{X!}\varphi_{\tr W}\pi_X^*.
}\]
Combining this with Lemma~\ref{lem:hyp_red}, this shows that the hyperbolic localization $h$ defined here is compatible with th hyperbolic localization of Braverman, Finkelberg, and Nakajima \cite[(5.6.4)]{BFN} under the dimensional reduction. The later is in turn known \cite[Theorem~1.6.1.(2) and \S~5.6]{BFN} to be equal to the stable envelope of Maulik and Okounkov  \cite{MO,NakPCMI,SV2}. Therefore, in particular, $\Delta^h$ agrees with $\bold{\Delta}^{\vec{c}}$ when specialized to $\SH^{\bold{c}^{k}}$, $k=1, 2, 3$.

\subsection{Concluding the proofs}\label{subsec:coprod_form_conclud}
By the definition of the central extended algebras, for each $k=1,2,3$, there is a specialization map $\SH^{\vec{c}}\to \SH^{\bold{c}^{(k)}}$, with $\SH^{\bold{c}^{(k)}}$ isomorphic to the algebra $\SH^{\bold{c}}$ defined in \cite{SV2}. 
The direct sum of these three specialization maps embeds the algebra $\SH^{\vec{c}}$ into $ \bigoplus_{k=1, 2, 3}\SH^{\bold{c}^{(k)}}$. Moreover, this embedding fits into  the upper part of the following diagram. 
\[
\xymatrix@C=1em@R=1.5em{
%---line 1-----
 \bigoplus_{i\in \bbN, k=1, 2, 3} \C\{c_i^{(k)}\} \ar[r] \ar@{=}[d]& 
\SH^{\vec{c}} \ar[r] \ar[rdd]^(0.3){\Delta^{\vec{c}}}\ar@{_{(}->}[d]& \SH \ar[r] \ar@{_{(}->}[d]&0\\
%---line 2-----
 \bigoplus_{i\in \bbN, k=1, 2, 3} \C\{c_i^{(k)}\} \ar[r] & 
\bigoplus_{k=1, 2, 3}\SH^{\bold{c}^{(k)}} \ar[r] \ar[rdd]_(0.3){ \oplus_{k}\Delta^{\bold{c}^{(k)}}} & \bigoplus_{k=1, 2, 3}\SH \ar[r] &0\\
%---line 3-----
   & & (\SH^{\vec{c}})^{\otimes 2} \ar[r] \ar@{_{(}->}[d]& \SH^{\otimes 2}\ar@{_{(}->}[d]\\
   %---line 4-----
    & & \bigoplus_{k=1, 2, 3}(\SH^{\bold{c}^{(k)}})^{\otimes 2} 
    \ar[r] &\bigoplus_{k=1, 2, 3} (\SH^{\otimes 2})
}
\]
In other words, $\SH^{\vec{c}}$ is characterized as the universal central extension of $\SH$, which specializes into $\SH^{\bold{c}^{(k)}}$ for each $k=1,2,3$.

From the formula of $\Delta^{\vec{c}}$ on the generators of $\SH^{\vec{c}}$ in Proposition~\ref{prop:coprod_formula}
specializes to 
 $\Delta^{\bold{c}^{(k)}}: \SH^{\bold{c}^{(k)}}\to (\SH^{\bold{c}^{(k)}})^{\otimes 2}$ by specializing $\bold{c}^{(k')}\mapsto 0$, for $1\leq k\neq k'\leq 3$. 
These three maps $\Delta^{\bold{c}^{(k)}}$ are well-defined algebra homomorphisms, which agree on the specialization $\SH$. Therefore, by the aforementioned universal property, $\Delta^{\vec{c}}$ given in Proposition~\ref{prop:coprod_formula} is a well-defined algebra homomorphism. By the same universal property,  $\Delta^{\vec{c}}$ is coassociative. This in particular proves Proposition~\ref{prop:coprod_formula}.

To prove Proposition~\ref{prp:coprod_hyperbolic}, it suffices to prove the commutativity of the left part of
 the following diagram of the actions of $\SH^{\vec{c}}$. \[
\xymatrix@R=1em @C=1em{
%----Line 1----
&\End(\oplus_{\vec{r}} V_{\vec{r}}) \ar[rd]\ar[dd]_(0.3){h}& \\
%----Line 2----
\SH^{\vec{c}} \ar@{^{(}->}[rr] \ar@{^{(}->}[ru] \ar[dd]^{\Delta^{\vec{c}}}& 
& \End(W) \ar[dd]_{h_W}\\
%----Line 3----
& \End(\oplus_{\vec{r}} V_{\vec{r}})\otimes \End(\oplus_{\vec{r}} V_{\vec{r}})\ar[rd]&\\
%----Line 4----
\SH^{\vec{c}}\otimes \SH^{\vec{c}}\ar@{^{(}->}[ur] \ar[rr]&& \End(W)\otimes \End(W)
}
\]
where $W: = (\oplus_{r_1} V_{r_1, 0, 0})\bigoplus (\oplus_{r_2} V_{0, r_2, 0})\bigoplus (\oplus_{r_3} V_{0, 0, r_3})$. 
By definition, $\Delta^h$ from \S~\ref{subsec:R_matrix} automatically makes the entire diagram commutative. Therefore we need to show that $\Delta^h$ agrees with $\Delta^{\vec{c}}$.
By \S~\ref{subsec:dim_red}, $\Delta^{\vec{c}}$ makes the lower square (containing $\Delta^{\vec{c}}$ and $h_W$) commutative. 
Thus, $\Delta^{\vec{c}}$ and $\Delta^h$ agree when restricting on $\SH^{\bold{c}^{k}}$, $k=1, 2, 3$. 
By the same argument as above,  $\Delta^{\vec{c}}$ is determined by the three specializations $\Delta^{\bold{c}^{(k)}}, k=1, 2, 3$. 
Therefore, $\Delta^{\vec{c}}$ and $\Delta^h$ agree.

%%%%%%%%%%%%%%%%%%%%%%%%%%%%%%%%%%%%%%%%%%%%%%%%%%%%%
{\bf The authors:}

M.~Rap\v{c}\'{a}k- { Perimeter Institute for Theoretical Physics, 31 Caroline St N, Waterloo, ON N2L 2Y5, Canada},
{mrapcak@perimeterinstitute.ca}

Y.~Soibelman-
{ Department of Mathematics, Kansas State University, Manhattan, KS 66506, USA},
{soibel@math.ksu.edu}

Y.~Yang- {School of Mathematics and Statistics, The University of Melbourne, 813 Swanston Street, Parkville VIC 3010, Australia},
{yaping.yang1@unimelb.edu.au}

 G.~Zhao-
 {School of Mathematics and Statistics, The University of Melbourne, 813 Swanston Street, Parkville VIC 3010, Australia};
{address until July 2018: Institute of Science and Technology Austria,
Am Campus, 1,
Klosterneuburg 3400,
Austria},  
{gufangz@unimelb.edu.au}

%%%%%%%%%%%%%%%%%%%%%%%%%%%%%%%%%%%%%%%%%%%%%%%%%%
\end{document}